\documentclass[11pt]{report} 

\def \dfll {\leaders \hbox to 1em {\hss.\hss}\hfill}
\def\beq {\begin{equation}}
\def\eeq {\end{equation}}
\newtheorem{theorem}{Theorem}
\newtheorem{propos}[theorem]{Proposition}
\newtheorem{cor}[theorem]{Corollary}
\newtheorem{conjecture}[theorem]{Conjecture}

\newtheorem{defi}{Definition}
\newtheorem{qes}{Question}
\newtheorem{lemma}{Lemma}
\def\blm{\begin{lemma}}
\def\elm{\end{lemma}}
\def\bdf{\begin{defi}}
\def\edf{\end{defi}}
\def\btm{\begin{theorem}}
\def\etm{\end{theorem}}
\def\bpp{\begin{propos}}
\def\epp{\end{propos}}
\def\bQ {\begin{qes}}
\def\eQ {\end{qes}}
\def\btm{\begin{theorem}}
\def\etm{\end{theorem}}
\def\ben{\begin{enumerate}}
\def\een{\end{enumerate}}
\def\into{\hookrightarrow}
\def\ep{\ \hfill{\rule {2.5mm}{2.5mm}}\smallskip}

\newcommand{\lbl}[1]{\label{#1}}

\def \isto {\widetilde{\longrightarrow}}

\def \Cob {{\cal C}ob_3}

\def\Sp{{\rm Sp}}
\def\SL{{\rm SL}}

\def\opc{{\scriptstyle =}}

\def\drulefill{$\mathord\opc\mkern-3mu\cleaders
                  \hbox{$\mkern-2mu\mathord\opc\mkern-2mu$}\hfill\mkern-2mu \mathord\opc$}

\newcommand{\dov}[1]{\vbox{\ialign{##\crcr\noalign
             {\kern-3pt\nointerlineskip}\drulefill \crcr\noalign
             {\kern0pt\nointerlineskip}
             $\hfil\displaystyle{#1}\hfil$\crcr}}}

\def\zer#1{\vbox{\ialign{##\crcr\noalign
             {\kern-3pt\nointerlineskip}
             {$\hfil\;\,\scriptscriptstyle\oslash\hfil\!$} \crcr\noalign
             {\kern.5pt\nointerlineskip}
             $\hfil\displaystyle{#1}\hfil$\crcr}}\!}

\def\ste#1{\vbox{\ialign{##\crcr\noalign
             {\kern-3pt\nointerlineskip}
             {$\hfil\;\,\scriptscriptstyle\bowtie\hfil\!$} \crcr\noalign
             {\kern.5pt\nointerlineskip}
             $\hfil\displaystyle{#1}\hfil$\crcr}}\!}

\newcommand{\TO}[2]{\stackrel {\mbox{#1}}{\hbox to #2pt{\rightarrowfill}}}

\def\thrafill{$\mathsurround=0pt \mathord- \mkern-6mu 
\cleaders\hbox{$\mkern-2mu
\mathord- \mkern-2mu$}\hfill \mkern-6mu\mathord\twoheadrightarrow$}
\newcommand {\onto} [1]{\hbox to #1pt{\thrafill}}

\newcommand{\emptystuff}[1]{}

\usepackage{amsfonts,
            amssymb}

\input{epsf.sty}

\newcommand{\Z}{{\mathbb Z}}
\newcommand{\D}{{\mathbb D}}

\newcommand{\R}{{\mathbb R}}
\newcommand{\Cc}{{\mathbb C}}

\newcommand{\Q}{{\mathbb Q}}
\newcommand{\E}{{\mathbb E}}
\newcommand{\Ss}{{\mathbb S}}
\newcommand{\Aa}{{\mathbb A}}
\newcommand{\Pp}{{\mathbb P}}

\newcommand{\kk}{{\mathbb K}}

\newcommand{\posit}[2]
{\raise -1.4ex\hbox{${\textstyle #1}\atop {\stackrel{\uparrow}{#2}}$}}

\newcommand{\head}[1]{
\smallskip 

\begin{center}{\large \sc #1}\end{center}

\nopagebreak}

\def \lz  {\langle}
\def \rz  {\rangle}

\def \th {\theta}
\def \tb {\bar{\theta}}

\def\Hh{\mbox{\bm${\cal H}$\ubm}}

\newcommand{\ext}[1] {\mbox{\raisebox{.4ex}{$\bigwedge^{\!#1}$}}\mkern-1mu}

\def\Ii {\mbox{\raise .4 ex\hbox{$\int$}$\!\! I$}}

\def \bm {\boldmath}
\def \ubm {\unboldmath}
\newcommand{\Qq}
{\mbox{\bm${\cal Q}$\ubm} \,}

\begin{document} 

\vspace*{1cm}

\begin{center}

\section*{Homology TQFT's and  the Alexander-Reidemeister 
Invariant of 3-Manifolds via Hopf Algebras and Skein Theory }

\bigskip

\medskip

{\large Thomas Kerler}\\
\medskip

May 2001
 \vspace*{1.2cm}

\bigskip

\end{center}

{\small \noindent{\bf Abstract :}
We develop an explicit  skein theoretical algorithm to compute 
the Alexander polynomial of a 3-manifold from a
surgery presentation employing the methods used in 
the construction of quantum invariants of 3-manifolds. 
As a prerequisite we establish and prove a rather unexpected 
equivalence between the topological quantum field
theory constructed by Frohman and Nicas using the
homology of $U(1)$-representation varieties
on the one side and the combinatorially constructed Hennings-TQFT based on the 
quasitriangular Hopf algebra ${\cal N}=\Z/2\ltimes \ext *\R^2$ on the
other side. We find that both TQFT's are $\SL(2,\R)$-equivariant functors and, 
as such, are isomorphic.  The
$\SL(2,\R)$-action in the Hennings construction comes from the natural
action on $\cal N$ and in the case of the Frohman-Nicas theory from the
Hard-Lefschetz decomposition of the $U(1)$-moduli spaces given that they 
are naturally
K\"ahler. The  irreducible components of this TQFT, corresponding  to simple 
representations of $\SL(2,\Z)$ and $\Sp(2g,\Z)$, thus yield a large family of 
homological TQFT's by taking sums and products. 
We give several examples of TQFT's and invariants
that appear to fit  into this family, such as Milnor and Reidemeister Torsion,  
Seiberg-Witten theories, Casson type theories for homology circles \'a la Donaldson, 
higher rank gauge theories following Frohman and Nicas, and the $\Z/p\Z$ reductions of
Reshetikhin-Turaev theories over the cyclotomic integers $\Z[\zeta_p]$. We also
conjecture that the Hennings TQFT for quantum-${\mathfrak s}{\mathfrak l}_2$   
is the product of the Reshetikhin-Turaev TQFT and such a homological TQFT. 
 \footnote{
2000 Mathematics Subject Classification: Primary 57R56; Secondary 14D20, 16W30, 17B37, 18D35, 57M27.}

}
\medskip

\begin{center}
{\large \sc Contents}\
\smallskip

\parbox[t]{12cm}{
1. Introduction\dotfill\pageref{S1}

2. Topological Quantum Field Theory\dotfill\pageref{S2}

3. The Frohman-Nicas TQFT for $U(1)$\dotfill\pageref{S3}

4. The Mapping Class Groups and their Actions on Homology\dotfill\pageref{S4}

5. Hennings TQFT's\dotfill\pageref{S5}

6. The Algebra $\cal N$\dotfill\pageref{S6}

7. The Hennings TQFT for $\cal N$\dotfill\pageref{S7}

8. Skein theory for ${\cal V}_{\cal N}$\dotfill\pageref{S8}

9. Equivalence of ${\cal V}_{\cal N}^{(2)}$ and ${\cal V}^{FN}$\dotfill\pageref{S9}

10. Hard-Lefschetz Decomposition\dotfill\pageref{S10}

11. Alexander-Conway Calculus for 3-Manifolds\dotfill\pageref{S11}

12. Lefschetz compatible Hopf Algebra Structures on $H^*(J(\Sigma))$\dotfill\pageref{S12}

13. More Examples of Homological TQFT's and Open Questions\dotfill\pageref{S13}
}
\end{center}

\head{1. Introduction}\lbl{S1}

In recent years much energy has been put into  finding new ways to
 describe and compute  classical invariants of 
3-manifolds using the tools and structures developed in the relatively new area
of quantum topology. In this paper we will establish another such relation between
quantum and classical invariants, which emerged in quite different guises in recent
research in 3-dimensional topology.

The classical invariant of a 3-manifold $M$ we are interested in here is its
Alexander polynomial $\Delta(M)\in\Z[H_1(M)]$. It is closely related and in most cases identical 
to the Reidemeister Milnor Torsion $r(M)$, see \cite{Mil61} and \cite{Tur76}. More recently,
Meng and Taubes \cite{MengTaub} show that this invariant is also equal to the 
Seiberg Witten invariant for 3-manifolds. Turaev  \cite{Tur98} proves a refined version of this
theorem by comparing the behavior of both invariants under surgery.

On the side of the quantum invariants we consider the formalism used for the  
Hennings invariant of 3-manifolds \cite{Hen96}. This invariant is motivated by 
and follows the same principles as the Witten-Reshetikhin-Turaev invariant,
which is  developed in  \cite{Wi89},  \cite{ResTur91} and \cite{TurVi92}, in the sense
that it 
assigns algebraic data to a surgery presentation for $M$. The innovation of the Hennings
approach is that  it starts
directly from a possibly non-semisimple  Hopf algebra ${\cal A}$ rather than its 
semisimple representation theory. This formalism is refined by Kauffman and Radford in
 \cite{KauRad95}. Also Kuperberg \cite{Kup96} gives a construction that  assigns data 
directly from a Hopf algebra to a Heegaard presentation of  $M$.

In this article we discover and explain in detail the relation between the Hennings theory 
for a certain 
8-dimensional Hopf algebra $\cal N$ and the  (reduced) Alexander polynomial 
$\Delta_{\varphi}(M)\in \Z[t,t^{-1}]$ for the cyclic covering given by an epimorphism
$\varphi:\pi_1(M)\to\Z$. As a consequence we have at our disposal 
the entire combinatorial machinery of the Hennings formalism 
 in order to evaluate the Alexander polynomial from surgery
diagrams. Particularly, we are able to develop from this an efficient skein theoretical
algorithm. 
The method of relating these two very differently defined theories is based 
itself on a quite unexpected equivalence of more refined structures. 

More precisely, it turns out that underlying  both invariants is the structure of a 
topological quantum field theory (TQFT). The notion of a TQFT, which can be thought
of as a fiber functor on a category of cobordisms, was first cast into a mathematical 
axiomatic framework by Atiyah \cite{Ati88}. Typically (or by definition)
 all quantum invariants extend to TQFT's on 3-manifolds with boundaries. 
In the case of the semisimple theories generalizing the Witten-Reshetikhin-Turaev
invariant these TQFT's are described in great detail in \cite{Tur94}. In our
context we need the non-semisimple version as it is worked out for the Hennings invariant
in \cite{Ker96} and in full generality in \cite{KerLub00}.

On the side of the classical invariants Frohman and Nicas \cite{FroNic92} 
managed to give an interpretation of the Alexander polynomial of knot complements  
in the setting of TQFT's. In particular, they 
construct a TQFT ${\cal V}^{FN}$, which  assigns to every surface $\Sigma$ as a vector
space the cohomology  ring $H^*(J(\Sigma))$
of the $U(1)$-representation variety $J(\Sigma)=Hom(\pi_1(\Sigma),U(1))$.
The morphisms are constructed in the style of the Casson invariant from the intersection
numbers of representation varieties for a given Heegaard splitting of a cobordism. The
Alexander polynomial is thus given as the Lefschetz trace over  ${\cal V}^{FN}(C_{\Sigma})$,
where $\Sigma$ is an arbitrary Seifert surface and $C_{\Sigma}$ is the 3-dimensional
cobordisms from $\Sigma$ to itself, obtained by cutting away a neighborhood of $\Sigma$.

  The unexpected upshot is that this functor  ${\cal V}^{FN}$ is isomorphic to the  
Hennings TQFT ${\cal V}_{\cal N}$ for the non-semisimple Hopf algebra 
${\cal N}\cong \Z/2\ltimes \ext *\R^2$. The realization of the abelian gauge    
field theory by a specific Hopf algebra 
is not at all obvious since ${\cal V}^{FN}$ and ${\cal V}_{\cal N}$ are
defined in entirely different ways. In fact the isomorphism between these functors
on the vectors spaces mixes up the degrees of exteriors algebras in still 
puzzling ways.  For these reason the proof is rather explicit and computational. 

Nonetheless, it can be seen quite easily that it is not possible to realize
${\cal V}^{FN}$ as a semisimple theory.  
 Particularly,  ${\cal V}^{FN}$ represents  Dehn twists  by 
matrices of the form $1+N$ where $N$ is nilpotent. Furthermore, the invariant vanishes
on $S^1\times S^2$. Yet, in the semisimple
theories from \cite{Tur94} Dehn twists are represented by  
semisimple matrices $D$ with
$D^n=1$ and the invariant on $S^1\times S^2$ is never zero.

Once ${\cal V}^{FN}$ and thus the Alexander polynomial $\Delta_{\varphi}$
are translated into the language of the Hennings formalism for the 
Hopf algebra ${\cal N}$ we are  in the position to develop a skein theory 
for the computation of $\Delta_{\varphi}$. The skein identities reflect 
algebraic relations in  ${\cal N}$. We derive from this a step by step recipe 
for the computation of the Alexander polynomial.

Another intriguing feature of the two TQFT's is that both of them admit natural 
equivariant $\SL(2,\R)$-actions that have very different origins but are,
nevertheless, intertwined by the isomorphism between them. In the case of ${\cal V}^{FN}$ the 
$\SL(2,\R)$-action on  $H^*(J(\Sigma))$
 is given by the Hard Lefschetz decomposition of  the cohomology ring 
that arises from a K\"ahler structure on $J(\Sigma)$. For ${\cal V}_{\cal N}$
this action is derived from an $\SL(2,\R)$-actions on ${\cal N}$ as a Hopf algebra. 
As a consequence   $H^*(J(\Sigma))$
carries a nonstandard ring-structure induced by that of ${\cal N}^{\otimes g}$,
which, as opposed to the standard one,
is compatible with the Hard Lefschetz $\SL(2,\R)$-action. 
\medskip

Let us summarize the content and the
main results of this paper in better order and  detail. In Section~2  
we recall relevant 
notions that characterize 
 topological quantum field theories, such as (non)semisimplicity. Section~3 reviews
the  construction of the functor ${\cal V}^{FN}$
of Frohman and Nicas and its values  on basic cobordisms. In Section~4 we
describe a convenient set of
generators of the mapping class groups as combinations of Dehn twists and tangles,  
and determine 
their actions on homology. Section~5 introduces the basic rules for the construction of
a Hennings TQFT as well as a method that allows us to construct TQFT's even from
{\em non-modular}
 Hopf algebras or categories. In Section~6 we give the precise definition of  $\cal N$ as 
a quasi triangular Hopf algebra in the sense of Drinfel'd together with the 
$\SL(2,\R)$-action on it. The vector spaces and the basic morphisms of the associated Hennings 
TQFT are computed in Section~7 using standard tangle presentations. 
We prove $\SL(2,\R)$-covariance and single out an index 2
 subcategory of framed cobordisms that naturally
yields a {\em real} valued TQFT. For later applications
we also determine  the {\em categorical}
 Hopf algebra that is canonically associated to this TQFT. 
The nilpotent braided structure of $\cal N$ is then used in Section~8 to develop
a skein theory for the evaluation of tangle diagrams. The pivotal equivalence
of TQFT's that relates this theory to the Alexander polynomial is given by a natural 
isomorphism of  functors as follows. This is
 proven in Section~9 by explicit comparison of generating morphism. 
 
\begin{theorem}\lbl{thm-main}
There is an $\SL(2,\R)$-equivariant isomorphism  
$$
\xi\;:\;\;\; {\cal V}_{\cal N}^{(2)}
\;\;\stackrel{\bullet\,\,\cong}{-\!\!\!-\!\!\!-\-\!\!\!\!\!\longrightarrow}
\;\;{\cal V}^{FN}\;\;,
$$
where both TQFT's are ``non-semisimple'', 
$\Z/2\Z$-projective functors from the category $\Cob^{\bullet}$
of  surfaces with one boundary component and relative cobordisms to the category
of  real $\SL(2,\R)$-modules. 
\end{theorem}

The Hard Lefschetz $\SL(2,\R)$ action on the cohomology of the $U(1)$ moduli spaces and 
its covariance with ${\cal V}^{FN}$ are described more precisely in Section~10. The fact 
that $\xi$ is an $\SL(2,\R)$-equivariant transformation
 is proven. Moreover, we describe the canonical decompositions of the TQFT and the 
Alexander polynomial according to their dual $\SL(2,\R)$-representations. The summands
are irreducible TQFT's for which the mapping class groups are represented by  fundamental
weight representations of the symplectic groups $\Sp(2g,\Z)$. In 
Section~11 we use the equivalence from Section~9 and the skein theory for tangles
from Section~12 to lay out an explicit algorithm, based on a skein theory that extends
the Alexander-Conway calculus, for the computation of  $\Delta_{\varphi}(M)$. 

\begin{theorem}\label{them-skeintheory}
 Let $\cal L$ be a framed link and $\cal Z\subset \cal L$ a distinguished 
component that has zero framing and algebraic linking number zero with all other components.
Let $M_{\cal L}$
be the 3-manifold obtained by surgery along ${\cal L}$ and $\varphi_{\cal Z}:\pi_1(M)\to\Z$
the linking number with $\cal Z$. 

Then $\Delta_{\varphi_{\cal Z}}(M_{\cal L})\in\Z[t,t^{-1}]$ 
can be computed systematically as follows:
\begin{itemize}
\item Use the skein relations from Proposition~\ref{propos-Yskein} to unknot the special
strand $\cal Z$.
\item Put the new configuration into a standard form as depicted in Figure~\ref{fig-stan},
yielding a  tangle $\cal T$.
\item Use the skein relations from Theorem~\ref{thm-skein} and framing relations from
Figure~\ref{fig-frame} to decompose ${\cal T}^{\#}$ into elementary diagrams as
described in in Theorem~\ref{thm-solve}.
\item Translate the elementary tangle diagrams into Hopf algebra diagrams as in 
(\ref{eq-symhopfeva}). 
\item Go through the steps of Proposition~\ref{prop-HopfAlexeval} to assign 
polynomials to each component of a diagram. 
\item Take products over components and sums over elementary diagrams.
\end{itemize}

\end{theorem}

The calculus described here for the evaluation of tangle diagrams is precisely the
one used to compute the morphisms for the TQFT functors from Theorem~\ref{thm-main}
via tangle surgery presentations of cobordisms.

Another application of the equivalence  established in
Theorem~\ref{thm-main} arises from the observation 
that every TQFT ${\cal V}$ on $\Cob^{\bullet}$ naturally implies a braided
Hopf algebra structure ${\cal H}_{\cal V}$ on ${\cal N}_0:={\cal V}(\Sigma_{1,1})$. 
Now, the cohomology ring $H^*(J(\Sigma_g, U(1)))\cong\ext *H_1(\Sigma_g)$ already has
a canonical structure ${\cal H}_{ext}$  of a $\Z/2$-graded Hopf algebra induced by
the group structure on $J(\Sigma_g, U(1))$. It is easy to see that 
${\cal H}_{ext}$ is {\em not} compatible with the Lefschetz $\SL(2,\R)$-action.
However, the  braided Hopf algebra  
structure ${\cal H}_{{\cal V}^{FN}}$ inherited from the TQFT's in 
Theorem~\ref{thm-main} is naturally $\SL(2,\R)$-variant, and, furthermore,
 equivalent to ${\cal H}_{ext}$: 

\begin{theorem}\lbl{thm-structure}
For any choice of an integral Lagrangian decomposition, $H_1(\Sigma_g,\Z)=\Lambda\oplus\Lambda^*$,
and volume forms, $\omega_{\Lambda}\in\ext g\Lambda$ and $\omega_{\Lambda^*}\in\ext g\Lambda^*$, 
the space
$H^*(J(\Sigma_g))$ admits a canonical structure ${\cal H}_{\Lambda}$ of a $\Z/2$-graded
Hopf algebra. It coincides with the braided Hopf algebra structure induced
by ${\cal V}^{FN}$ and is isomorphic to the canonical structure ${\cal H}_{ext}$. 

In particular,
$(H^*(J(\Sigma_g)),{\cal H}_{\Lambda})$ is commutative and cocommutative in the graded
sense, with unit $\omega_{\Lambda^*}$, integral $\omega_{\Lambda}$, and primitive elements given by
$a\wedge\omega_{\Lambda^*}$ and $i^*_z\omega_{\Lambda^*}$ for $a\in H_1(\Sigma)$ and $z\in H^1(\Sigma)$.

The   structure  ${\cal H}_{\Lambda}$ is, furthermore, compatible with the 
Hard-Lefschetz $\SL(2,\R)$-action. Specifically, this action is the Howe dual to 
the action of $\SL(g,\Z)$ on the Lagrangian subspace in the group of Hopf
automorphisms:
$$
\SL(2,\R)_{Lefsch.}\times \SL(\Lambda)\,\;\subset \;\, {\rm GL}(2g,\R)= Aut(H^*(J(\Sigma_g)),{\cal H}_{\Lambda})
$$ 
\end{theorem} 

In Section~13 we   discuss the appearance of these TQFT's in other
contexts. To this end let us denote by  
${\cal V}^{(j)}$ the irreducible component of ${\cal V}^{FN}$ dual to the $j$-dimensional
$\SL(2,\R)$-representation. A detailed  description of it is given in Theorem~\ref{cor-dec}.
Choose for a closed 3-manifold $M$ with Betti number
$\beta_1(M)\geq 1$ a surjection $\varphi:H_1(M)\onto{13}\Z$
(which would be canonical for homology circles as given by 0-surgeries on knots). A series
of  invariants for the pair $(M,\varphi)$ can now be constructed  by choosing any two-sided,
embedded surface $\Sigma \subset M$ that is dual to $\varphi$, and considering the cobordism
$C_{\Sigma}:\,\Sigma\to\Sigma$ obtained by removing an open tubular neighborhood of $\Sigma$
from $M$. The {\em $j$-th (fundamental) Alexander Character} is now defined to be  the integer
\beq\label{eq-AlexMom}
\Delta^{(j)}_{\varphi}(M)\;\;=\;\;trace\Bigl({\cal V}^{(j)}(C_{\Sigma})\Bigr)\;,
\eeq
which is easily seen to depend only on $\varphi$ but not the choice of $\Sigma$. 
Besides the Alexander Polynomial also two other invariants invariant $I^{SW}$ and
$I^{DC}$ depending this data have been constructed by Donaldson in \cite{Don99}
from a Seiberg-Witten Theory and an $SO(3)$-Casson-type gauge theory respectively.
Let us also denote by $\lambda_L$ the Lescop Invariant \cite{Les}. As specified in
the next theorem all of these invariants are in fact linear combinations
of the (fundamental) Alexander Characters. 

\begin{theorem}[{\rm\small mostly corollaries to} {\cite{FroNic92},  
\cite{Don99},\cite{Les},\cite{KerKyoto}}]\label{thm-relations}
\begin{eqnarray}
\Delta_{\varphi}(M)\;\;&=&\;\quad\sum_{j\geq 1}\,[j]_{-t}\cdot \Delta^{(j)}_{\varphi}(M)\label{eq-AlexChar}\;
\\
&&\nonumber\\ 
I^{DC}_{\varphi}(M)\;\;&=&\;\;\;\sum_{j\geq 2}\,{{j+1}\choose 3}\cdot \Delta^{(j)}_{\varphi}(M)\label{eq-DCChar}\;
\\
&&\nonumber\\
I^{SW}_{d,\varphi}(M)\;\;&=&\;\;\sum_{j\geq d+2}\,\left[\!\!{\left[
\Bigl({\frac{j-d}2}\Bigr)^2\right]}\!\!\right]\cdot
\Delta^{(j)}_{\varphi}(M)\label{eq-SWChar}
\\
&&\nonumber\\
\lambda_{L}(M)\;\;&=&\;\;\sum_{j\geq 1} (-1)^{j-1}\frac {j(2j^2-3)}{12}\cdot
\Delta^{(j)}_{\varphi}(M)\;\label{eq-LesChar}
\end{eqnarray}
\end{theorem}
Here we denoted $[j]_q=\frac {q^j-q^{-j}}{q-q^{-1}}$ and by $\bigl[\![x]\!\bigr]$ the
largest integer $\leq x$. We further review in how far the higher $PSU(n)$ knot invariants 
$I_{k,n,\varphi}^{FN}$
of Frohman and Nicas \cite{FroNic94} come out to be 
polynomial expressions in the Alexander Characters.
As products of characters are associated to tensor products of TQFT's and their decompositions
into irreducible components it is natural to  consider the corresponding higher,
irreducible Alexander Characters $\Delta^{(\gamma)}$. We conjecture that the  
$I_{k,n,\varphi}^{FN}$ are linear combinations of the $\Delta^{(\gamma)}$ with 
coefficients in ${\mathbb N}\cup\{0\}$ as it is the case for $I^{DC}$ and $I^{SW}$. 


\emptystuff{
 From this sequence of TQFT's  we are able to construct a large family of homological 
TQFT functors ${\cal V}^{(P)}$ for
suitable polynomials $P$ in variables $x_0, x_1, \ldots$ by taking tensor products and direct sums
accordingly. 

 One example of such a TQFT ${\cal V}^{(Y)}$ with $Y=x_0+x_3-x_5-x_8+x_{10}+\ldots\,$
taken over $\Z/5$ turns out describe the lowest order contribution
of the Reshetikhin Turaev invariant taken over the cyclotomic integers $\Z[\zeta_5]$ at least 
in small genera and very likely in general \cite{Kerr=5}. A rather interesting application that
emerges from that is found in in joint work with  Gilmer \cite{GilKer}. Namely, that non triviality
of the Alexander polynomial evaluated at a 5-th root of unity implies that a 3-manifold disconnects
if a surface with more than one component is removed from it. 
 
 The polynomial homological TQFT's appear to be also isomorphic to ones
 constructed by gauge theoretical means
using an approach similar to the Casson invariant or extending the methods
of Seiberg-Witten theory as described by Donaldson. Furthermore,  Frohman and Nicas consider
generalizations to higher rank Lie groups. In all cases the Alexander polynomial appears
as the dominant invariant, suggesting the corresponding decomposition of the TQFT into
 the basic functors ${\cal V}^{(j)}$. 
}


Moreover, we explain the irreducible  $p$-modular reductions  $\dov{\cal V}^{(j)}_p$
over ${\mathbb F}_p=\Z/p\Z$ of the ${\cal V}^{(j)}$ relate to the irreducible factors
of the $\Z[\zeta_p]\to{\mathbb F}_p$ of the
Reshetikhin Turaev TQFT's at a $p$-th root of unity $\zeta_p$. 
 We finally give evidence that  the TQFT from Theorem~\ref{thm-main}
is essentially the missing tensor factor that relates the semisimple and the non-semisimple
TQFT constructions for $U_q({\mathfrak s}{\mathfrak l}_2)$
following Reshetikhin Turaev and Hennings respectively. 

\medskip

\paragraph{Acknowledgements:} I am  indebted to Charlie Frohman for making me aware 
of \cite{FroNic92}  and explaining  \cite{FroNic94} to me. I thank 
Bernhard Kr\"otz for discussions about Howe pairs, and David Radford for helping me find
his example in \cite{Rad76}. Thanks also to Andrew Nicas and Hans Boden for their
interest and  Pierre Deligne, Daniel Huybrechts, and Manfred Lehn
for discussions about Lefschetz decompositions in the higher rank case.
Finally, I want to thank Razvan Gelca, Pat Gilmer, Jozef Przytycki 
for opportunities to speak about this paper, when it was still in its
early stages. 
\newpage 

\head{2. Topological Quantum Field Theory}\lbl{S2}
We start with the definition of  a TQFT as a functor as proposed by  
Atiyah \cite{Ati88}, largely suppressing a more detailed discussion 
of the tensor structures.

For every integer, $g\geq 0$, choose a compact, oriented
 model surface, $\Sigma_g$, of genus $g\,$, and to a tuple of 
integers $\underline g=(g_1,\ldots,g_n)$ associate the ordered union
$\Sigma_{\underline g}:=\Sigma_{g_1}\sqcup\ldots\sqcup\Sigma_{g_n}$.
A {\em cobordism}  
is a collection, ${\bf M}=(M,\phi_{\#},\Sigma_{g_{\#}})$, of the following:

 A compact, oriented 3-manifold, $M$, whose boundary  is divided into 
two components 
$\partial M=-\partial_{in} M\sqcup\partial_{out} M$,  
two standard surfaces $\Sigma_{\underline g_{in}}$ and 
$\Sigma_{\underline g_{out}}$,   and two orientation preserving 
homeomorphisms $\phi_{in}:\Sigma_{\underline g_{in}}
\,\isto\, \partial_{in} M\,$  and $\phi_{out}:\Sigma_{\underline g_{out}}
\,\isto\, \partial_{out} M\,$. 
 
We say two cobordisms, ${\bf M}$ and ${\bf M}'$, are equivalent if
they have the same "in" and "out" standard surfaces, and  
there is a homeomorphism $h:\,M\,\isto \, M'\,$, such that 
 $h\circ \phi_{\#}=\phi_{\#}'\,$. 

Let $\Cob$ be the category of cobordisms in dimension 2+1, which has
the standard surfaces as objects and equivalence classes of cobordisms
as morphism. The composition of morphisms is defined via gluing over
boundary components using the coordinate maps to the same standard
surfaces.  In addition, $\Cob$ has a tensor
product given by disjoint unions of surfaces and cobordisms.

A {\em Topological Quantum Field Theory} (TQFT)   is  a functor,
${\cal V}:\Cob\longrightarrow{\rm Vect}(\kk)$, from the category of cobordisms
to the category of vector spaces over a field $\kk\,$. 
\medskip

  Let us recall next some generalizations of the definition given in
\cite{Ati88} that will be relevant for our purposes. 
By $\Cob^{2fr}$ we denote the category of 2-framed cobordisms, where we 
fixed some standard framings on the model surfaces $\Sigma_g\,$, see \cite{Ker99}. A 
{\em 2-framed TQFT} is now a functor
 ${\cal V}:\Cob^{2fr}\longrightarrow{\rm Vect}(\kk)$. The category of 2-framed  
cobordisms can be understood as a central extensions 
\beq\label{eq-cobext}
\,1\;\to\;\Z\;\longrightarrow\;\Cob^{2fr}\;\longrightarrow\;\Cob\;\to \;1\,
\eeq 
of the 
ordinary cobordism category, if restricted to connected cobordisms. Hence, an
irreducible 2-framed TQFT yields a {\em projective TQFT} 
since $\Z$ is presented as a scalar. 
See \cite{Ker99} for further descriptions of this extension in terms of signatures of
bounding 4-manifolds.

For a group, $G$, we introduce the notion of a {\em $G$-equivariant TQFT}.
It is a functor, 
${\cal V}:\Cob\longrightarrow G-{\rm mod}_{\kk}$, from the category of 
cobordisms to the category of finite dimensional $G$-modules over a field $\kk$. This means
that the linear map associated to any cobordism commutes with 
the action of $G$ on the vector spaces of the respective boundary
components. 

Recall also from \cite{Ker98} that a {\em half-projective} or {\em non-semisimple} 
TQFT is one in which functoriality is weakened and replaced by 
the composition law ${\cal V}(MN)=0^{\mu(M,N)}{\cal V}(M){\cal V}(N)\,$.
Here $\mu(M,N)=b(MN)-b(M)-b(N)\in\Z^{+,0}\,$, where $b(M)$ is the 
number of components of $M$ minus half the number of components of
$\partial M$. Note that $0^0=1$. 

We often call a cobordism for which  all (rational) homology comes from the homology
of the boundary {\em (rationally) homologically trivial} (r.h.t). More precisely, we mean by this 
that $i_*:H_1(\partial M,\Q)\to H_1(M,\Q)$ is onto. Typical examples of r.h.t.
cobordisms are the ones in (\ref{eq-defpsi}) and (\ref{eq-hdlcob})
below and closed, rational homology spheres. Examples of cobordisms that are not 
r.h.t. are any connected sums with closed manifolds $M$ with $\beta_1(M)\geq 1$. 
We find the following vanishing property:
\begin{lemma}[\cite{Ker98}]\lbl{lm-vanish}
If ${\cal V}$ is a non-semisimple TQFT, then for any cobordism $M$, 
$$
{\rm if}\quad {\cal V}(M)\, \neq \, 0\qquad{\rm then}\qquad 
M\;\;\;\mbox{is r.h.t.}
$$
\end{lemma} 

We further introduce $\Cob^{\bullet}$, the  category of cobordisms, for which 
the surfaces are connected and have exactly one boundary component. As objects we
thus use model surfaces $\Sigma_{g,1}$, such that $\Sigma_{g+1,1}$ is obtained
from $\Sigma_{g,1}$ by gluing in a torus, $\Sigma_{1,2}$, with two boundary 
components. Thus, we have a presentation  
\begin{equation}\lbl{eq-modelsurf}
\Sigma_{g,1}\;=\;\underbrace{\Sigma_1\#\ldots\#\Sigma_1\#\Sigma_{1,1}}_g
\qquad\qquad
\mbox{with inclusions}\qquad\quad \Sigma_{g,1}\;\subset\;\Sigma_{g+1,1}\;.
\end{equation}
 Instead of ordinary cobordisms we then consider {\em relative} ones. We finally introduce 
categories of cobordisms with 
combinations of these properties such as  $\Cob^{2fr, \bullet}$, the category of
2-framed, relative cobordisms. 
\medskip

 For any homeomorphsim,  $\psi\in Homeo^+(\Sigma_g)$, of a surface to itself 
we define the cobordism 
\beq\lbl{eq-defpsi}
{\bf I}_\psi\;\;=\;\;(\Sigma_g\times[0,1], id\sqcup \psi, \Sigma_g\sqcup \Sigma_g)\quad. 
\eeq
The morphism $[{\bf I}_{\psi}]$
depends only on the isotopy class $\{\psi\}$ of $\psi$,  and the resulting map
$\Gamma_g\to Aut(\Sigma_g):\{\psi\}\mapsto [{\bf I}_{\psi}]$ from the mapping class 
group to the group of invertible
cobordisms on $\Sigma_g$ is an isomorphism, see \cite{KerLub00}. Consequently,
every TQFT defines a representation  of the mapping class group 
$\Gamma_g\to GL({\cal V}(\Sigma_g)):\,
\{\psi\}\mapsto {\cal V}([{\bf I}_{\psi}])\,$.

Moreover, let us introduce  special cobordisms
\beq\lbl{eq-hdlcob}
{\bf H_{g}^+}\;:=\;=( H_{g}^+,id\sqcup id, \Sigma_g\sqcup \Sigma_{g+1})\;,
\eeq
where $H^+_g$ is obtained by  adding a 
full 1-handle to the cylinder $\Sigma_g\times[0,1]$  
at two discs in $\Sigma_g\times 1$.
This is done in a way compatible with the choice of the model surfaces in
equation~(\ref{eq-modelsurf}). Another cobordism $H_{g}^-$ is built by gluing in
a 2-handle into the thickened surface $\Sigma_{g+1}\times[0,1]$ along a
curve  $b_{g+1}$, which lies in the added torus from (\ref{eq-modelsurf})
and  has geometric intersection number one with the meridian
of the 1-handle added by  $H_g^+$. From this we obtain a cobordism 
${\bf H_g^-}=( H_g^-,\Sigma_{g+1}\sqcup\Sigma_g)$ in opposite 
direction, with the property that ${\bf H_{g}^-}\circ {\bf H_{g}^+}$ is 
equivalent to the identity.

Basic Morse theory implies a Heegaard decomposition as follows for any cobordism 
\beq\lbl{eq-Heegaard}
{\bf M}\;\;\;\cong\;\;\;{\bf H_{g_2}^-}\circ {\bf H_{g_2+1}^-}\circ\ldots
\circ {\bf H_{N-1}^-}\circ {\bf I}_{\psi}\circ {\bf H_{N-1}^+}\circ\ldots\circ 
{\bf H_{g_1+1}^+}
\circ {\bf H_{g_1}^+}\;,
\eeq
where $\psi\in Homeo^+(\Sigma_N)$. Hence, a  TQFT is completely determined by
the induced representations of the mapping class groups and the maps 
${\cal V}([{\bf H}^+_g])$ and ${\cal V}([{\bf H}^-_g])$. Therefore, any two TQFT's
coinciding on  the basic generators from (\ref{eq-defpsi}) and  (\ref{eq-hdlcob})
have to be equal. 

\head{3. The Frohman-Nicas TQFT for $U(1)$}\lbl{S3}

Let us review the basic steps in the construction of the topological
quantum field theory ${\cal V}^{FN}$ as given in \cite{FroNic92}
via intersection theory of $U(1)$-representation varieties:

For a compact, connected manifold $X$  its  $U(1)$-representation variety is
defined as
\beq\lbl{eq-defJX}
J(X)\quad:=\quad {\rm Hom}(\pi_1(X), U(1))\quad \cong \quad H^1(X,U(1))\;.
\eeq
Observe that $J(X)$ is  a manifold of dimension $\beta_1(X)$. Specifically, it is
a torus if $H_1(X,\Z)$ is torsion free, and a discrete group if $\beta_1(X)=0$. 

The vector space associated  to a surface $\Sigma_{\underline g}$ is given
by ${\cal V}^{FN}(\Sigma_{\underline g})
= H^*(J(\Sigma_{g_1})\times\ldots\times J(\Sigma_{g_N}) ,\R)$. 

We consider first  cobordisms, $M$, between  surfaces, $\partial_{in}M$ 
and $\partial_{out}M$,  that are  rationally homologically trivial in the
sense of Section~2. 
In this case the map $j:\,J(M)\to J(\partial_{in}M)\times J(\partial_{out}M)$ is a half
dimensional immersion. Thus the top form $\pm[J(M)]$ defines (up to sign)
a middle dimensional 
homology class in $H_*(J(\partial_{in}M),\R)\otimes H_*(J(\partial_{out}M),\R)$.
Using Poicar\'e Duality and the coordinate maps of the cobordism, 
the latter space is isomorphic to the space of linear maps
from ${\cal V}^{FN}(\Sigma_{\underline g_{in}})$ to 
${\cal V}^{FN}(\Sigma_{\underline g_{out}})$. ${\cal V}^{FN}(M)$, for a  
 homologically  trivial cobordism $M$, is now the linear map  associated to 
 $j_*(\pm[J(M)])$.

In the general case Frohman and Nicas define   ${\cal V}^{FN}(M)$ via
a Heegaard splitting of $M$ as in (\ref{eq-Heegaard}), and consider the intersection 
number of representation
varieties of the elementary thick surfaces with handles separated by the Heegaard surface. 
In the case where $H_1(\partial M,\R)\to H_1(M,\R)$ is not onto, i.e., $M$ is
not homologically trivial, these varieties
no longer transversely intersect so that ${\cal V}^{FN}(M)=0$. 
 
Regarding the composition structure ${\cal V}^{FN}$ has a couple of nonstandard
properties. For one, functoriality fails to hold when $M$ and $N$ are homologically
trivial but $M\circ N$ is not. Moreover, the orientation of the classes $\pm[J(M)]$ and cycles
cannot be chosen  consistently with composition so that a sign-projectivity persists.
Recall, however, that a 2-framed TQFT is really defined on the $\Z$-extensions of cobordisms
given in (\ref{eq-cobext}).

\begin{lemma}\lbl{lm-FNpart}
${\cal V}^{FN}$ is  a non-semisimple, $\Z/2\Z$-projective TQFT in the sense of Section~2.
\end{lemma} 

The mechanism by which the universal $\Z$-extension is factored into  a $\Z/2\Z$-extension
is explained further for the quantum theory in Lemma~\ref{lm-mod2prop} and 
Proposition~\ref{pp-ATQFT} of Section~7. At least indirectly, we have thus related the
orientation ambiguities in \cite{FroNic92} to the usual framing ambiguities of quantum
theories. 

Now, in the $U(1)$ case $J(X)$ has a group structure itself, which induces a 
coalgebra structure  on the cohomology ring so that $H^*(J(X))$ is endowed
with a canonical Hopf algebra structure ${\cal H}_{ext}$. If $H_1(X)$ is torsion free 
then $H^*(J(X))$ is connected and we obtain a natural isomorphism
$H^*(J(X))\cong\ext * H_1(X)$ of $\Z/2$-graded Hopf algebras, and
$H_1(X)$ is the space of primitive elements. Hence, we can write for the vector spaces:
\begin{equation}
{\cal V}^{FN}(\Sigma_{\underline g})\;=\;\ext * H_1(\Sigma_{\underline g})\;.
\end{equation}
The representation of the mapping class group $\Gamma_{\underline g}$
on this space is given by the obvious action 
\beq\lbl{eq-FNmapcg}
{\cal V}^{FN}([{\bf I_{\psi}}])\;\;=\;\;\ext * [\psi]\qquad
\forall \{\psi\}\in\Gamma_g\;.
\eeq
Here, $[\psi]\in \Sp(H_1(\Sigma_{\underline g}))$ is the natural, induced action on homology.
For a connected surface $\Sigma_g$ we have the associated short exact sequence
\beq\lbl{eq-TGSseq}
1\;\to\;{\cal J}_g\;\longrightarrow\;\Gamma_g\;\stackrel{\psi\mapsto[\psi]}
{-\!\!\!-\!\!\!-\!\!\!-\!\!\!\longrightarrow}\;\Sp(2g,\Z)\;\to\;1\;,
\eeq
where ${\cal J}_g$ is the Torelli group. 

Let  
${\bf H_{g}^+}$ be the cobordism as defined in (\ref{eq-hdlcob}), and 
let $[a_{g+1}]$ be a generator of  
$ker(H_1(\Sigma_{g+1},\Z)\to H_1(H_{g+},\Z))$ seen as an element 
of $H_1(\Sigma_{g+1},\R)$. It is represented by the meridian $a_{g+1}$ of
the added handle. 
In a slight variation of the Frohman Nicas formalism we see that the associated
linear map is given as 
\beq\lbl{eq-FNmcg}
{\cal V}^{FN}({\bf H_{g}^+}):\ext * H_1(\Sigma_g)\longrightarrow\ext * H_1(\Sigma_{g+1})\;:
\quad \alpha\;\mapsto\;i_*(\alpha)\wedge [a_{g+1}]\;.
\eeq
Here we use the fact that  $H_1(\Sigma_{g,1})=H_1(\Sigma_g)$ so that the 
inclusion of surfaces in (\ref{eq-modelsurf}) implies also an inclusion
$i_*:H_1(\Sigma_g)\subset H_1(\Sigma_{g+1})$.

Let ${\bf H_{g}^-}$ be the cobordism obtained by gluing a 2-handle along $b_{g+1}$
as defined above. We note that 
$H_1(\Sigma_{g+1})=H_1(\Sigma_{g})\oplus\lz[a_{g+1}],[b_{g+1}]\rz$ so that 
$\ext * H_1(\Sigma_{g+1})$ is the direct sum of spaces 
$V_1\oplus V_a\oplus V_b\oplus V_{a\wedge b}$ where 
$V_{x}=[x_{g+1}]\wedge \ext *H_1(\Sigma_{g})$. The linear map associated in \cite{FroNic92}
to ${\bf H_{g}^-}$ acts on $V_a$ as 
\beq\lbl{eq-FNgcm}
{\cal V}^{FN}({\bf H_{g}^-}): V_a\longrightarrow\ext * H_1(\Sigma_{g})\;:
\quad i_*(\alpha)\wedge [a_{g+1}]\;\mapsto\;\alpha
\eeq
and is zero on all other summands.



\head{4. The Mapping Class Groups and their Actions on Homology}\lbl{S4}
The mapping class group
$\Gamma_{g,1}=\pi_0(Homeo^+(\Sigma_{g,1}))$ on a model surface
$\Sigma_{g,1}$ is generated by the right handed 
Dehn twists along oriented curves $a_j$, $b_j$, 
and $c_j$,   as depicted in
Figure~\ref{fig-modsurf}. We denote them by capital letters
 $A_j,\,B_j,\,C_j\, \in\Gamma_{g,1}$ respectively.
In fact we only need the $A_2$ of the $A_j$'s to generate 
$\Gamma_{g,1}$. A presentation of $\Gamma_{g,1}$ in 
these generators is given by Wajnryb \cite{Waj83}. 
\begin{figure}[ht]
\begin{center}
\leavevmode
\epsfxsize=7cm
\epsfbox{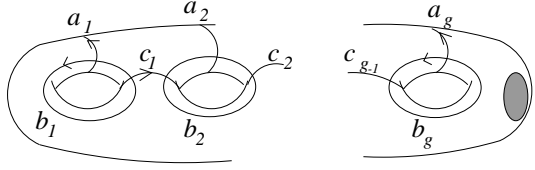}
\end{center}
\caption{Curves on $\Sigma_{g,1}$}\lbl{fig-modsurf}
\end{figure}
For our purposes we prefer the  set $\{A_j, D_j, S_j\}$ of generators
defined as follows: 
\begin{equation}\lbl{eq-newgen}
D_j\;:=\;A_j^{-1}A_{j+1}^{-1}C_j\qquad\mbox{and}\qquad S_j\;:=\;A_jB_jA_j\qquad\quad
\mbox{for}\;\;j=1,\ldots,g \;.
\end{equation}
In \cite{MatPol94} a tangle presentation of $\Gamma_{g,1}$ is given using the results
in \cite{Waj83}. The same presentation results from the tangle presentation of 
$\Cob^{2fr,\bullet}$ in \cite[Proposition 14]{Ker99}, which extends to the central
extension $1\to\Z\to\Gamma_{g,1}^{2fr}\to\Gamma_{g,1}\to 1$ that stems from the 
2-framing of cobordisms. The framed tangles associated to our preferred generators
are given in Figures~\ref{fig-A-tgl}, \ref{fig-D-tgl}, and \ref{fig-S-tgl}.
We use an empty circle to  indicate a right handed $2\pi$-twist on the framing 
of a strand as in Figure~\ref{fig-A-tgl}, and a full circle for a left handed one as in 
Figure~\ref{fig-Rigg-tgl}. Note, that the extra 1-framed circle in
Figure~\ref{fig-S-tgl} does not change the 3-cobordism in $\Cob^{\bullet}$  but 
shifts its 2-framing in $\Cob^{2fr, \bullet}$  by one. 

\begin{figure}[ht]
\begin{center}
\leavevmode
\epsfysize=2.2cm
\epsfbox{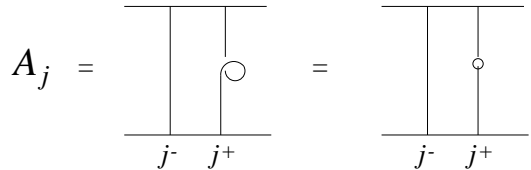}
\end{center}
\caption{Tangle for $A_j$}\lbl{fig-A-tgl}
\end{figure}

\begin{figure}[ht]
\begin{center}
\leavevmode
\epsfysize=2.2cm
\epsfbox{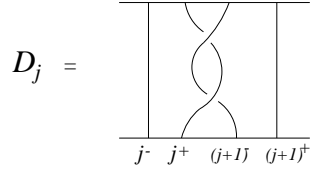}
\end{center}
\caption{Tangle for $D_j$}\lbl{fig-D-tgl}
\end{figure}

\begin{figure}[ht]
\begin{center}
\leavevmode
\epsfysize=2.2cm
\epsfbox{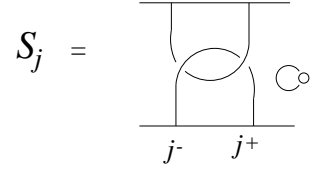}
\end{center}
\caption{Tangle for $S_j$}\lbl{fig-S-tgl}
\end{figure}
$\Gamma_{g,1}^{2fr}$ can then be thought of as the sub-group of tangles generated
by these diagrams, modulo isotopies, 2-handle slides, the $\sigma$-move and the
Hopf link move, see \cite{Ker99}.  
\medskip

For later purposes we give the explicit action of these generators on
$H_1(\Sigma_g,\Z)=H_1(\Sigma_{g,1},\Z)$ in the sense of (\ref{eq-TGSseq}). 
Suppose $p,\,f\,\subset\Sigma_{g,1}$ 
are two transverse, oriented curves. We denote by $P$ the Dehn twist along $p$,
by $[P]\in \Sp(2g,\Z)$ its action on homology, 
and by $[p]$ and $[f]$ the respective homology classes. We have
\begin{equation}\lbl{eq-DehnHom}
[P].[f]\quad=\quad[f]\,+\,([p]\cdot[f])[p]\;.
\end{equation}
Here $([p]\cdot[f])\in\Z$ is the algebraic intersection number of $p$ with $f$,
counting $+1$ for  a crossing if the tangent vectors of $p,f$ form an oriented
basis and $-1$ if the basis has opposite orientation. 

A basis for $H_1(\Sigma_g)$ is given by 
$\{[a_1],\ldots, [a_g],[b_1],\ldots, [b_g]\}$,
 and intersection numbers can be read
off Figure~\ref{fig-modsurf}. 
For example $a_j$ intersects $b_j$ in only one point, where 
$[a_j]\cdot[b_j]= +1$ since $b_j$ follows $a_j$ counter clockwise at the 
crossing. Hence 
\begin{equation}\lbl{eq-A-hom}
[A_j].[b_j]=[b_j]+[a_j]\qquad \mbox{and}\quad [A_j].[x]=[x]\quad  
\mbox{for all other basis vectors.}
\end{equation}
Similarly, we have that $[C_j]$ only acts on $[b_j]$ and $[b_{j+1}]$
with $[C_j].[b_j]=[b_j]+[c_j]$ and $[C_j].[b_{j+1}]=[b_{j+1}]-[c_j]$. 
Substituting $[c_j]=[a_j]-[a_{j+1}]$, and using the definition  of $D_j$ in
(\ref{eq-newgen}) and  (\ref{eq-A-hom}) we compute
\begin{equation}\lbl{eq-D-hom}
[D_j].[b_j]=[b_j]-[a_{j+1}]\quad \mbox{and}\quad [D_j].[b_{j+1}]=[b_{j+1}]-[a_{j}] \;,
\end{equation}
and, again,  $[D_j].[x]=[x]$ for all other basis vectors $[x]$ of $H_1(\Sigma_1,\Z)$.
Finally, we find $[B_j].[a_j]=[a_j]-[b_j]$ so that
\begin{equation}\lbl{eq-S-hom}
[S_j].[a_j]=-[b_j] \quad \mbox{and}\quad [S_j].[b_j]=[a_{j}] 
\end{equation}
and $[S_j].[x]=[x]$ elsewise. 

The above action can be identified with specific generators  of
the Lie algebra ${\mathfrak s}{\mathfrak p}(2g,\R)$ as follows:
$$
[A_j]\;
=\;I_{2g}\;+\;E_{j,-j}\;=\;I_{2g}\;+\;e_{2\epsilon_j}\;=\;\exp(e_{2\epsilon_j})
$$
\begin{equation}\lbl{eq-sproots}
[B_j]\;=\;I_{2g}\;-\;E_{-j,j}\;=\;I_{2g}\;-\;f_{2\epsilon_j}\;=\;\exp(-f_{2\epsilon_j})
\end{equation}
$$
[D_j]\;=\;I_{2g}\;-\;E_{j,-(j+1)}-E_{j+1,-j}\;=
\;I_{2g}\;-\;e_{\epsilon_j+\epsilon_{j+1}}\;=\;\exp(-e_{\epsilon_j+\epsilon_{j+1}})\\
$$
The conventions and notations for the weights $\epsilon_j$ and the 
matrices $E_{i,j}$ are taken from \cite[Chapter 2.3]{GooWal98}. Hence,
the natural representation on ${\rm Sp}(2g,\Z)$ clearly lifts   to the
fundamental representation of ${\rm Sp}(2g,\R)$.

Finally, there is an 
 ${\rm Sp}(2g,\Z)$-invariant 2-form, which is unique up to signs
and given in our basis as:
\begin{equation}\lbl{eq-defsymform}
\omega_g\quad:=\quad \sum_{j=1}^g[a_j]\wedge[b_j]\qquad\in\;
\ext 2 H_1(\Sigma_g)\,=\,H^2(J(\Sigma_g))\;.
\end{equation}
It is identical to twice the K\"ahler metric form in
$H^2(J(\Sigma_g))$, see Section~10 and \cite{GriHar78}.

\head{5. Hennings TQFT's}\lbl{S5}
In \cite{Hen96} Hennings describes a calculus
that allows us to compute an invariant, ${\cal V}_{\cal A}^H(M)$,
for a closed 3-manifold, $M$, starting from a
surgery presentation, $M=S^3_{\cal L}$,  by a framed link, ${\cal L}\subset S^3$, and 
 a quasitriangular Hopf algebra $\cal A$. It is obtained by inserting and moving elements
of $\cal A$ along the strands of a projection of $\cal L$ and evaluating them against
integrals. This procedure was  refined by Kauffman and 
Radford \cite{KauRad95} permitting
 unoriented links and simplifying the evaluation and proofs substantially. 
${\cal V}^H_{\cal A}$ turns out to be a special case of the invariant given by  
Lyubashenko \cite{Lub95}, which is constructed from general abelian categories. 
In \cite[Theorem 14]{Ker96} we generalize the Hennings procedure to 
tangles and cobordisms and thus construct a topological 
quantum field theory ${\cal V}^H_{\cal A}$ for any modular Hopf algebra $\cal A$. 
In turn ${\cal V}^H_{\cal A}$ is derived as a special case of
the general TQFT construction by Lyubashenko and 
the author in \cite{KerLub00}. 

The TQFT in  \cite{Ker96} was formulated as
a contravariant functor, ${\cal V}_{\cal A}^{*}:\,Cob_3^{\bullet}\to Vect(\kk)$,
where ${\cal V}_{\cal A}^{*}(\Sigma_{g,1})={\cal A}^{\otimes g}$. In this section we
will  give  the rules for construction for the covariant version,  defined by 
${\cal V}_{\cal A}(M)=(f^{\otimes g})^{-1}({\cal V}^*_{\cal A}(M))^*f^{\otimes g}\,$,  
where  $f:{\cal A}\to{\cal A}^*:\, x\mapsto \mu(S(x)\_\_)$.  We generalize 
\cite{Ker96} further by allowing Hopf algebras, $\cal A$, that are not modular,
at the expense of reducing the vector space by a canonical projection.

Let $M$ be a 2-framed cobordism between two model surfaces,
$\Sigma_{g_1}$ and $\Sigma_{g_2}$.
As in \cite{Ker99} we associate to the homeomorphism class of $M$  an equivalence 
class of framed tangle diagrams. The projection of a  representative tangle,  $T_M$, in
$\R\times [0,1]$  has $2g_1$  endpoints $\,1^-<1^+<2^-<\ldots<g_1^-<g_1^+\,$
in the  top line $\R\times 1$ and $2g_2$ endpoints  $\,1^-<1^+<2^-<\ldots<g_2^-<g_2^+\,$
in the bottom line  $\R\times 0$.  Besides closed components ($\cong S^1$) the tangle can
have components with boundary ($\cong[0,1]$). An interval component, $J$, of the tangle 
can either run between points $j^-$ and $j^+$ at the top line or between 
 $j^-$ and $j^+$ at the bottom line.
 As a forth possibility we admit pairs of components, $I$ and $J$, of
which each starts at the top line and ends at the bottom line and cobords a pair
$\{j^-,j^+\}$ to a pair $\{k^-,k^+\}$. The equivalences of tangles are generated 
by isotopies, 2-handle slides (second Kirby move) over closed components,  the 
addition and removal of an isolated Hopf link, in which one component has 0-framing, and
additional boundary moves, called $\sigma$- and $\tau$-Moves, see  \cite{Ker99}. 
For later purposes we also depict here the  $\sigma$-Move: 

\beq\label{eq-boundmove}
\epsfbox{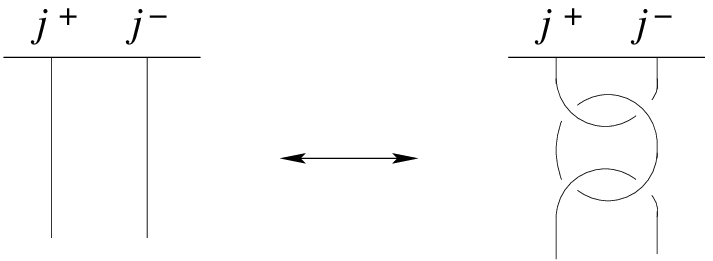}
\eeq

The next ingredient is a unimodular, ribbon Hopf algebra, ${\cal A}$,  
in the sense of \cite{Rad94}, over a
perfect  field $\kk$ with $char(\kk)=0$. 
In particular,  ${\cal A}$ is a {\em quasitriangular} Hopf algebra
as introduced by Drinfel'd \cite{Dri90}. This means there exists   an element 
${\cal R}=\sum_je_j\otimes f_j\,\in{\cal A}^{\otimes 2}$, called the {\em R-matrix},
which fulfills several natural conditions. As in \cite{Dri90} we define the element
$u=\sum_jS(f_j)e_j$, which implements the square of the antipode $S$ by $S^2(x)=uxu^{-1}$.
A {\em ribbon} Hopf algebra is now a quasitriangular Hopf algebra with a group like
element, $G$, such that $G$ also implements $S^2$ and  $G^2=uS(u)^{-1}$. 
From this we define the ribbon element $v:=u^{-1}G$, which is central in
${\cal A}$. Furthermore, it satisfies the equation 
\begin{equation}\lbl{eq-M-matrix}
\,{\cal M}\;\;=\;\;{\cal R}^{\dagger}{\cal R}\;\;=\;\;\Delta(v^{-1})v\otimes v\;\;,
\end{equation}
where $(a\otimes b)^{\dagger}=b\otimes a$ is the transposition of tensor factors.

 Now, any finite dimensional Hopf algebra contains a {\em right integral},
which is an element $\mu \in {\cal A}^*$ characterized by the equation:
\beq\lbl{eq-int-def}
( \mu\otimes id_{\cal A})(\Delta(x))\;=\;1\cdot\mu(x)
\eeq
Its existence and uniqueness (up to scalar multiplication) has been proven
in \cite{LarSwe69}. The adjective  ``unimodular''
implies that 
\beq\lbl{eq-muSinvar}
\mu(xy)\;\;=\;\;\mu(S^2(y)x)\qquad\qquad\mbox{and}\qquad\qquad\mu(S(x))\;\;=\;\;\mu(G^2x)\;,  
\eeq
see \cite{Rad94}. For the remainder of this article  we  will also
 assume  the following normalizations:
\beq\lbl{eq-balnorm}
\mu\otimes\mu({\cal M})\;\;=\;\;1\qquad\mbox{and}\qquad\;\;\mu(v)\mu(v^{-1})=1
\eeq

The next step in the Hennings procedure is to replace the tangle projection $T_M$ 
with distinguished over and under crossings by a formal linear combination of copies of 
the  projection $T_M$ in which we do not distinguish between over and under crossings 
but decorate 
segments of the resulting planar curve with elements of $\,{\cal A}\,$. Specifically, 
we replace an over crossing by an indefinite crossing and insert at the two incoming
pieces the elements occurring in the $R$-matrix, and similarly for an under crossing,
  as indicated in the following diagrams.

\beq\lbl{fig-Dec-cross}
\setlength{\unitlength}{0.25mm}
\begin{picture}(470,50)
\put(20,5){\line(1,1){40}}
\put(20,45){\line(1,-1){17}}
\put(60,5){\line(-1,1){17}}

\put(85,30){}
\put(80,25){\vector(1,0){26}}

\put(125,20){$\displaystyle \sum_j$}

\put(165,50){\line(1,-1){50}}
\put(165,0){\line(1,1){50}}
\put(180,44){$\scriptscriptstyle e_j$}
\put(175,40){\circle*{10}}
\put(205,40){\circle*{10}}
\put(214,38){$\scriptscriptstyle f_j$}


\put(280,5){\line(1,1){17}}
\put(320,45){\line(-1,-1){17}}
\put(320,5){\line(-1,1){40}}

\put(333,25){\vector(1,0){27}}

\put(375,20){$\displaystyle \sum_j$}

\put(425,50){\line(1,-1){50}}
\put(425,0){\line(1,1){50}}
\put(435,10){\circle*{10}}
\put(405,15){$\scriptscriptstyle S(e_j)$}
\put(465,10){\circle*{10}}
\put(482,13){$\scriptscriptstyle f_j$}
\end{picture}
\eeq

The elements on the segments of 
the planar diagram can then be moved along the connected components according to
the following rules. 

\beq\lbl{eq-elem-mv}
\setlength{\unitlength}{0.25mm}
\begin{picture}(440,85)(0,-15)
\put(20,5){\line(0,1){60}}
\put(20,20){\circle*{10}}
\put(27,21){$x$}
\put(20,50){\circle*{10}}
\put(27,51){$y$}

\put(37,35){$=$}

\put(60,5){\line(0,1){60}}

\put(60,35){\circle*{10}}
\put(67,36){$xy$}

\put(135,35){\oval(40,40)[b]}
\put(115,35){\line(0,1){20}}
\put(155,35){\line(0,1){20}}
\put(155,42){\circle*{10}}
\put(162,43){$S(x)$}

\put(190,33){$=$}

\put(235,35){\oval(40,40)[b]}
\put(215,35){\line(0,1){20}}
\put(255,35){\line(0,1){20}}
\put(215,42){\circle*{10}}
\put(222,43){$x$}

\put(300,55){\line(1,-1){40}}
\put(300,15){\line(1,1){40}}
\put(330,25){\circle*{10}}
\put(340,25){$x$}

\put(358,35){$=$}

\put(380,55){\line(1,-1){40}}
\put(380,15){\line(1,1){40}}
\put(390,45){\circle*{10}}
\put(397,48){$x$}

\end{picture}
\eeq

Finally, every diagram can be untangled using the local moves given below,
and the usual planar third Reidemeister move. In particular, undoing a closed
curve in the diagram  
yields an extra overall  factor $G^d$, where $G$ is the  group like element defined
above and $d$ the Whitney number of the curve. 

\beq\lbl{fig-st-isom}
\setlength{\unitlength}{0.25mm}
\begin{picture}(400,100)(0,-15)
\put(30,30){\oval(40,40)[b]}
\put(10,30){\line(2,5){6}}
\put(50,30){\line(-2,5){6}}
\put(16,45){\line(1,1){34}}
\put(44,45){\line(-1,1){34}}

\put(72,43){$=$}

\put(120,30){\oval(40,40)[b]}
\put(100,30){\line(0,1){46}}
\put(140,30){\line(0,1){46}}
\put(140,55){\circle*{10}}
\put(147,56){$G$}


\put(220,40){\line(0,1){20}}
\put(250,40){\line(0,1){20}}
\put(220,60){\line(2,5){6}}
\put(250,60){\line(-2,5){6}}
\put(220,40){\line(2,-5){6}}
\put(250,40){\line(-2,-5){6}}
\put(250,99){\line(-1,-1){24}}
\put(220,99){\line(1,-1){24}}
\put(250,1){\line(-1,1){24}}
\put(220,1){\line(1,1){24}}

\put(270,50){$=$}

\put(290,4){\line(0,1){92}}
\put(320,4){\line(0,1){92}}


\end{picture} 
\eeq 

The assignments that result from this for the left and right ribbon $2\pi$-twists are 
summarized in Figure~\ref{fig-Rigg-tgl}. Note, that in the assignment on the 
right hand side the full circle on the left side stands for a left handed twist
for the framing, while the fat dot on the right hand side indicates a decoration of the strand
by the element $v^{-1}$.

\begin{figure}[ht]
\begin{center}
\leavevmode
\epsfxsize=13cm
\epsfbox{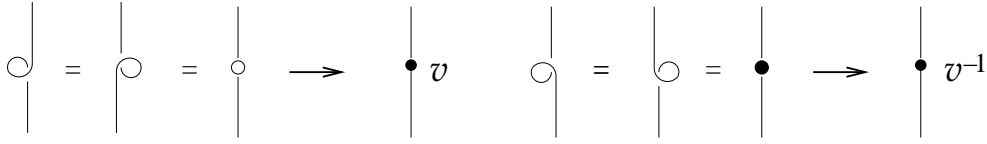}
\end{center}
\caption{Twist Assignments}\lbl{fig-Rigg-tgl}
\end{figure}

It is clear that after application of these types of manipulations to any decorated diagram 
we eventually obtain a set of disjoint, planar curves  which can be one of four types.
For each of these types we describe next the evaluation rule that leads to the definition
of a linear map ${\cal V}^{\#}(T_M)$:

Components of the first type are closed circles decorated with one element $a_i\in{\cal A}\,$ on the
right side. To this we associate the number $\mu(a_i)\in\kk\,$. 

Next,  we may have an arc at the bottom line of the diagram connecting points
$p_k'$ and $q_k'$ with one decoration 
$\,b_k\in{\cal A}\,$ at the left strand. To this to we associate 
the vector $\,b_k\in{\cal A}^{(k)}\,$ in the $k$-th copy of 
the tensor product $\,{\cal A}^{\otimes g_2}\,$.

Thirdly, for an arc at the top line between points $p_j$ and $q_j$ with decoration
$c_j\in{\cal A}\,$ on the right we assign the linear form $l_{c_j}:{\cal A}^{(j)}\to\kk\,$
given by $l_{c_j}(x)=\mu(S(x)c_j)$
 on the $j$-th copy of 
the tensor product $\,{\cal A}^{\otimes g_1}\,$. 

Finally, we may have pairs of straight strands that connect a pair $\,\{p_j,q_j\}\,$ to the pair
$\,\{p_k',q_k'\}\,$,  carrying decorations, $a$ and $b$. In case the strands are parallel,
that is, one connects $p_j$ to $p'_k$ and the other $q_j$ to $q'_k$, we assign a linear
map $T_{a,b}\,:{\cal A}^{(j)}\to {\cal A}^{(k)}\,$ between the $j$-th copy of 
$\,{\cal A}^{\otimes g_1}\,$ to the $k$-th copy of $\,{\cal A}^{\otimes g_2}\,$,
by $T_{a,b}(x)=axS(b)\,$. 

If the connecting strands cross over we apply in addition the endomorphism $K(x)=G^{-1}S(x)$ on the
$k$-th copy ${\cal A}^{(k)}$ for a crossing right at the bottom line. It is quite useful to summarize 
these rules also pictorially as follows: 
\beq\lbl{pic-Rclos}
\setlength{\unitlength}{0.25mm}
\begin{picture}(250,70)(0,-5)
\put(35,35){\circle{50}}
\put(60,35){\circle*{8}}
\put(65,37){$a_i$}

\put(130,37){\vector(1,0){20}}

\put(200,35)
{$\;\mu(a_i)\;$}
\end{picture}
\eeq
\beq\lbl{pic-Rbot}
\setlength{\unitlength}{0.25mm}
\begin{picture}(400,80)(0,-15)
\put(5,10){\rule{15mm}{1mm}}
\put(12,-10){$p_k'$}
\put(52,-10){$q_k'$}
\put(35,35){\oval(40,40)[t]}
\put(15,35){\line(0,-1){25}}
\put(55,35){\line(0,-1){25}}
\put(15,24){\circle*{10}}
\put(21,25){$b_k$}

\put(130,30){\vector(1,0){20}}

\put(200,25)
{$b\,:\kk\longrightarrow {\cal A}^{(k)} \,\,\,:\;\,\,\;1\,\mapsto\,\,b_k\,$}
\end{picture}
\eeq
\beq\lbl{pic-Rtop}
\setlength{\unitlength}{0.25mm}
\begin{picture}(500,93)(0,-5)
\put(5,60){\rule{15mm}{1mm}}
\put(12,80){$p_j$}
\put(52,80){$q_j$}
\put(35,35){\oval(40,40)[b]}
\put(15,35){\line(0,1){25}}
\put(55,35){\line(0,1){25}}
\put(55,42){\circle*{10}}
\put(62,43){$c_j$}

\put(130,53){\vector(1,0){20}}

\put(200,40)
{$l_{c_j}\,:{\cal A}^{(j)}\longrightarrow \kk\,\,\,:\;\,\,\;x\,\mapsto\,\,\mu(S(x)c_j)\,$}
\end{picture}
\eeq
\beq\lbl{pic-Rthru}
\setlength{\unitlength}{0.25mm}
\begin{picture}(500,112)(0,-15)
\put(20,5){\line(0,1){60}}
\put(20,35){\circle*{10}}
\put(27,36){$a$}

\put(12,80){$p_j$}
\put(52,80){$q_j$}

\put(12,-10){$p_k'$}
\put(52,-10){$q_k'$}
\put(60,5){\line(0,1){60}}

\put(60,35){\circle*{10}}
\put(67,36){$b$}

\put(130,33){\vector(1,0){20}}

\put(200,27){$T_{a,b}\,:{\cal A}^{(j)}\longrightarrow {\cal A}^{(k)}\,\,\,:\;\,\,\;x\,\mapsto\,\,
a x S(b)\,$}
\end{picture}
\eeq

\setlength{\unitlength}{0.25mm}
\begin{picture}(400,80)(0,-15)
\put(5,10){\rule{15mm}{1mm}}
\put(12,-10){$p_k'$}
\put(52,-10){$q_k'$}
\put(15,48){\line(1,-1){37}}
\put(55,48){\line(-1,-1){37}}

\put(130,30){\vector(1,0){20}}

\put(200,25){$K\,:{\cal A}^{(k)}\longrightarrow {\cal A}^{(k)}
\,\,\,:\;\,\,\;x\,\mapsto\,\, G^{-1}S(x)\,$}

\end{picture}

From these rules for evaluating diagrams 
 we obtain a linear map $\,{\cal A}^{\otimes g_1}\to{\cal A}^{\otimes g_2}\,$
for any decorated planar tangle. For a given tangle $T_M$ we denote by $\,{\cal V}^{\#}(T_M)\,$ the 
sum of all of these maps associated to the sum of decorated diagrams for $T_M$. 
Thus,  if we consider, for simplicity, a tangle $T_M$ without
components of the fourth type, and denote by $a^{\nu}_i\,$,
$b^{\nu}_j\,$ and $c^{\nu}_k\,$ the respective elements of the $\nu$-th summand of 
the same untangled curve of $T_M$, this linear map can be expressed as 
$$
{\cal V}^{\#}(T_M)\;\;:=\;\;\sum_{\nu\, }\mu(a_1^{\nu})\ldots \mu(a_N^{\nu})\, 
b_1^{\nu}\otimes\ldots\otimes b_{g_2}^{\nu}\, l_{a_1^{\nu}}\otimes\ldots\otimes 
l_{a_{g_1}^{\nu}}\;\;. 
$$
For tangles with strand pairs that connect top and bottom pairs we insert the operators
$T_{a,b}\,$ in the respective positions.
\medskip

\begin{lemma}\lbl{lm-Vprop}
The linear maps ${\cal V}^{\#}(T_M)$ are well defined, (covariantly)
functorial under the composition of tangles, and they commute with   
 the adjoint action of ${\cal A}$ on ${\cal A}^{\otimes g}$. They are also invariant under isotopies
and the following moves: 
\begin{enumerate}
\item 2-handle slides of any type of strand over a closed component of  $\,T_M\,$
\item Adding/removing an isolated Hopf link for which one component has 0-framing 
and the other framing 0 or 1. 
\end{enumerate}
\end{lemma}

{\em Proof:} The fact that the construction procedure for a given diagram is unambiguous
is almost straight forward, except that one has to pay attention to the positioning of the 
resulting elements. Details for closed links can be found in \cite{Kau94}. 
 Functoriality is easily checked from the rules of construction. The fact that the maps
are ${\cal A}$-equivariant follows from the fact that it is a special
case of the categorical construction in \cite{KerLub00}, and the fact that 
$f:{\cal A}\to{\cal A}^*$ intertwines the adjoint with the  coadjoint action. 
 Invariance under isotopies follows, as in \cite{Hen96} or \cite{KauRad95}, from the
properties of the $R$-matrix of a quasitriangular Hopf algebra. In the same articles
the 2-handle slide is directly related to the defining equation (\ref{eq-int-def})
of the right integral, see also \cite{Lub95} for the categorical version of the 
argument. Invariance under the Hopf link moves is a direct consequence of the normalizations
in (\ref{eq-balnorm}), since they imply that  the Hennings invariants on the Hopf links are all one. 
\ep
\medskip

In order to describe the reduction procedure that allows us to define a TQFT also
for non-modular Hopf algebras 
we introduce the operators associated to the diagrams in Figure~\ref{fig-Smat},
the left being isotopic to the one in Figure~\ref{fig-S-tgl}.
\begin{figure}[ht]
\begin{center}
\leavevmode
\epsfysize=2.2cm
\epsfbox{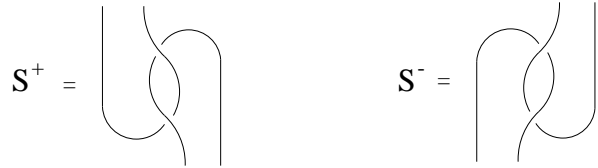}
\caption{$S^{\pm}$-Transformations}\lbl{fig-Smat}
\end{center} 
\end{figure}
The double crossing is replaced by the elements $m^{+}_j, n^{+}_j$ from 
${\cal M}\,=\,\sum_jm_j^+\otimes n_j^+\,$, as defined in (\ref{eq-M-matrix}). 
The transformation $S^+:{\cal A}\to{\cal A}$ is readily worked out to be
\begin{equation}\lbl{eq-S-matrix}
S^+(x)\;=\;\sum_j\mu(S(x) m_j^+) n_j^+\qquad.
\end{equation}
The formula for $S^-$ follows analogously, substituting ${\cal M}$ for
${\cal M}^{-1}=\sum_jm_j^-\otimes n_j^-$. We consider next the result $\Pi$ 
of stacking the two tangles in Figure~\ref{fig-Smat} on top of each other:
\begin{lemma}\lbl{lm-Pi}
Let  $\Pi:=S^+\circ S^-=S^-\circ S^+$, and denote
$\Pi^{(j)}=1\otimes\ldots1\otimes \Pi\otimes 1\ldots\otimes 1$, with $\Pi$
occurring in the $j$-th tensor position.  
\begin{enumerate}
\item $\Pi$ is an idempotent that commutes with the adjoint action of $\cal A$.
\item ${\cal V}^{\#}(T_M)\Pi^{(j)}={\cal V}^{\#}(T_M)$ if the $j$-th top index pair in 
$T_M$ is attached to a top ribbon in $T_M$. (Analogously for bottom ribbons).
\item $\Pi^{(k)}{\cal V}^{\#}(T_M)={\cal V}^{\#}(T_M)\Pi^{(j)}$ if $T_M$ has 
a through pair
connecting the $j$-th top pair to the $k$-th bottom pair. 
\end{enumerate}
\end{lemma}

{\em Proof:} For {\em 1.} note that the picture for $\Pi$ consists of two arcs
that are connected by a circle. Stacking $\Pi$ on top of itself we obtain the picture 
for $\Pi^2$ by functoriality in 
Lemma~\ref{lm-Vprop}. The resulting tangle is the chain of circles $C_j$ and 
arcs $A_{t/b}$ depicted on the left of Figure~\ref{fig-pi-idem}. By  
{\em 1.} of Lemma~\ref{lm-Vprop} we may use 2-handle slides to manipulate
this picture. We first slide $C_1$ over $C_3$, and then $A_b$ over $C_2$.
The result is the tangle for $\Pi$ and a separate Hopf link. The value of
the latter, however, is 1 by (\ref{eq-balnorm}). Hence, $\Pi^2=\Pi$. 
\begin{figure}[ht]
\begin{center}
\leavevmode
\epsfysize=2.2cm
\epsfbox{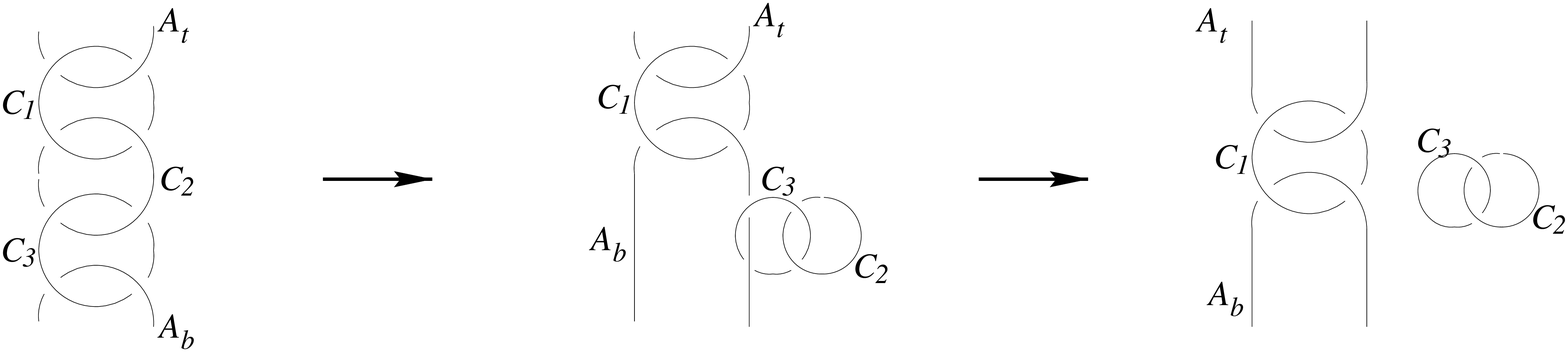}
\caption{$\Pi$ is idempotent}\lbl{fig-pi-idem}
\end{center} 
\end{figure}

Equivariance with respect to the action of ${\cal A}$ is immediate from 
Lemma~\ref{lm-Vprop}. 

For {\em 2.} we repeat an argument from \cite{KerLub00}. Suppose $\tau$ is a top
component and $\eta$ any band connecting two intervals $I_i$ in $\tau$ in an
orientation preserving way. To this we associated the surgered diagram in which
the component $\tau$ is replaced by the union $\tau_{\eta}$ of three components.
They are obtained by cutting away the intervals $I_i$ from $\tau$ and inserting
the other two edges of $\eta$ at the endpoints  $\partial I_i$ as indicated in
Figure~\ref{fig-toppath}. Furthermore,  we insert a 0-framed annulus $A$ around $\eta$. 
\begin{figure}[ht]
\begin{center}
\leavevmode
\epsfysize=2.2cm
\epsfbox{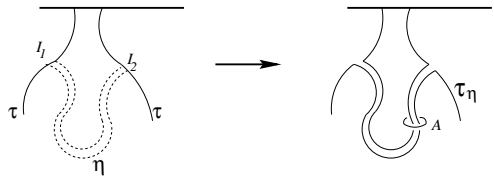}
\caption{$\eta$-Surgery}\lbl{fig-toppath}
\end{center} 
\end{figure}
Sliding any other component over $A$ at an arbitrary point along $\eta$
has the effect of just moving it through $\eta$ at this point. Moreover, we
can slide a $\pm 1$-framed annulus $K$ over $A$ so that it surround the two 
parallel strands in $\tau_{\eta}\,$, and then slide the two strands over $K$. The 
effect is the same as putting a $2\pi$-twist into $\eta$. These two operation
allow us to move any band $\eta$ to any other band $\eta'$ such that $\tau_{\eta}$
and $\tau_{\eta'}$ are related by a sequence of two handle slides. 

Now,  adding the picture of $\Pi$ to the top-component $\tau$ of a  tangle $T_M$ 
is the same as surgering $\tau$ along a straight band parallel and close to
the interval between the attaching points of $\tau$ at the top line. We replace
this $\eta$ by a small planar arc at $\tau$ separate from the rest of the 
tangle. Surgery along this corresponds to linking  a Hopf link to $\tau$, 
as $C_2\cup C_3$ is linked to $A_b$ in the middle of Figure~\ref{fig-pi-idem}, 
and consequently  can be removed by the same argument. 
The proofs for the formulas for bottom and through strands are entirely analogous. 
\ep 

Set $\Pi^{\#}=\Pi^{\otimes g}$, when acting on ${\cal A}^{\otimes g}$. It follows
now easily from Lemma~\ref{lm-Pi} that 
${\cal V}^{\#}(T_M)\Pi^{\#}=\Pi^{\#}{\cal V}^{\#}(T_M)$ for all $T_M$. Thus 
each ${\cal V}^{\#}(T_M)$ maps the image of $\Pi^{\#}$ to itself so that we can 
define the restriction 
\begin{equation}\lbl{eq-redtqft}
{\cal V}(T_M):={\cal V}^{\#}(T_M)\Bigl |_{im(\Pi^{\#})}:\;\;
{\cal V}_{\cal A}(\Sigma_{g_1,1})\;\longrightarrow\;{\cal V}_{\cal A}(\Sigma_{g_2,1})\;\;,
\end{equation}
where the vector spaces are given as 
\beq\lbl{eq-A0def}
{\cal V}_{\cal A}(\Sigma_{g,1})=\Pi^{\#}({\cal V}^{\#}(\Sigma_{g}))={\cal A}_0^{\otimes g}
\qquad\mbox{with} \qquad {\cal A}_0=\Pi({\cal A})\;\;\;.
\eeq
\begin{theorem}\lbl{thm-henn}
The assignment $\cal V$ as given in (\ref{eq-redtqft}) yields a 
well defined, 2-framed, relative, ${\cal A}-equivariant$ 
topological quantum field theory
$$
{\cal V}_{\cal A}\;:\quad \Cob^{2fr,\bullet}\;\longrightarrow\;{\cal A}{\rm -mod}_{\kk}\quad\subset
\quad Vect(\kk)\quad. 
$$
Using the invariance functor ${\rm Inv}={\rm Hom}(1,\_):{\cal A}{\rm -mod}\to
Vect(\kk)$ we  obtain an ordinary  2-framed TQFT for closed surfaces as
$$
{\cal V}_{\cal A}^0:={\rm Inv}\circ{\cal V}_{\cal A}\;:
\quad \Cob^{2fr}\;\longrightarrow\; 
\quad Vect(\kk)\quad. 
$$
\end{theorem}

{\em Proof:}
We recall from \cite[Proposition 12]{Ker99} that two presentations, $T_M$ and $T'_M$,
of a framed, relative cobordism $M\in \Cob^{2fr,\bullet}$ are related by the 
moves described in Lemma~\ref{lm-Vprop} and the so called $\sigma$-moves, which
consist of adding the picture of $\Pi$ to a pair of points at the top or bottom
line of the diagram. From ${\cal V}(T_M)\Pi^{(j)}={\cal V}^{\#}(T_M)\Pi^{\#}\Pi^{(j)}
={\cal V}^{\#}(T_M)\Pi^{\#}$ we see that ${\cal V}(T_M)$ is invariant under this
move. Hence, ${\cal V}(T_M)$ only depends on the cobordism represented by $T_M$
and we can write ${\cal V}_{\cal A}(M):={\cal V}(T_M)$.  

Due to the equivariance of  $\Pi$ also ${\cal A}_0$ from (\ref{eq-A0def}) is invariant 
under the adjoint action of ${\cal A}$, and the restricted maps commute with
the action of ${\cal A}$ as well. Functoriality of ${\cal V}$ follows from
functoriality of ${\cal V}^{\#}$ and the fact that $\Pi^{\#}$ commutes with 
${\cal V}^{\#}$.

Since each ${\cal V}(M)$ commutes with the action of ${\cal A}$ they also
map the ${\cal A}$-invariant subspaces 
${\cal V}^0(\Sigma_{g}):={\rm Inv}({\cal V}(\Sigma_{g,1}))$ to themselves. 
This implements the additional $\tau$-move \cite{Ker99} needed to represent
cobordisms between closed surfaces. 
\ep

\head{6. The Algebra $\cal N$}\lbl{S6}
The Hopf algebra $\cal N$ we will define in this section is  the same 
as the algebra  $A_2$  described by Radford in Example~1 of Section~4.1 
in \cite{Rad76}.  The quasitriangular structure that we endow  $\cal N$ 
with is essentially distilled from the one of $U_{-1}({\mathfrak s}{\mathfrak l}_2)$.

Let $\E\cong\R^2$ be the Euclidean plane, and consider 
the 8-dimensional algebra 
\beq\lbl{eq-Adef}
{\cal N}\;:=\;\Z/2\ltimes \ext * \E\;.
\eeq
The generator of $\Z/2$ is denoted by $K$, with $K^2=1$, and we write $x^K=KxK$
for any $x\in{\cal N}$. We thus have relations $w'w=-ww'$ and $w^K:=KwK=-w$ 
for all $w,w'\in\E$.

\blm\lbl{lm-Ahopf}
$\cal N$ is a Hopf algebra with coproducts 
\beq\lbl{eq-coprod}
\Delta(K)=K\otimes K\quad\mbox{ and}\qquad \Delta(w)\,=\,w\otimes 1\,+\,K\otimes w\;\;
\forall w\in\E\;.
\eeq
\elm 
{\em Proof:} The fact that $\Delta:{\cal N}\to{\cal N}^{\otimes 2}$ is a 
coassociative homomorphism
is readily verified. The antipode is given by 
\beq\lbl{eq-antipode}
S(K)\,=\,K\qquad\mbox{and}\qquad S(w)=-Kw,\;\;\forall w\in\E\;. 
\eeq
\ep

We note the following formulas for the adjoint action and antipode: 
\beq\lbl{eq-somefor}
ad(w)(x)=wx-x^Kw\qquad, \qquad S^2(x)=x^K\qquad\forall x\in{\cal N}, w\in\E
\eeq
 
Let us pick a non-zero element $\rho\in\ext 2\E\subset{\cal N}$, and for this define
a form $\mu_0\in{\cal N}^*$ as follows:
\beq\lbl{eq-Ainteg}
\mu_0(\rho)\;=\;1\;\,,\qquad\quad \mu_0(K\rho)=0\;,\qquad\mbox{and}
\eeq
$$
\mu_0(K^{\delta}x)=0\;,\qquad\forall x\in\ext j\E\;\;,\;\;\mbox{whenever}
\;\;j,\delta\in\{0,1\}\;\;.
$$
\blm\lbl{lm-Ainteg}\ \ \ 
$\mu_0$ is a right (and left) integral on ${\cal N}$. Moreover,  
\beq\lbl{eq-intIN}
\lambda_0:=  (1+K)\rho\qquad\mbox{with}\quad\mu_0(\lambda_0)=1
\eeq 
is a two sided integral in ${\cal N}$.
\elm
 
{\em Proof:} Straightforward verification of (\ref{eq-int-def}). The defining
equation for a two sided integral in ${\cal N}$ is 
$x\lambda_0=\lambda_0 x=\epsilon(x)\lambda_0$, which is also readily found.
\ep

Next, we fix a basis $\{\th,\tb\}$ for $\E$.
We define an $R$-matrix, ${\cal R}\in{\cal N}\otimes{\cal N}$, by the formula 
\beq\lbl{eq-ARmat}
{\cal R}\quad:=\quad \Bigl(1\otimes 1\,+\,\th\otimes K\tb\Bigr)\cdot {\cal Z}\;, 
\qquad\mbox{where}\quad {\cal Z}:=\frac 12 \sum_{i,j=0}^1(-1)^{ij}K^i\otimes K^j
\eeq
 
\blm\lbl{lm-Rmat}
The element ${\cal R}$ makes ${\cal N}$ into a quasitriangular Hopf algebra.

Moreover, ${\cal N}$ is a ribbon Hopf algebra with unique balancing element $G=K$. 

\elm
{\em Proof:} Quasitriangularity follows from a  
straightforward verification of the axioms in \cite{Dri90}.
We compute the special element  $u^{-1}=\sum_jf_jS^2(e_j)=K(1+\tb\th)$ for which
$uS(u)^{-1}=uu^{-1}=1$ so that $G=K$ is a valid and unique choice. The ribbon element
is then given by 
\beq\lbl{lm-Aribb}
v\;:=\;1+\rho\qquad\mbox{with}\qquad \rho:=\tb\th
\eeq
\ep


For the monodromy matrix, as defined in (\ref{eq-M-matrix}),  we obtain:
\beq\lbl{eq-Amono}
{\cal M} \;= \; 1\,+\,K\tb\otimes \th \,+\,\th K\otimes \tb-
\rho\otimes\rho\;\;.
\eeq
Setting $T=K\tb\otimes \th \,+\,\th K\otimes \tb$ we compute $T^2=-2\rho\otimes\rho\,$ and 
$T^3=0\,$ 
so that ${\cal M}=\exp(T)$. Hence we can also compute $p$-th powers of the monodromy matrix:
\beq\lbl{eq-qAmono}
{\cal M}^p\;=\;\exp(pT)\;=\;1\,+\,pT\,+\,\frac {p^2}2 T^2\;.
\eeq

With $\mu_0$ as defined in (\ref{eq-Ainteg}), and  for $\rho$ as in (\ref{lm-Aribb})
we find $\mu_0\otimes\mu_0({\cal M})=
\mu_0(v)\mu_0(v^{-1})=-1$. Hence, in order to fulfill (\ref{eq-balnorm}) we need to use the 
renormalized integrals
\beq\lbl{eq-renormint}
\mu\;=\;i\mu_0\;\;, \qquad\;\;\lambda=\frac 1 i\lambda_0\;\;,
\qquad\;\; \mbox{with}\;\;\;i=\sqrt{-1}\;\;.
\eeq
For these choices we compute the $S^{\pm}$-transformations assigned to (\ref{eq-S-matrix}) as follows:
\beq\lbl{eq-ASmat}
\begin{array}{lll}
\frac 1 iS^{\pm}(w)=\mp w\;\;\;\forall w\in\E &\qquad \frac 1 i S^{\pm}(\rho)=1\\
\frac 1 iS^{\pm}(Kx)=\,0\;\;\;\forall x\in\ext * \E &\qquad \frac 1 i S^{\pm}(1)=-\rho\;.\\
\end{array}
\eeq
This implies that the projector $\Pi$ from Lemma~\ref{lm-Pi} has   
kernel $ker(\Pi)=\{Kw\, :\, w\in\ext *\E\}$ and image
\beq\lbl{eq-APi}
{\cal N}_0\;\;=\;\;im(\Pi)\;\;\;=\;\;\ext * \E\;.
\eeq
From (\ref{eq-somefor}) we see that ${\cal N}_0$ acts trivially
on itself  so that the action of ${\cal N}$ factors through the obvious 
$\Z/2\Z={\cal N}/{\cal N}_0$-action.

Finally, we note that $\SL(2,\R)$ acts on  $\E$ and, hence, also on $\cal N$,
assuming $K$ is $\SL(2,\R)$-invariant.  
\begin{lemma}\lbl{lm-ASLinv}  
$\SL(2,\R)$ acts on $\cal N$ 
 by Hopf algebra automorphisms. 

The ribbon element $v$, the monodromy $\cal M$, and the two integrals
 are invariant under this action. 
\end{lemma}

{\em Proof:} The fact that  $\SL(2,\R)$ yields algebra automorphsims is 
obvious by construction. Linearity of coproduct and antipode in $w$ 
in (\ref{eq-coprod}) and (\ref{eq-antipode}) imply that this is, in fact, 
a Hopf algebra homomorphism. $v$ and $\lambda$ are invariant
since $\SL(2,\R)$ acts trivially on $\E\wedge\E\,$. Invariance of $\cal M$
follows then from (\ref{eq-M-matrix}). 
\ep

Note,  that $\cal R$ itself is {\em not} $\SL(2,\R)$-invariant.

\bigskip

\head{7. The Hennings TQFT for $\cal N$}\lbl{S7}

From (\ref{eq-APi}) and (\ref{eq-redtqft}) we see that the vector spaces of the 
Hennings TQFT for the algebra from  (\ref{eq-Adef}) are given as
\beq\lbl{eq-AVvect}
{\cal V}_{\cal N}(\Sigma_g)\;\;:=\;\;\Bigl(\ext *\E\Bigr)^{\otimes g}
\qquad\mbox{with}\qquad {\rm dim}({\cal V}_{\cal N}(\Sigma_g))=4^g\;.
\eeq 
We now compute the action of the mapping class group generators from the 
tangles in Figures~\ref{fig-A-tgl}, \ref{fig-D-tgl}, and \ref{fig-S-tgl}. 
From the extended Hennings rules it is clear that
 the pictures for both $A_j$ and $S_j$ result in actions only on the  $j$-th
factor in the tensor product in (\ref{eq-AVvect}). For $A_j$ we use the presentation
from Figure~\ref{fig-A-tgl} and the rules from Figure~\ref{fig-Rigg-tgl} and
(\ref{pic-Rthru}) to obtain the linear map $\Aa(x):=x\cdot v$.

The extra 1-framed circle in Figure~\ref{fig-S-tgl} results in an extra factor
$\mu(v)=i$, since an empty circle corresponds to an insertion of $v$. The action on 
the $j$-th factor is thus given by an application of $\Ss:=iS^{+}\bigl |_{{\cal N}_0}$
so that
\beq\lbl{eq-SSform}
\Ss(\rho)=-1\;,\qquad\quad \Ss(1)=\rho\;,\qquad{\rm and}\qquad \Ss(w)=w\;,\quad\forall w\in\E\;\;.
\eeq 

Similarly,  $D_j$ acts only on the $j$-th and the $(j+1)$-st
factors of ${\cal N}_0^{\otimes g}$. From  (\ref{pic-Rthru}) and the formula 
for ${\cal M}^{-1}$ in (\ref{eq-Amono}) we compute for the action on these two factors
\beq\lbl{eq-D-2act}
\D\;:\;{\cal N}_0^{\otimes 2}\,\to\,{\cal N}_0^{\otimes 2}, \qquad
x\otimes y \;\mapsto\; x\otimes y +  x\th\otimes \tb y -
x\tb\otimes \th y -  x\rho \otimes \rho y\;. 
\eeq
The generators of the mapping class group $\Gamma_g$ are thus represented as follows:
\beq\lbl{eq-AVgen}
\begin{array}{ll}
{\cal V}_{\cal N}({\bf I}_{A_j})\;\;=\;\;I^{\otimes j-1}\otimes \Aa\otimes I^{\otimes g-j}\;,
 &\qquad \quad {\cal V}_{\cal N}({\bf I}_{S_j})\;\;=\;\;I^{\otimes j-1}\otimes \Ss\otimes I^{\otimes g-j}\\
&\\
\qquad\qquad\mbox{and}&
{\cal V}_{\cal N}({\bf I}_{D_j})\;\;=\;\;I^{\otimes j-1}\otimes \D\otimes I^{\otimes g-j-1}\;.\\
\end{array}
\eeq
Let us also compute the linear  maps associated to the cobordisms ${\bf H}_g^{\pm}$
from (\ref{eq-hdlcob}). Their tangle presentations follow from \cite{Ker99} and 
have the forms given in Figure~\ref{fig-tglhdl}. 

\begin{figure}[ht]
\begin{center}
\leavevmode
\epsfysize=3.2cm
\epsfbox{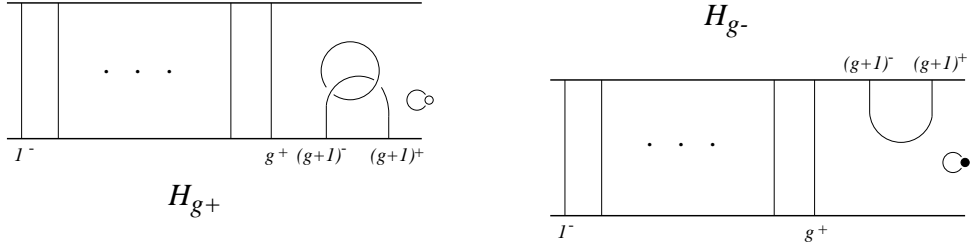}
\caption{Tangles for Handle additions}\lbl{fig-tglhdl}
\end{center} 
\end{figure}

We included $\pm 1$-framed circles to adjust the 2-framings of ${\bf H}_g^{\pm}$. 
A 0-framed circle around a
strand has the effect of inserting $\lambda= S^+(1)=\frac 1 i\rho$. In this
normalization we find with $\rho=i\Pi\lambda$ and (\ref{pic-Rbot}) that 
\beq\lbl{eq-Ahdl+map}
{\cal V}_{\cal N}({\bf H}_g^{+})\;:\;\;\alpha\;\mapsto\;\;\alpha\otimes \rho
\qquad \quad \forall \alpha \in{\cal N}_0^{\otimes g}\;.
\eeq
Similarly, we obtain from (\ref{pic-Rtop}) that 
\beq\lbl{eq-Ahdl-map}
{\cal V}_{\cal N}({\bf H}_g^{-})\;:\;\;\alpha\otimes x\;\mapsto\;\;\mu_0(x)\alpha
\qquad \quad \forall \alpha \in{\cal N}_0^{\otimes g}, \; x\in {\cal N}_0\;,
\eeq
where $\mu_0$ is as in (\ref{eq-Ainteg}). We note the following:

\begin{lemma}\lbl{lm-ASL2inv}
The generators in (\ref{eq-AVgen}), (\ref{eq-Ahdl+map}), and (\ref{eq-Ahdl-map})
intertwine the $\SL(2,\R)$-action on ${\cal N}_0^{\otimes g}$.
\end{lemma}

{\em Proof:} The fact that $\Aa$ and $\D$ commute with the $\SL(2,\R)$-action follows from
invariance of $v$ and $\cal M$. From (\ref{eq-ASmat}) we see that $\Ss$ is scalar on the 
non-invariant part, and thus commutes as well. Finally, $\rho$ and $\mu_0$ are clearly invariant.

\ep

For $g\geq 0$ set $\chi_g:=S_g\circ\ldots\circ S_1$, $h^+_g:= {\bf H}_{g-1}^+\circ \ldots 
\circ {\bf H}_{0}^+$, and $h^-_g:={\bf H}_{0}^-\circ \ldots 
\circ {\bf H}_{g-1}^-$. We define a standard closure of a 2-framed 3-cobordism as the closed
3-manifold 
\beq\lbl{eq-closM}
\lz M\rz\;\;:=\;\;h^-_{g_2}\circ\chi_{g_2}\circ M\circ \chi^{-1}_{g_1}\circ h^+_{g_1}\;\cup\,D^3\;.
\eeq
If $M$ is represented by a tangle $T$ we obtain, similarly, a link $\lz T\rz$. We introduce 
the following function from the class of 2-framed cobordisms into $\Z/2$:
\beq\lbl{eq-number}
\varphi(M)\quad:=\quad \beta_1(\lz M\rz)\;+\;{\rm sign}(\lz T\rz)\;\;{\rm mod}\; 2\;,
\eeq
where $\beta_j$ denotes the $j$-th Betti number. 
We further denote by $\Cob^{22fr,*}\subset \Cob^{2fr,*}$ the subset of  all cobordisms $M$
with $\varphi(M)=0$, which we will call {\em evenly} 2-framed. 

\begin{lemma}\lbl{lm-mod2prop}
\begin{enumerate}
\item $\varphi(M)\;=\;|\lz T\rz|\;{\rm mod}\, 2$, where $|\lz T\rz|\,:=$ $\#$ components of $\lz T\rz$.
\item  $\varphi(M)\;=\;\#$ components of $ T$ not connected to the bottom line. 
\item ${\cal V}_{\cal N}(M)$ is real  if $\varphi(M)=0$ and imaginary for $\varphi(M)=1$.  
\item $\Cob^{22fr,*}$ is a subcategory.
\end{enumerate}
\end{lemma}

{\em Proof:} Let $W$ be the 4-manifold given by adding 2-handles to $D^4$ along 
$\lz T\rz\subset S^3$ so that $\lz M\rz=\partial W$, and let $L_T$ be the linking matrix of 
$\lz T\rz$. We have $\beta_2(W)= |\lz T\rz |=d_++d_-+d_0$, where $d_+$, $d_-$, and $d_0$ are
the number of eigenvalues of $L_T$ that are $>0$, $<0$, and $=0$ respectively. 
From the exact sequence 
$0\to H_2(\lz M\rz )\to H_2(W) \stackrel {L_T}{\longrightarrow} H^2(W)\to H_1(\lz M\rz)\to 0$
we find that $\beta_1(\lz M\rz)=d_0$, which implies {\em 1.} using ${\rm sign}(W)=d_+-d_-$.
{\em 2.} follows immediately from the respective tangle compositions. 

The possible components
not connected to the bottom line are strands connecting point pairs at  the top line
or closed components. From the rules (\ref{pic-Rclos}) through (\ref{pic-Rthru}) we see 
that these are just the types of components that involve an evaluation against $\mu=i\mu_0$.
All other parts of the Hennings procedure involve only real maps. Finally, {\em 4.} follows
from counting tangle components under composition. 
\ep

\begin{propos}\lbl{pp-ATQFT}
The Hennings procedure yields a relative, 2-framed, $\SL(2,\R)$-equivariant, half-projective TQFT
$$
{\cal V}_{\cal N}\;:\;\;\Cob^{2fr,\bullet}\;\longrightarrow\;\SL(2,\R)-mod_{\Cc}\;,
$$
which is $\Z/4$-projective on $\Cob^{\bullet}$. We have a restriction
$$
{\cal V}_{\cal N}^{(2)}\;:\;\;\Cob^{22fr,\bullet}\;\longrightarrow\;\SL(2,\R)-mod_{\R}\;,
$$
which is $\Z/2$-projective on $\Cob^{\bullet}$.
\end{propos}

{\em Proof:} From Lemma~\ref{lm-ASL2inv} we know that the generators of $\Gamma_g$ are
represented $\SL(2,\R)$-equivariantly, hence also $\Gamma_g$ itself. The decomposition in 
(\ref{eq-Heegaard}) and equivariace of the maps in (\ref{eq-Ahdl+map}) and  
(\ref{eq-Ahdl-map}) implies the same for  general cobordisms. That this TQFT is
half-projective follows from the fact that $\cal N$ is non-semisimple, or, equivalently,
that ${\cal V}_{\cal N}(S^1\times S^2)=\mu(1)=\varepsilon(\lambda)=0$, see \cite{Ker98}.
The projective phase of the TQFT is determined by the value $\mu(v)=i$ on the 1-framed
circle.

Lemma~\ref{lm-mod2prop}, {\em 3.} implies that ${\cal V}_{\cal N}^{(2)}$ maps into
the {\em real}  $\SL(2,\R)$-equivariant maps and modules. This reduces the 
ambiguity of multiplication with $i$ to a  sign ambiguity. 
\ep 

\medskip

 An important point of view in the TQFT constructions in \cite{KerLub00} 
is the existence of a categorical Hopf algebra, which can be understood as
the TQFT image of a topological Hopf algebra given as an object in $\Cob^{\bullet}$. 

To be more precise, in 
\cite{Yet97} and \cite{Ker96} $\Cob^{\bullet}$ is described as a braided tensor
category, and it is found 
that the object $\Sigma_{1,1}\in \Cob^{\bullet}$ is naturally identified as a
braided Hopf algebra in this category in the sense of \cite{Maj93} and \cite{LubMaj94}. 
Particularly, $\Sigma_{2,1}$  is identified with   $\Sigma_{1,1}\otimes\Sigma_{1,1}$
since the tensor product on $\Cob^{\bullet}$  is defined by sewing two 
surfaces together along a pair of pants. The 
multiplication and comultiplication are thus given by 
elementary cobordisms ${\bf M}:\Sigma_{2,1}\to\Sigma_{1,1}$ and 
${\bf \Delta} :\Sigma_{1,1}\to \Sigma_{2,1}$. Their tangle diagrams are worked out
explicitly in \cite{BKLT00}, and depicted in Figure~\ref{fig-MC-tgl}
with minor modifications in the conventions:

\begin{figure}[ht]
\begin{center}
\leavevmode
\epsfysize=2.2cm
\epsfbox{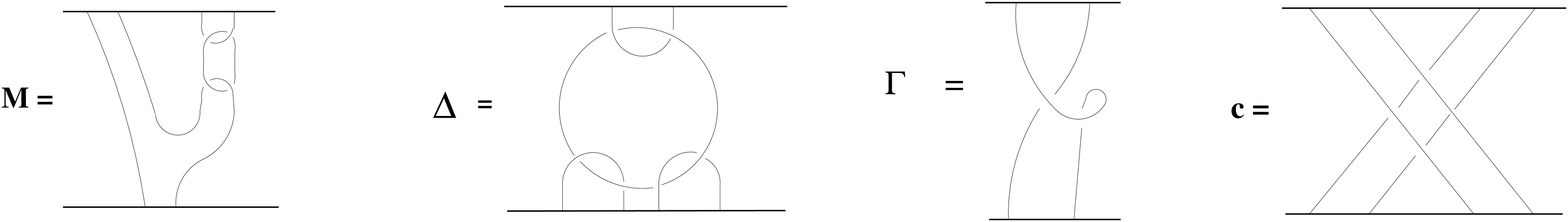}
\caption{Tangles for Mulitplications}\lbl{fig-MC-tgl}
\end{center} 
\end{figure}

Here ${\bf c}:\Sigma_{2,1}\to\Sigma_{2,1}$ is the braid isomorphism. 
The braided antipode is given by the tangle 
$\Gamma=(S^+)^2$, with $S^+$ as in Figure~\ref{fig-Smat}. 

\begin{lemma}
\lbl{lm-multHeeg}
The cobordisms {\bf M} and ${\bf \Delta}$ have the following 
Heegaard decompositions.
$$
{\bf M}\;=\;{\bf H}_2^-\circ {\bf I}_{D_1 \circ S_2} \qquad\mbox{and}
\qquad {\bf \Delta}\;=\;
 {\bf I}_{S_1\circ  D_1^{-1} \circ S_1^{-1}\circ S_2^{-1}}
\circ{\bf H}_2^+
$$
\end{lemma}

{\em Proof:} Verification by composition of the associated tangles.\ep

The explicit formulae for the linear maps associated to the generators
of the mapping class group and the handle attachments in Section~7 allow us now to 
compute the braided Hopf algebra structure induced on 
${\cal N}_0={\cal V}_{\cal N}(\Sigma_{1,1})$. We write 
$M_0:={\cal V}_{\cal N}({\bf M})$, $\Delta_0:={\cal V}_{\cal N}({\bf \Delta})$, 
$S_0:={\cal V}_{\cal N}(S_1^2)$, and $c_0:={\cal V}_{\cal N}({\bf c})$ for
the braided multiplication, comultiplication, antipode and braid isomorphism
respectively. 

\begin{lemma}\lbl{lm-super-Hopf}
The induced braided Hopf algebra structure on ${\cal N}_0$ is the canonical $\Z/2$-graded 
Hopf algebra with:
$$ 
M_0(x\otimes y)= xy\qquad\quad c_0(x\otimes y)=(-1)^{d(x)d(y)}y\otimes x\qquad\quad
\forall x,y\in {\cal N}_0
$$
$$
\mbox{and} \qquad \Delta_0(w)=w\otimes 1 + 1\otimes w \qquad\quad 
\Gamma_0(w)=-w \qquad\forall w\in\E\quad.
$$
In particular, ${\cal N}_0$ is commutative and cocommutative in the graded and braided sense,
${\cal N}_0\cong{\cal N}_0^*$ is self dual, $\SL(2,\R)$ still acts by Hopf automorphisms 
on ${\cal N}_0$, and $S_0$ is an  involutory homomorphism on ${\cal N}_0$.
\end{lemma}
 
{\em Proof:} For {\bf M} and ${\bf \Delta }$ insert the morphism associated to 
the generators in Lemma~\ref{lm-multHeeg}. The braid isomorphism is given via the Hennings
rules by acting with the operator $ad\otimes ad({\cal R})$ on 
${\cal N}_0^{\otimes 2}$ and then permuting the factors. It is easy to see that
$ad\otimes ad({\cal Z})$ acts on $x\otimes y$ by multiplying $(-1)^{d(x)d(y)}$,
where $d(x)$ is the $\Z/2$-degeree of $x$ in ${\cal N}_0$. Moreover, we we know 
that the adjoint action of ${\cal N}_0$ on itself is trivial so that
the term $\th\otimes K\tb$ in the second factor of ${\cal R}$ in (\ref{eq-ARmat})
does not contribute.
\ep

\smallskip

\medskip

\head{8. Skein theory for ${\cal V}_{\cal N}$}\lbl{S8}

The skein theory of the Hennings calculus over ${\cal N}$ is mostly a consequence
 on the form 
$v=1+\rho$ of the ribbon element 
as in (\ref{lm-Aribb}). In the Hennings procedure we substitute a strand 
with decoration $\frac 1 i\rho$ by a dotted strand (with possibly more decorations)
as shown on the left of Figure~\ref{fig-rhostrand}. Observe from (\ref{eq-Amono}) that 
$$
{\cal M}^{\pm 1}(1\otimes \rho)=(1\otimes \rho)\qquad\mbox{and}
\qquad {\cal M}^{\pm 1}(\rho\otimes 1)=(\rho\otimes 1).
$$
This means that for a dotted strand we do not have to distinguish between over
and undercrossing with other strands as indicated on the right of Figure~\ref{fig-rhostrand}.
As a result such a strand can be 
disentangled from the rest of the diagram.  

\begin{figure}[ht]
\begin{center}
\leavevmode
\epsfbox{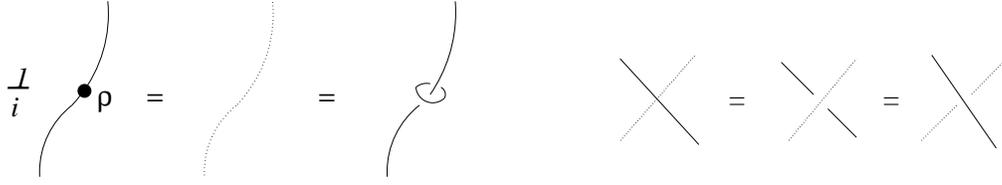}
\end{center}
\caption{Transparent $\rho$-decorated strand}\lbl{fig-rhostrand}
\end{figure}

The next additional ingredient in the calculus are symbols for 1-handles. They are used
in the bridged link calculus as described in \cite{Ker99} and \cite{KerLub00}. We indicate 
a pair of 1-surgery balls by pairs of coupons. The defining relation is the modification
move depicted on the left of Figure~\ref{fig-modif}. The move indicated on the right 
of Figure~\ref{fig-modif} and its reflections  is a standard consequence of the boundary
move from (\ref{eq-boundmove}). 

\begin{figure}[h]
\begin{center}
\leavevmode
\epsfbox{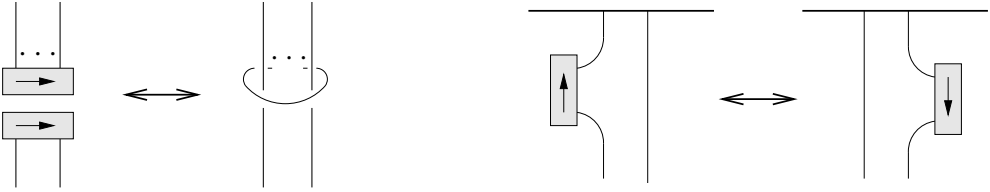}
\end{center}
\caption{Coupons for 1-handles}\lbl{fig-modif}
\end{figure}

Since $v^{k}=1+k\rho$ for $k\in\Z$ we find that the framing of any component can be changed
at the expense of introducing dotted lines. This translates to the diagrams in 
Figure~\ref{fig-frame}. 

\begin{figure}[h]
\begin{center}
\leavevmode
\epsfbox{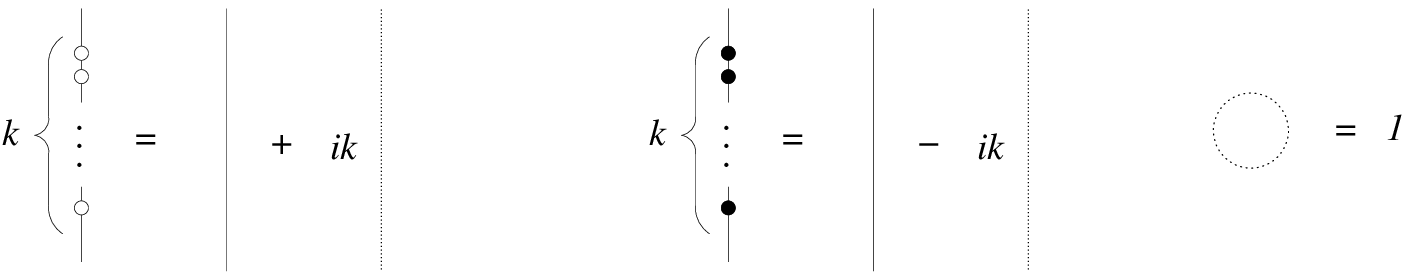}
\end{center}
\caption{Framing shift}\lbl{fig-frame}
\end{figure}

The skein relation is now obtained by applying  Figure~\ref{fig-frame} to
the Fenn Rourke move as in Figure~\ref{fig-FR}, see also\cite{MatPol94}. 
           
\begin{figure}[h]
\begin{center}
\leavevmode
\epsfysize=2.6cm
\epsfbox{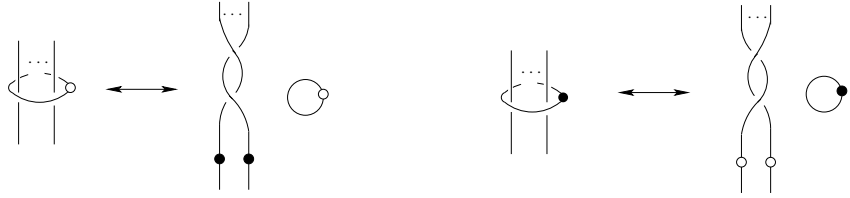}
\end{center}
\caption{Fenn Rourke Move}\lbl{fig-FR}
\end{figure}

\begin{lemma}\label{lm-skein}  For two strands belonging to two different components of a tangle
diagram we have the relation
\begin{center}
\leavevmode
\epsfbox{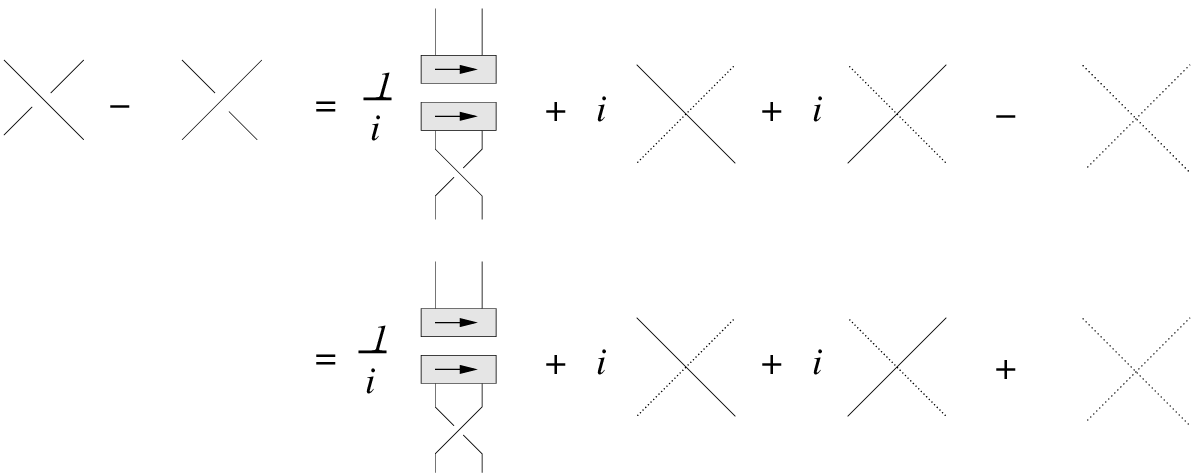}
\end{center}
For strands belonging to the same component of the tangle the relation is 
\begin{center}
\leavevmode
\epsfbox{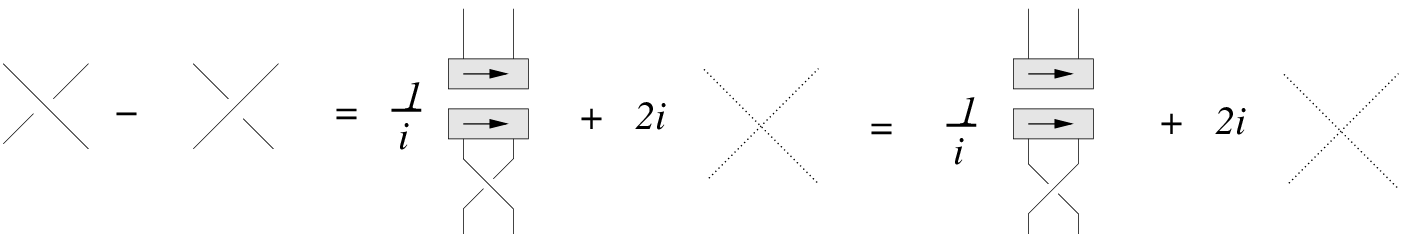}
\end{center}
\end{lemma} 

At this point it is convenient to extend the tangle presentations to general diagrams,
dropping the condition that a strand starting at a point $j^-$ has to end at a point
$j^+$ (or the corresponding condition for through strands). From such a general tangle 
diagram we can get to an admissible one by applying boundary moves (\ref{eq-boundmove})
at all intervals $[j^-,j^+]$. (This is in fact the original definition used in 
\cite{Ker99}.) We shall allow the occurrence of coupons but restrict ourselves to the cases
where exactly two strands enter (or exit) a coupon as in Lemma~\ref{lm-skein}. 

We also introduce two notions of components:
The first  is that of a {\em diagram component} $\cal X$  of a generalized tangle diagram.
It is given by a concatenation of curve segments, coupons that have two strands going in
on one side, and intervals  $[j^-,j^+]$ connecting a strand  ending in $j^-$ with the one
ending in $j^+$. 

The second is a  {\em strand component}, which is also a collection of curves that can be 
joined in two ways.  As before curves that end in two sides of the same interval 
$[j^-,j^+]$ belong to the same strand component, as well as  curves exiting and entering a coupon
pair that would be connected under application of Figure~\ref{fig-modif}.

We have the following rules for manipulating the coupons:

\begin{lemma} \label{lm-1handlemoves}
In the following equivalences the labels $A, B,\ldots$ indicate which coupons form a pair. 
\begin{enumerate}
\item 1-handles can be slid over other 1-handles, through a boundary interval, and hence 
anywhere along a strand component. 
\beq\label{eq-1handlemove}
\epsfbox{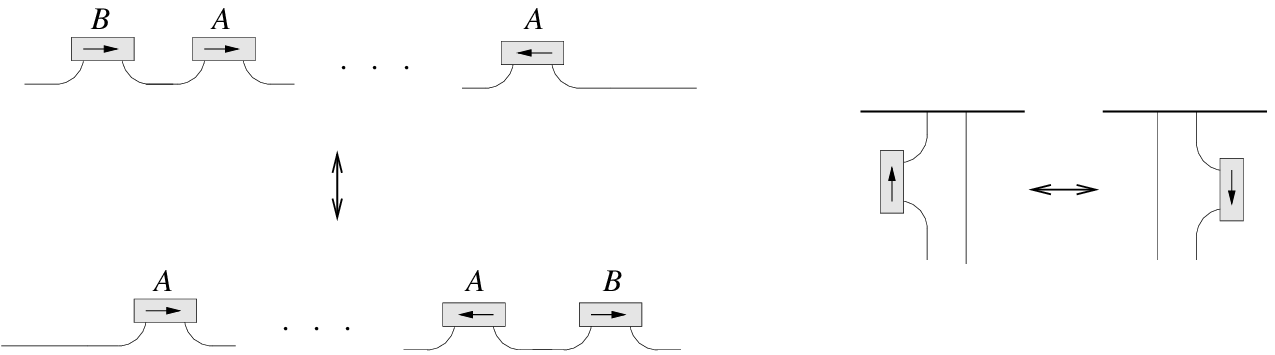}
\eeq
\item If in a diagram the coupons of a pair belong to different diagram components the entire
diagram does not contribute, i.e., is evaluated as zero.  Hence only diagrams contribute in which
the diagram components coincide with strand components. 
\beq\label{eq-pairvanish}
\epsfbox{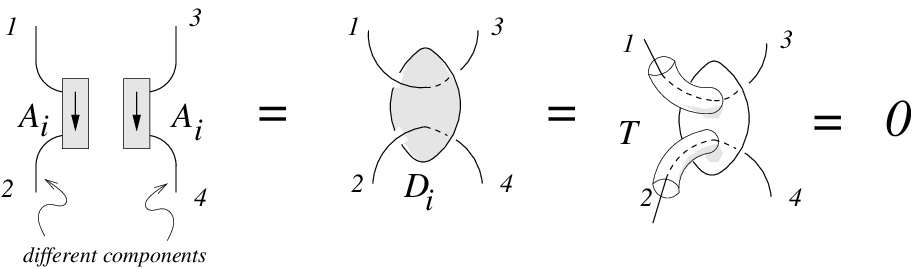}
\eeq
\item Direct 1-handle cancellation: If coupons with the same label are adjacent on the same 
side of a strand they can be canceled: 
\beq\label{eq-1handlecanc}
\epsfbox{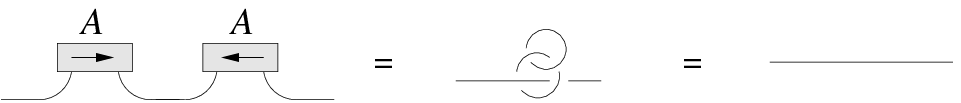}
\eeq
\item Opposite 1-handle cancellation: If coupons with the same label are adjacent on opposite
sides of a strand the strand is replaced by a dotted strand and the evaluation gains a 
factor of 4. 
\beq\label{eq-1handleoppcanc}
\epsfbox{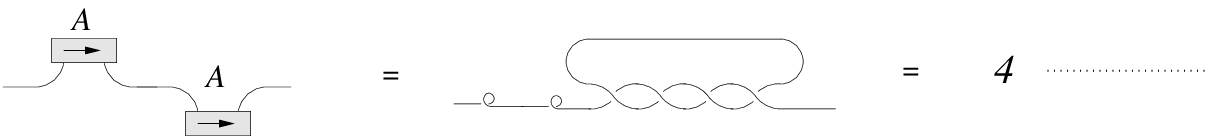}
\eeq
\item If a generalized tangle diagram contains a coupon configuration as indicated the entire 
diagram is evaluated as zero.
\beq\label{eq-coupontorus}
\epsfbox{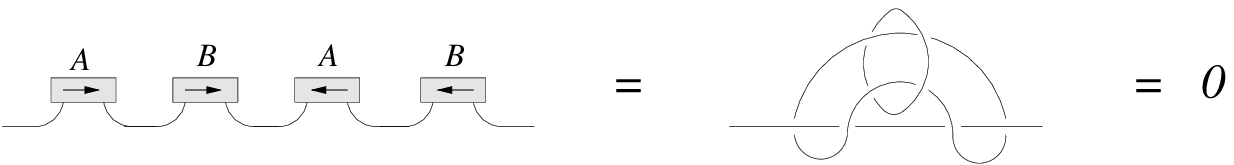}
\eeq
\end{enumerate} 
\end{lemma}

{\em Proof:} The slide of $B$ over the pair $A$ in (\ref{eq-1handlemove}) translates to
a simple isotopy if we apply the move in 
Figure~\ref{fig-modif} to the $A$-pair. Similarly, the slide
through a boundary interval is given by an isotopy conjugated by a $\sigma$-move as in
(\ref{eq-boundmove}).

For {\em b)}
let  $\cal X$ be a diagram component that contains coupons $A_1,\ldots, A_n$ whose partner
lie on different diagram components. Performing boundary moves we can make $\cal X$ to be a true inner
component. Furthermore, we can eliminate the other coupons on $\cal X$ that occur in pairs
by undoing the modification from Figure~\ref{fig-modif}. The component $\cal X$ is now a 
closed curve interrupted only by coupons $A_1,\ldots, A_n$. We undo the modification also
for these and the corresponding annuli added in the move 
 bound discs that we denote by $D_1,\ldots, D_n$. Note, that the arcs of
$\cal X$ all end in only one side of a disc $D_j$ since the strands emerging from the other
side belong to a different component. We can thus surger the discs along the arcs, as shown in 
(\ref{eq-pairvanish}), so that we 
obtain a torus $T$ with $n$  holes 
$\partial T=\partial D_1\sqcup\ldots\sqcup\partial D_n$
which misses all other parts of the tangle. After surgery along the annuli the
torus $T$ can be capped off so that we have found a non-separating surface 
inside the represented cobordism. Since we are dealing with a non-semisimple
TQFT this implies that the associated linear map is zero.

The direct cancellation in (\ref{eq-1handlecanc}) follows by applying Figure~\ref{fig-modif}.
In the resulting configuration in the middle of (\ref{eq-1handlecanc}) the Hopf link can 
be slid off and removed. 

The opposite cancellation in (\ref{eq-1handleoppcanc}) and the remodification from 
Figure~\ref{fig-modif} give the tangle in the middle. Now consider in general a 
straight strand that is entangled with an annulus with $2p$ positive crossings as
in (\ref{eq-annuluswinding}). 
\beq\label{eq-annuluswinding}
\epsfbox{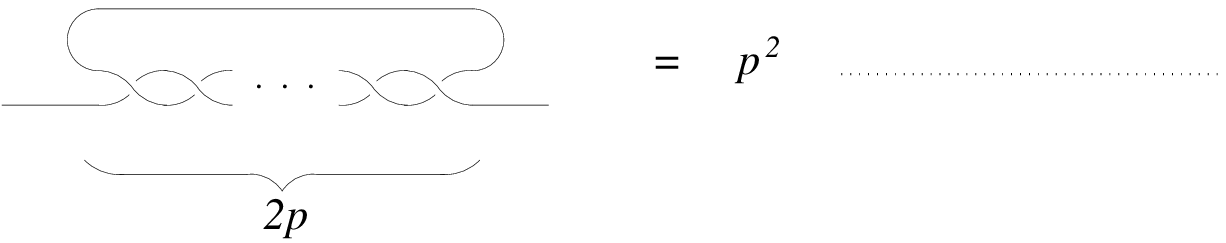}
\eeq
Using  the formula in (\ref{eq-qAmono}) we find by applying the Hennings procedure and
 evaluating the elements on the annulus against the integral
that the resulting element on the open strand
is 
$$
\mu\otimes id({\cal M}^p)\;\;=\;\;\frac {p^2} i \,\rho
$$
which with Figure~\ref{fig-rhostrand} implies the claim. 

Finally, we also 
reexpress the  coupons in 
in (\ref{eq-coupontorus}) by a tangle. As before non-semisimplicity 
of the TQFT implies that a diagram containing such a subdiagram is always zero. For example the
0-framed annulus clearly bounds a surface disjoint from the rest of the link so that the 
cobordism contains a non separating surface. 

\ep

We  now combine the previous two lemmas in the following skein relations without coupons. 

\begin{theorem} \label{thm-skein}
For generalized tangle diagrams we have the following skein relations:

For crossings of strands of different components: 

\beq\label{eq-SKEIN}
\epsfbox{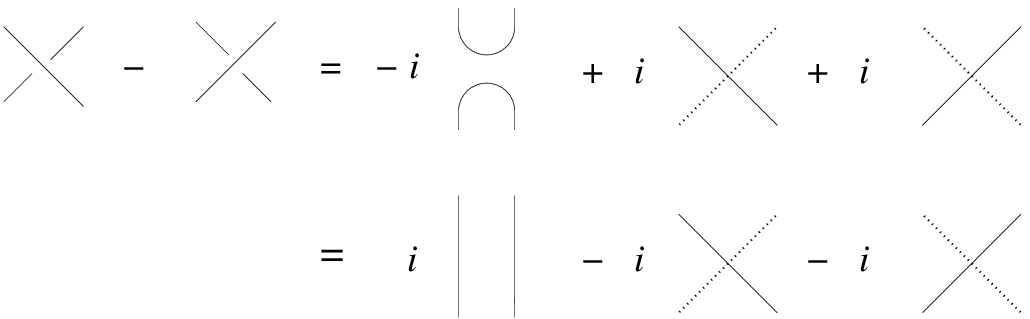}
\eeq

For crossing of strands of the same component we need to introduce an orientation on the component. 

\beq\label{eq-SKEINsame}
\epsfbox{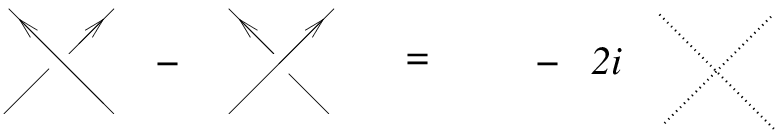}
\eeq

\end{theorem} 

  {\em Proof:}   The proof is given by moving the coupons in the skein relations of 
Lemma~\ref{lm-skein}  through the components using   Lemma~\ref{lm-1handlemoves}. 
  \ep

Note, that relation (\ref{eq-SKEINsame}) implies the relation for the Kauffman 
polynomial for $z=\frac 12$. However, the framing relations are quite different. 

Let $\widehat {\cal B}_g$ be  the group  of tangles in $2g$ strands generated 
by the braidings $\bf c$ of double strands  and the 
braided antipodes $\Gamma$ as in Figure~\ref{fig-MC-tgl} 
acting in different positions.  It is thus the image of the abelian extension  $B_g\ltimes {\Z/2}^g $
of the braid group. 

Moreover, let us introduce a few elementary generalized tangles $M_k: k\to 0$, $\varepsilon:1\to 0$
and $X_n:0\to 2n$ as depicted below. 

\beq\label{eq-planelem}
\epsfbox{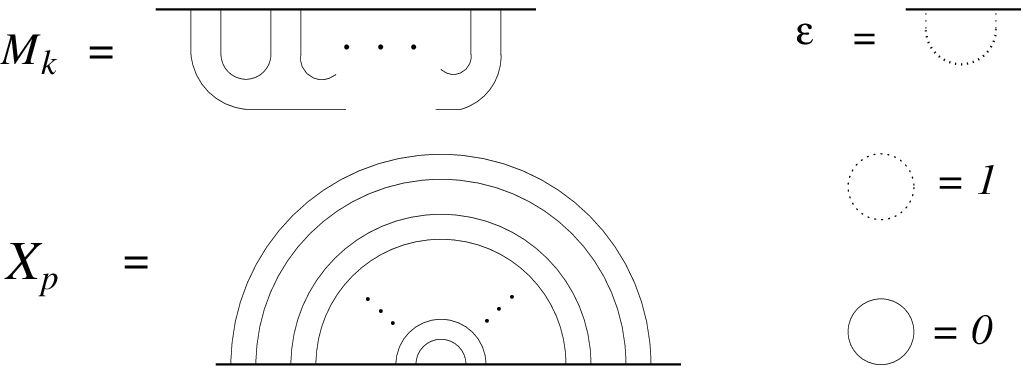}
\eeq

\begin{theorem} \label{thm-solve}
Every tangle $T:G\to 0$   with $2G$ starting (top)  points and {\em no} endpoints 
can be resolved via the skein relations in Theorem~\ref{thm-skein} into a combination
of tangles of the form 
$$
T\;=\;(M_{k_1}\otimes\ldots M_{k_r}\otimes \epsilon^{\otimes N})\circ B\;,
$$
with $B\in \widehat {\cal B}_G$ and $\sum_{i=1}^rk_i=G-N$. 
\end{theorem}

{\em Proof:}    We consider generalized tangles without coupons.  We proceed by induction on 
the number $m$ of connected components of $T$. We only count components that involve solid lines,
those with dotted lines reduce to a collection of $\varepsilon$-diagrams at the intervals belonging to
that component or closed dotted circles that do not contribute.
 Suppose now $T$ has only one component, which we equip with some orientation. 
Applying $\Gamma$'s  to the intervals we can arrange it that the strands enter an interval
$[j^-,j^+]$  at the left point $j^-$ and leave at the right one $j^+$. Furthermore, we can find a 
permutation of intervals so that the strand exiting $j^+$ enters at $(j+1)^-$, except for 
$G^+$, which is connected to $1^-$. Hence, by multiplying an element of $\widehat {\cal B}_G$ 
to $T$ we can assume that the endpoints of the intervals are connected to each other by strands
as they are for $M_G$. 

Next we note that the skein relation (\ref{eq-SKEINsame}) from Theorem~\ref{thm-skein} does not
change this connectivity property for the solid lines and any diagram with dotted lines
collapses to $\varepsilon$-diagrams.

For diagrams where equally labeled coupons are on the same components
there are three planar moves that allow us to
manipulate the arrangement of coupons. They are the 1-handle slide and the
1-handle cancellation depicted below, and the boundary flip as in Figure~\ref{fig-modif}. 
In fact it is easy to see that 
we have the skein relation $T=M_G+iw(T)\varepsilon^{\otimes G}$, where $w(T)$ is the generalization 
of the writhe number of the diagram as defined, for example, in  \cite{Lik97}. In case $G=0$ the diagram
$M_0$ is a closed solid circle which therefore makes the entire diagram zero. 

Assume now $T$ has $m$ components and the claim is true for all diagrams with $m-1$ components.
Pick one component  $C$ and apply an element of $\widehat {\cal B}_G$ such that the intervals included in 
this component are all to the left of the other intervals. Note,
 that the set of intervals that belongs to
$C$ may also be empty. Next apply the skein relations (\ref{eq-SKEIN})
from Theorem~\ref{thm-skein} to untangle $C$ from the other components. In each step of changing 
crossings of a strand of $C$ with the strand of another component $D$ we can choose the relation 
for  which the tangle that belongs to the first local diagram on the right side of the equation has
one component less since $C$ and $D$ are connected. The other diagrams on the right side also
have one less component since we do not count dotted lines. Hence, by induction, the error of changing 
a crossing between $C$ and another component can be resolved into elementary diagrams as claimed. 
After $C$ is untangled  we have expressed $T$, modulo elementary diagrams,
 in the form $C\otimes T'$ (juxtaposition) where $T'$
has $m-1$ components. Again each factor can be resolved independently by induction, and, hence, 
the whole
diagram since $\otimes$-products of elementary diagrams are again elementary. 

\ep

Next note that every tangle $R:\,g_1\to g_2$ is in fact of the form 
\beq\label{eq-RT}
R\;=\; (T\otimes id_{g_2})\circ (id_{g_1}\otimes X_{g_2})
\eeq
for some $T:\,g_1+g_2\to 0$.   Thus, in order to evaluate a general tangle diagram it suffices
by Theorem~\ref{thm-solve}  to specify the evaluations of  the  elementary tangles in (\ref{eq-planelem}). To this end we define the tensor
\beq\label{eq-defAtens}
A\;=\;\frac 1 i S\otimes 1\Delta(\rho)\;=\;\frac 1 i \Bigl( \rho\otimes 1\,+\,1\otimes\rho\,-\,\tb\otimes\th\,+\,\th\otimes\tb\Bigr)\;\;\in\,{\cal N}_0^{\otimes 2}\;.
\eeq

\begin{cor} \label{cor-eval}
Every diagram can be resolved into a sum of composites of diagrams in 
(\ref{eq-planelem}). The linear maps associated to them are 
\beq\label{eq-X1}
{\cal V}_{\cal N}(X_1)\,:\,\Cc\to {\cal N}_0^{\otimes 2}\,:\;1\,\mapsto
A\,=\,\sum_{\nu}x_{\nu}\otimes y_{\nu}
\eeq
\begin{eqnarray}\label{eq-Xn}
{\cal V}_{\cal N}(X_n)\;&=&\;
(1^{\otimes (n-1)}\otimes {\cal V}_{\cal N}(X_1)\otimes1^{\otimes (n-1)})\circ
{\cal V}_{\cal N}(X_{n-1})\,
\;
:\, \Cc\to {\cal N}_0^{\otimes 2n}\,
\\
&:&\;
1\,\mapsto
A_{\{n\}}\,
=\,\sum_{\nu_1,\ldots,\nu_n}x_{\nu_1}\otimes x_{\nu_2}\otimes \ldots\otimes x_{\nu_n}\otimes y_{\nu_n}
\otimes \ldots\otimes y_{\nu_2} \otimes y_{\nu_1}\nonumber
\end{eqnarray}
\beq\label{eq-Mn}
{\cal V}_{\cal N}(M_n):\;\; {\cal N}_0^{\otimes n}\to \Cc\;\;:\;\;\;\;a_1\otimes\ldots\otimes a_n\;\mapsto\;\mu(a_1\cdot\ldots\cdot a_n)
\eeq

Dotted circles can be removed and diagrams with solid circles do not contribute. 
\end{cor}

{\em Proof:} The formulae follow easily from the pictures in Figure~\ref{fig-MC-tgl} to which we 
assigned linear maps in Lemma~\ref{lm-super-Hopf}. Particularly, we find that the upside down
reflection of the multiplication tangle {\bf M} is mapped to the S-conjugate coproduct
\beq\lbl{eq-Scoprod}
\widetilde \Delta= i \Ss^{-1}\otimes\Ss^{-1}\Delta_0\Ss\,:\;{\cal N}_0\otimes{\cal N}_0\to
{\cal N}_0\;.
\eeq
The tangle $X_1$ is obtained by capping this off with an arc at the top, which corresponds to
the insertion of the unit. Hence, $A=\widetilde \Delta(1)$. The diagrams $M_p$ are easily 
identified as composites $M^p=(M\otimes 1^{\otimes(p-1)})\circ M^{p-1}$ capped off with
an arc at the  bottom, which is hence assigned to the $p$-fold multiplication followed
by an evaluation against the integral $\mu\in{\cal N}^*$. 
\ep
\medskip

Let us consider a few examples. One useful case is when the braid $B\in\widehat{\cal B}_n$ can 
be chosen trivially. Hence the contribution to the linear map for a tangle $R:g_1\to g_2$ is
given by a union of planar diagrams as depicted in (\ref{eq-plancomp}):
\beq\label{eq-plancomp}
\epsfbox{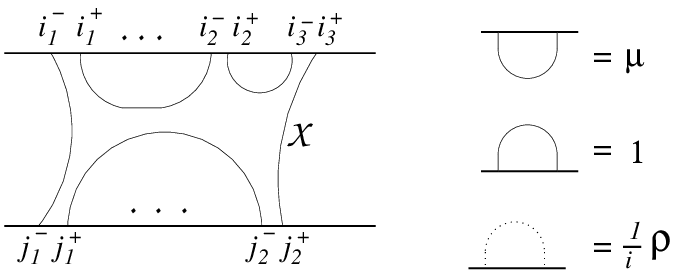}
\eeq
Define the map
\beq\label{eq-defCpq}
C_p^q\,=\,\widetilde\Delta^{q-1}\circ M_0^{p-1}\;:
\quad {\cal N}_0^{\otimes p}\longrightarrow {\cal N}_0^{\otimes q}\;,
\eeq
where the exponents denote the usual multiple products and coproducts. 
The linear map associated to a planar diagram is now the tensor product
of maps associated to the individual components of the diagram. For example, 
if we want to evaluate the linear map on a homogeneous vector $x_1\otimes\ldots\otimes x_{g_1}$   
an the diagram has 
a component with solid lines as in (\ref{eq-plancomp}) containing
top intervals $[i_1^-,i_1^+], \ldots, [i_p^-,i_p^+]$ and bottom intervals
$[j_1^-,j_1^+], \ldots, [j_q^-,j_q^+]$ we compute the vector 
$C_p^q(x_{i_1}\otimes\ldots \otimes x_{i_p})\in {\cal N}_0^{\otimes q}$ and insert the
entries in order into the positions $j_1,\ldots, j_q$ in ${\cal N}_0^{\otimes g_2}$.

With these rules the computation of the maps associated to the generators of the
mapping class group are readily carried out. For
 example we can evaluate the diagram for the $S$-transformation from 
Figure~\ref{fig-S-tgl}. We resolve the right most crossing by taking the 
skein relation in the first row in Proposition~\ref{thm-skein} but with every
diagram rotated  clockwise by $\frac \pi 2$. The result is 
$$
\Ss\;=\;id - \rho\otimes\mu_0 - 1\otimes \epsilon - 1\otimes \mu_0 + \rho\otimes \epsilon
$$
This  yields exactly  the formula from (\ref{eq-SSform}).

As another example we may consider the $C_1$ waist cycle in $\Sigma_2$. The diagram consists
of four parallel strands with a 1-framed annulus around the second and third. We apply 
Figure~\ref{fig-frame} and then  Figure~\ref{fig-modif} to this annulus. The resulting coupons
can be canceled. We find 
$${\cal V}_{\cal N}({\bf I}_{C_1})=id - iC_1^1\,.
$$
This implies the formula for the $D$-transformation from (\ref{eq-D-2act}). 
\smallskip

Finally, let us show how to use the skein calculus to find
 the precise formula for the invariant of a 2-framed
closed 3-manifold presented by a link ${\cal L}\subset S^3$. It is 
basically given by the order of the first integral homology. More precisely,
let
\beq\label{eq-defeta}
\eta(M)\;:=\;\left\{\begin{array}{cl}
\Bigl| H_1(M,\Z)\Bigr| & \mbox{for}\; \beta_1(M)=0\\
0 & \mbox{for}\;\beta_1(M)>0\end{array}\right.
\eeq

\begin{lemma}\label{lm-defeta} For a given framed link ${\cal L}\subset S^3$ and 
$\eta$ as in (\ref{eq-defeta}) we have 
$$
{\cal V}_{\cal N}(M_{\cal L})\;=\; i^{|{\cal L}|}det({\cal L}\cdot{\cal L})
=\;\pm i^{|{\cal L}|} \eta(M)
$$
\end{lemma}

{\em Proof:}
By 2-handle slides we can move ${\cal L}$ into a link ${\cal L}^{\delta}$ 
so that the linking form ${\cal L}^{\delta}\cdot {\cal L}^{\delta}$ is
diagonal and equivalent to the original one ${\cal L}\cdot{\cal L}$. Suppose 
$f_j$ is the framing number of the $j$-th component ${\cal L}_j^{\delta}$. 
From Figure~\ref{fig-frame} we see that 
$$
{\cal V}_{\cal N}({\cal L}^{\delta})= 
{\cal V}_{\cal N}({\cal L}^{\delta, -f_j}) + if_j
{\cal V}_{\cal N}({\cal L}^{\delta}-{\cal L}_j^{\delta})
$$
Here, ${\cal L}^{\delta, -f_j}$ is the link in which the framing of the $j$-th
component is shifted to zero. As a result the manifold represented by this
link has non-trivial rational homology. Since ${\cal V}_{\cal N}$ is  a 
non-semisimple theory this implies that 
${\cal V}_{\cal N}({\cal L}^{\delta, -f_j})=0$. Iterating the above identity
we find 
${\cal V}_{\cal N}({\cal L}^{\delta})=\prod_{j=1}^{|{\cal L}|}(if_j)
{\cal V}_{\cal N}(\emptyset)$. Clearly, $\prod_{j=1}^{|{\cal L}|}(f_j)$
is the determinant of the linking form of ${\cal L}^{\delta}$ and 
hence also the one of ${\cal L}$. 

\ep

\head{9. Equivalence of ${\cal V}_{\cal N}^{(2)}$ and ${\cal V}^{FN}$}\lbl{S9}

 In this section we compare the two topological quantum field theories ${\cal V}^{FN}$
described in Section~3 and ${\cal V}_{\cal N}^{(2)}$ constructed in Section~7. We already
found  a number of general properties that are shared by both theories:

By Lemma~\ref{lm-FNpart} and Proposition~\ref{pp-ATQFT} both theories are $\Z/2$-projective 
on $\Cob^{\bullet}$ and non-semisimple, fulfilling the property of Lemma~\ref{lm-vanish}. 
The $\Z/2$-projectivity is due to ambiguities of even 2-framings in the case of 
${\cal V}_{\cal N}^{(2)}$ and ambiguities of orientations in the case of ${\cal V}^{FN}$.
The non-semisimple half-projective property results in the case of ${\cal V}^{FN}$ from
representation varieties that are transversely disjoint, and in the case of 
${\cal V}_{\cal N}^{(2)}$  from the  nilpotency
of the integral $\lambda\in\cal N$. 
Further common features are the dimensions of vector spaces $(=4^g)$,
actions of $\SL(2,\R)$, see Section~9, and the fact that ${\cal J}_g$ lies in the kernel
of the mapping class group representations. 

We construct now an explicit isomorphism between ${\cal V}^{FN}$ and ${\cal V}_{\cal N}^{(2)}$. 
Let ${\cal Q}=\ext *\lz a,b\rz$ be the  exterior algebra over $\R^2$ with basis $a,b\in\R^2$.
We obtain a canonical isomorphism, which is defined on monomial elements as follows: 
\begin{equation}\lbl{eq-defi}
i_*\;:\quad {\cal Q}^{\otimes g}\;\isto\;\ext * H_1(\Sigma_g)\;\;:\qquad 
q_1\otimes \ldots\otimes q_g\mapsto i_1(q_1)\wedge\ldots\wedge i_g(q_g) \;, 
\end{equation}
where $i_j:{\cal Q}\,\isto\,\ext *\lz [a_j], [b_j]\rz$ is the canonical map sending $a$ and $b$ 
to $[a_j]$ and $[b_j]$ respectively. Next, 
we define an isomorphism between ${\cal Q}$ and ${\cal N}_0$, seen as linear spaces, 
by the following assignment of basis vectors: 
\begin{equation}\lbl{eq-phidef}
\begin{array}{lcc}
&\quad \phi(1)=b \quad
&\quad\phi(\tb\th)=a\quad\\
\phi\;:\quad{\cal N}_0\;\isto\;{\cal Q}\qquad\;\;\mbox{with}\;\;&&\\
 &\qquad \phi(\th)=a\wedge b \quad &\quad \phi(\tb)=1\quad \\
\end{array}\;.
\end{equation}
Note, that this map has odd $\Z/2$-degree and is, in particular,
 not an algebra homomorphism. From (\ref{eq-phidef}) we infer 
directly the following identities:
\beq\lbl{eq-}
\phi(\th x)=-\phi(x)\wedge a \qquad\qquad \phi(x \th )=a\wedge  \phi(x)\\
\eeq
\beq\lbl{eq-AS1phi}
\phi(\Aa x)\;=\;[A_1]\phi(x)\qquad\qquad \phi(\Ss x)\;=\;[S_1]\phi(x)\\
\eeq
Here, $\Aa$ and $\Ss$ are as in (\ref{eq-AVgen}), and $[A_1]$ and $[S_1]$ are the 
maps on $H_1(\Sigma_1)$ as in (\ref{eq-A-hom}) and (\ref{eq-S-hom}).

Moreover, let us introduce a sign-operator $(-1)^\Lambda$ on ${\cal Q}^{\otimes g}$ 
defined on monomials by 
\begin{equation}\lbl{eq-Lambda}
(-1)^{\Lambda_g}(q_1\otimes \ldots\otimes q_g)\quad=\quad
(-1)^{\lambda_g(d_1,\ldots,d_g)}q_1\otimes \ldots\otimes q_g\;.
\end{equation}
The function $\lambda_N$ is defined in the $N$-fold  product of 
 $\Z/2$'s as follows:
\begin{equation}\lbl{eq-lambdadef}
\lambda_N:\;(\Z/2)^N\,\to\,\Z/2\qquad\mbox{with}\quad
\lambda_N(d_1,\ldots,d_N)\;=\;\sum_{i<j}d_i(1-d_j)\;,
\end{equation}
where $d_j={\rm deg}(q_j)\,{\rm mod}\, 2\,$. Consider now the following isomorphism
of vector spaces. 
\begin{equation}\lbl{eq-defxi}
 \xi_g:=\;i_*\circ(-1)^{\Lambda_g}\circ\phi^{\otimes g}\quad:
\qquad{\cal N}_0^{\otimes g}\,\;\isto\,\;\ext *H_1
\end{equation}
Given a linear map, $F:{\cal N}^{\otimes g_1}\to {\cal N}^{\otimes g_2}$, we write 
$(F)^{\xi}:=\xi_{g_2}\circ F\circ \xi_{g_1}^{-1}$ for the respective map on homology.
Moreover, we denote by  ${\bf L}^{(k)}_x$ the operator on ${\cal N}^{\otimes g}$  
that multiplies the $k$-th factor in the tensor product by $x$ from the left,
and by ${\bf R}^{(k)}_x$ the respective operator for multiplication from the right. We
compute: 
\beq\lbl{eq-LRdef}
\begin{array}{cl}
\quad &({\bf L}^{(k)}_{\th})^{\xi}(\alpha\wedge u_k\wedge\beta)\;\;=\;\;(-1)^{g-k+s+1}\,\,
\alpha\wedge a_k\wedge u_k\wedge\beta\;,\\
\\
\mbox{and}\qquad &({\bf R}^{(k)}_{\th})^{\xi}(\alpha\wedge u_k\wedge\beta)\;\;=\;\;(-1)^{g-k+s}\,\,
\alpha\wedge u_k\wedge a_k \wedge\beta\;,
\end{array}
\eeq
where $s=\sum_{j=1}^g d_j$ is the total degree of $\alpha\wedge u_k\wedge\beta$, 
$\alpha\in\ext *\lz a_1,\ldots, b_{k-1}\rz$,  and $\beta\in\ext *\lz a_{k+1},\ldots, b_{g}\rz$. 

\begin{lemma}\lbl{lm-commut}
For every  standard generator $G\in \{A_j, D_j, S_j\}$, we have 
$$
({\cal V}_{\cal N}({\bf I}_{G}))^{\xi}\;=\;\ext * [G]\;\;,
$$
where $[G]$ denotes as before the action on homology. 
\end{lemma}

{\em Proof:} For the $A_j$ and $S_j$ this follows readily from (\ref{eq-AS1phi}), and 
the fact that $[A_j]$ and $[S_j]$ do not change the  degrees $d_j$ and hence commute with 
$(-1)^{\Lambda_g}$. 

The operator in (\ref{eq-D-2act}) decomposes into $\D=\D^0+\D^1$, where 
$\D^0=id - {\bf R}_{\rho}\otimes {\bf L}_{\rho}$ and 
$\D^1={\bf R}_{\th}\otimes {\bf L}_{\tb}-{\bf R}_{\tb}\otimes {\bf L}_{\th}$.
Now $\D^0$ does not change the $\Z/2$-degree of both factors, and $\D^1$ flips the
degree of both factors. One readily verifies that
$$
\lambda_g(\ldots, 1-d_j, 1-d_{j+1},\ldots)-\lambda_g(\ldots, d_j, d_{j+1},\ldots)=d_j + d_{j+1}\qquad
{\rm mod}\,2
$$
$$
\mbox{so that}\qquad
{\cal V}_{\cal N}({\bf I}_{D_j})^{\xi}\;\;=\;\;({\cal V}_{\cal N}^0({\bf I}_{D_j}))^{\zeta}\,+\,
(-1)^{d_j+d_{j+1}}({\cal V}_{\cal N}^1({\bf I}_{D_j}))^{\zeta}\quad
$$
$$\qquad\quad=\;\;(I^{\otimes j-1}\otimes (\D^0)^{\phi^{\otimes 2}}\otimes I^{\otimes g-j-1})^{i_*}\;+\;(-1)^{d_j+d_{j+1}}(I^{\otimes j-1}\otimes (\D^1)^{\phi^{\otimes 2}}
\otimes I^{\otimes g-j-1})^{i_*}
$$

Here, $\zeta_g=i_*\circ \phi^{\otimes g}$ and 
${\cal V}_{\cal N}^i({\bf I}_{D_j})$ is the operator with $\D^i$ in $j$-th position.
Since $\zeta_g=\zeta^{\otimes g}_1$ the $\zeta$-conjugate maps only act on the generators
$\{a_j, b_j, a_{j+1}, b_{j+1}\}$ the action is the same for all positions $j\,$.
Observe that also $[D_j]$ acts only on the homology generators $\{a_j, b_j, a_{j+1}, b_{j+1}\}$.
It is, therefore, enough to prove the relation for $g=2$ and ${\cal V}_{\cal N}({\bf I}_{D_1})
=\D$. 

Now, from (\ref{eq-D-2act}) it is obvious
that ${\cal V}_{\cal N}({\bf I}_{D_j})$ commutes
with ${\bf L}^{(j)}_{\th}$ and ${\bf R}^{(j+1)}_{\th}$. Moreover, it is easy to see that
$\ext *[D_j]$, as given in (\ref{eq-D-hom}), commutes with 
$({\bf L}^{(j)}_{\th})^{\xi}$ and $({\bf R}^{(j+1)}_{\th})^{\xi}$ from (\ref{eq-LRdef}). 
Specifically, we use that $\ext *[D_j]$ does not change the total degree, and acts trivially on
$a_j$ and $a_{j+1}$. It thus suffices to check 
\beq\lbl{eq-12Drel}
\ext 2[D_1]\circ\zeta_2(x_1\otimes x_2)=\zeta_2\circ \D^0(x_1\otimes x_2) +(-1)^{d_1+d_2}
\zeta_2\circ\D^1 (x_1\otimes x_2)
\eeq
 with $d_i={\rm deg}(\phi(x_i))$, and only for $x_i\in\{1,\tb\}$. For example, for $x_1=x_2=1$,
with $d_1+d_2=0$, we find from (\ref{eq-D-2act}) and (\ref{eq-D-hom}) that
$$
\begin{array}{lcl}
\zeta_2 \circ \D(1\otimes 1)&=&\zeta_2(
1\otimes 1  +  \th\otimes \tb  -
\tb\otimes \th  -  \rho \otimes \rho)\;\\
&=& b_1\wedge b_2 + a_1\wedge b_1 - a_2\wedge b_2 - a_1\wedge a_2\\
&=& (b_1-a_2)\wedge (b_2-a_1)  \;=\;\ext 2[D_1](b_1\wedge b_2)\;=\;\ext 2[D_1](\zeta_2(1\otimes 1))\\
\end{array}
$$
We also compute for the case $x_1=\tb$ and $x_2=1$, with $d_1+d_2=1$:
$$
\begin{array}{lcl}
\zeta_2 \circ (\D^0-\D^1)(\tb\otimes 1)&=&\zeta_2(\tb\otimes 1 -\tb\th\otimes\tb)
\;=\;b_2 - a_1\\
&=& \ext 2[D_1](b_2)\;\;\;=\;\;\;\ext 2[D_1](\zeta_2(\tb\otimes 1))
\end{array}\;.
$$ 
The other two cases follow similarly.
\ep

As the $\{A_j, D_j, S_j\}$ generate $\Gamma_g$ we conclude from Lemma~\ref{lm-commut} and
(\ref{eq-FNmapcg})  that $({\cal V}_{\cal N}({\bf I}_{\psi}))^{\xi}
={\cal V}^{FN}({\bf I}_{\psi})$ for all $\psi\in\Gamma_g\,$.

Let us also consider the maps associated by both functors to the
handle additions  ${\bf H}_g^{\pm}$. We note that
$$
\lambda_{g+1}(d_1,\ldots,d_g,1)=\lambda_g(d_1,\ldots,d_g)
$$
so that we find from (\ref{eq-Ahdl+map}), (\ref{eq-FNmcg}) and (\ref{eq-phidef}) that
$({\cal V}_{\cal N}({\bf H}_g^+))^{\xi}={\cal V}^{FN}({\bf H}_g^+)$. Similarly,
(\ref{eq-Ahdl-map}), (\ref{eq-FNgcm}) and (\ref{eq-Ainteg}) imply 
$({\cal V}_{\cal N}({\bf H}_g^-))^{\xi}={\cal V}^{FN}({\bf H}_g^-)$. 
Using the Heegaard 
decomposition (\ref{eq-Heegaard}) we finally infer equivalence:
\begin{propos}\lbl{pp-Vequiv}
The maps $\xi_g$ defined in (\ref{eq-defxi}) give rise to an isomorphism
$$
\xi\;:\; {\cal V}_{\cal N}\;\;\;\;\stackrel{\bullet\,\cong}
{-\!\!\!-\!\!\!-\!\!\!\longrightarrow}\;
\;\;\;{\cal V}^{FN}\;\;
$$
of relative, non-semisimple, $\Z/2$-projective functors from $\Cob^{\bullet}$ to $Vect(\kk)$. 
 \end{propos}

\head{10. Hard-Lefschetz decomposition and Invariants}\lbl{S10}
 
The tangent bundle over the  moduli space $J(\Sigma_g)$ is trivial with fiber
$H^*(\Sigma_g,\R)$ so that its cohomology ring is naturally $\ext *H_1(\Sigma_g,\R)$. 
There is an almost complex structure on $J(\Sigma_g)$ given by a map ${\sf J}$ with
 ${\sf J}^2=-1$ in the cohomology. It is given by ${\sf J}.[a_j]=-[b_j]$ and 
${\sf J}.[b_j]=[a_j]$.
With the  K\"ahler form
$\omega_g\in H^2(J(\Sigma_g))$ defined in (\ref{eq-defsymform}) 
it is also a K\"ahler manifold. The dual K\"ahler metric provides us with a Hodge star
$\star:\ext jH_1(\Sigma_g)\to\ext {2g-j}H_1(\Sigma_g)$ for a given volume form
$\Omega\in\ext {2g} H_1(\Sigma_g)$ by the equation 
$\alpha\wedge\star\beta=\lz\alpha,\beta\rz\Omega$. 
Specifically, the $2g$ generators  $\{[a_1],\ldots, [b_g]\}$ of $H_1(\Sigma_g)$, 
with volume form $\Omega=[a_1]\wedge\ldots\wedge [b_g]$ the Hodge star is given by  
$\star(a_1^{1-\epsilon_1}\wedge\ldots\wedge b_g^{1-\epsilon_{2g}})
=(-1)^{\lambda_{2g}(\epsilon_1,\ldots,\epsilon_{2g})}a_1^{\epsilon_1}\wedge\ldots\wedge 
b_g^{\epsilon_{2g}}$, where $\lambda_{2g}$ is as in (\ref{eq-lambdadef}).

As a K\"ahler manifold $H^*(J(\Sigma_g))$ admits an $\SL(2,\R)$-action,  
see for example \cite{GriHar78},  given for the standard generators  
$E, F, H\,\in\,{\mathfrak s}{\mathfrak l}_2(\R)$ by 
\beq\lbl{eq-defsl2act}
H\alpha:=(j-g)\alpha\quad\; \forall \alpha\in \ext jH_1(\Sigma_g)\;,
\qquad\quad E\alpha:=\alpha\wedge\omega_g\;,\qquad\quad F:=\star\circ E\circ\star^{-1}
\eeq
\begin{lemma}\lbl{lm-FN-equiv}
The functor ${\cal V}^{FN}$ is $\SL(2,\R)$-equivariant with respect to the action
in (\ref{eq-defsl2act}).
\end{lemma}

{\em Proof:}  Commutation with $H$ follows from counting degrees. Since $\omega_g$ is invariant
under the $\Sp(2g,\R)$-action, $E$ commutes with the maps in (\ref{eq-FNmapcg}),
and since $\omega_g\wedge [a_{g+1}]=[a_{g+1}]\wedge \omega_{g+1}$ also with the ones 
in  (\ref{eq-FNmcg}) and  (\ref{eq-FNgcm}). Finally, as all maps ${\cal V}^{FN}(M)$
are isometries with respect to $\lz.,.\rz$ they also commute with $F$. \ep

In order to finish the proof of Theorem~\ref{thm-main} we still need to show that
the $\xi_g$ are $\SL(2,\R)$-equivariant as well. The fact that $H$ commutes with $\xi_g$ 
is again a matter of counting  degrees. We have $E=\sum (E_1^{(i)})^{i_*}$,
 where $E_1^{(i)}$ acts
on the $i$-th factor of ${\cal Q}^{\otimes g}$ 
by $q\mapsto E_1(q)= q\wedge a\wedge b$. Since $E$ does not change
degrees we find that $E^{\xi}=\sum (E^{(i)})^{\phi^{(i)}}$, where
$(E^{(i)})^{\phi^{(i)}}$ acts on the $i$-th factor by $E_1^{\phi}$.
We find $E_1^{\phi}(\tb)=\th$, and 
$E_1^{\phi}(1)=E_1^{\phi}(\th)=E_1^{\phi}(\tb\th)=0$, which yields precisely
 the desired action of $E$ 
on ${\cal N}_0$. The conjugate action of $\star$ on ${\cal N}_0^g$ is as follows: 
\beq\lbl{eq-xistar}
\star^{\xi}\;:\;\; x_1\otimes\ldots x_g\;\;\mapsto \; (-1)^{\sum_{i<j}d_id_j}(\star x_1)
\otimes\ldots\otimes (\star x_g)\qquad\;\forall x_j\in{\cal N}_0\;, 
\eeq	
where $\star \th=\tb$, $\star\tb=\th$, $\star\tb\th=1$, and $\star 1=-\tb\th$. 
From this we see  that $F^{\xi}$ acts on each factor by 
$F_1^{\phi}(\th)=\tb$, and 
$F_1^{\phi}(1)=F_1^{\phi}(\tb)=F_1^{\phi}(\tb\th)=0$, as required.

With Lemma~\ref{lm-FN-equiv} and equivariance of $\xi_g$ 
we have thus completed the proof of Theorem~\ref{thm-main}.
Henceforth, we will use the simpler notation ${\cal V}={\cal V}^{FN}={\cal V}_{\cal N}$
 \ep

\medskip

The $\SL(2,\R)$-action implies a Hard-Lefschetz decomposition \cite{GriHar78} as follows
\beq\lbl{eq-Lefsch}
H^*(J(\Sigma_g))\quad\cong\quad \bigoplus_{j=1}^{g+1} V_j\otimes W_{g,j}\;.
\eeq
Here, $V_j$ is the irreducible ${\mathfrak s}{\mathfrak l}_2$-module 
with ${\rm dim}(V_j)=j\,$, and 
\beq\lbl{eq-W}
W_{g,j}\;:=\;\{u\in\ext {g-j+1}H_1(\Sigma_g)\,:\;\omega_g\wedge u=0\}
\eeq	
is the space of {\em isotropic} vectors of degree $(g-j+1)$, or, equivalently, the
space of ${\mathfrak s}{\mathfrak l}_2$-highest weight vectors of weight $(j-1)$. On each of these
spaces we have an action of the mapping class groups from (\ref{eq-FNmapcg})
factoring through $\Sp(2g,\R)$.
\begin{theorem}[\cite{GooWal98} Chapter 5.1.8]\lbl{thm-GW}
Each $W_{g,j}$ is an irreducible $\Sp(2g,\R)$-module with fundamental highest weight
$\varpi_{g-j+1}$ and dimension 
$$
{\rm dim}(W_{g,j})={2g\choose g-j+1}\,-\, {2g\choose g-j-1}
$$
In particular, the pair of subgroups 
$$
\SL(2,\R)\,\times\,\Sp(2g,\R)\;\;\subset\;\;{\rm GL}(H^*(J(\Sigma_g)))
$$ 
forms a Howe pair, that is, the two subgroups are exact commutants of each other.
\end{theorem}
The fundamental weights are given as in \cite{GooWal98} by 
$\varpi_{k}=\epsilon_1+\ldots +\epsilon_{k}$ with $\epsilon_j$ as
in  (\ref{eq-sproots}).

In the decomposition into irreducible TQFT's the one for $j=1$ associated to the
trivial $\SL(2,\Cc)$ representation plays a special role for invariants of closed
manifolds. 

For any invariant, $\tau$, of closed 3-manifolds there is a standard ``reconstruction'' 
of TQFT vector spaces 
as follows. We take the formal $\kk$-linear span ${\mathfrak C}_g^+$ of 
cobordisms $M:\emptyset\to\Sigma_g$
and ${\mathfrak C}_g^-$ of cobordisms $N:\Sigma_g\to\emptyset$. We obtain a pairing
${\mathfrak C}_g^-\times {\mathfrak C}_g^+\to\kk:\,(N,M)\to \tau(N\circ M)$. If 
${\mathfrak N}^+_g\subset {\mathfrak C}_g^+$ is the null space of this pairing
we define ${\cal V}^{\tau-rec}(\Sigma_g)={\mathfrak C}_g^+/{\mathfrak N}^+_g$.
 For generic $\tau$ these vector spaces are infinite dimensional. The exception is when $\tau$
stems from a TQFT. In this case 
${\cal V}^{\tau-rec}(\Sigma_g)^*={\mathfrak C}_g^-/{\mathfrak N}^-_g$, and the linear map 
${\cal V}^{\tau-rec}(P)$ associated to a cobordism $P$ is reconstructed from its
matrix elements $\tau(N\circ P\circ M)$. This construction, which basically imitates the
GNS construction of operator algebras, is folklore since the emergence of TQFT's and 
appears, for example, in \cite{Tur94}.

\begin{theorem}\lbl{cor-dec}
\ben
\item
The TQFT functor from Theorem~\ref{thm-main} decomposes into a direct sum
$$
{\cal V}\;=\;\;\bigoplus \R^{j}\otimes {\cal V}^{(j)}\;\;\;=\;\;\;
{\cal V}^{(1)}\,\oplus\, \R^2\otimes {\cal V}^{(2)}
\,\oplus\, \R^3\otimes{\cal V}^{(3)}\,\ldots \,
$$
of irreducible TQFT's with multiplicities. 
\item 
The associated vector space for each TQFT is ${\cal V}^{(j)}(\Sigma_g)=W_{g,j}$
so that ${\cal V}^{(j)}(\Sigma_g)=0$ whenever $j>g+1$. 
In particular, for any closed 3-manifold $M$ and $j>1$ we have ${\cal V}^{(j)}(M)=0$ so that
${\cal V}(M)={\cal V}^{(1)}(M)$. 
\item 
The vector spaces associated to the invariant $\pm\eta$ from (\ref{eq-defeta}) 
are finite dimensional. The 
reconstructed $\Z/2$-projective TQFT is ${\cal V}^{\eta-rec}={\cal V}^{(1)}$ with dimensions
${\rm dim}({\cal V}^{\eta-rec}(\Sigma_g))=
{\rm dim}(W_{g,1})=\frac 2{g+2}{2g+1\choose g}\,$.
\een
\end{theorem}

{\em Proof:} The fact that the TQFT's decompose in the prescribed manner follows 
from the $\SL(2,\R)$-covariance. Irreducibility of each ${\cal V}^{(j)}$, meaning there
are no proper sub-TQFT's, results from the fact that each $\Sp(2g,\Z)$ representation
is irreducible so that in a sub-TQFT the vector spaces for each $g$ are either 
${\cal V}^{(j)}(\Sigma_g)$ or 0. Since the handle maps yield non-zero maps between 
these vector spaces if one space is non-zero none of them can be. The reconstructed TQFT
must be a quotient TQFT of ${\cal V}^{(1)}$, which is, however, irreducible. 
Hence, they are equal.

\ep 


\emptystuff{
From the irreducible TQFT's in Theorem~\ref{cor-dec} we can construct a much larger class of TQFT's,
which appear to be  related to higher rank gauge theories, 
as follows. Let ${\cal P}^+
\subset {\Z}^{0,+}[x_1,x_2,\ldots]$ be the set of formal power series 
$$
P(x_1,x_2,\ldots)=\sum_{k=1}^{\infty}\,\sum_{n_1, \ldots , n_k=1}^{\infty}c_{n_1,n_2,\ldots,n_k}
x_1^{n_1}x_2^{n_2}\ldots x_k^{n_k}\;,
$$
such that all $c_{n_1,n_2,\ldots,n_k}\in {\Z}^{0,+}$ are non-negative integers, and for fixed $k$
only finitely many $c_{n_1,n_2,\ldots,n_k}$ are non-zero. To every such $P$ we associate a
TQFT by the formula
\beq\label{eq-defpolyTQFT}
{\cal V}^{(P)}\;=\;\bigoplus_{k=1}^{\infty}\bigoplus_{n_1, \ldots , n_k=1}^{\infty}
\, \R^{c_{n_1,n_2,\ldots,n_k}}\otimes ({\cal V}^{(1)})^{\otimes n_1}\otimes 
({\cal V}^{(2)})^{\otimes n_2}\otimes \ldots\otimes ({\cal V}^{(k)})^{\otimes n_k}
\;.
\eeq
For example ${\cal V}^{FN}={\cal V}^{(F)}$, where $F(x_1,x_2,\ldots)=\sum_j j x_j$. 
The restriction on the coefficients together with the second part of Theorem~\ref{cor-dec}
implies that all vector spaces are finite dimensional. 
\begin{lemma}\label{lm-polyTQFT}
The TQFT functor  ${\cal V}^{(P)}$ is well defined for every $P\in {\cal P}^+$ . 
\end{lemma} 
}


Let us finally give an alternative proof of Lemma~\ref{lm-defeta} using the language in
which the Frohman Nicas invariant is constructed. 

We present $M$  by a Heegaard splitting 
$M_{\psi}=h^-_g\circ{\bf I}_{\psi}\circ h^+_g$, as defined in 
(\ref{eq-Heegaard}) and (\ref{eq-closM}).
 The invariant is given as
the matrix coefficient of $\ext g[\psi]$ for the basis vector
${\cal V}(h^+_g)=[a_1]\wedge[a_2]\wedge\ldots\wedge[a_g]$. If we denote by
$[\psi]_{aa}$ the $g\times g$-block of $[\psi]$ acting on the Lagrangian
subspace spanned by the $[a_i]$'s this number is just ${\rm det}([\psi]_{aa})$.
At the same time, the Mayer-Vietoris sequence for $M_{\psi}$ shows that
$[\psi]_{aa}$ is a presentation matrix for the group $H_1(M_{\psi},\Z)$
so that the order of $H_1(M_{\psi},\Z)$ is, indeed, given by $\pm{\rm det}([\psi]_{aa})$.
  
\ep

\head{11. Alexander-Conway Calculus for 3-Manifolds}\lbl{S11}

Let $M$ be a 3-manifold   with an  epimorphism
$\varphi:H_1(M,\Z)\twoheadrightarrow \Z$.  We recall the definition of
the  {\em (reduced) Alexander polynomial}
$\Delta_{\varphi}(M)$, as it is given in the case of 
knot and link complements for example in \cite{BurZie}.

Let $\widetilde M\to M$
be the cyclic cover associated to $\varphi$ and view $H_1(\widetilde M)$  as
a $\Z[t,t^{-1}]$-module with $t$ acting by Decktransformation. Let 
$E_1\subset \Z[t,t^{-1}]$ be the first elementary ideal generated by the 
$n\times n$ minors of an $n\times m$ presentation matrix $A(t)$ of $H_1(\widetilde M)$. 
Then $\Delta_{\varphi}(M)$ is the generator of the smallest principal idea containing
$E_1$, or, equivalently, the g.c.d. of the $n\times n$ minors of a presentation
matrix. Particularly, if $A(t)$ is a square matrix $\Delta_{\varphi}(M)=det(A(t))$ and
if $n>m$, i.e., there are more rows than columns,  $\Delta_{\varphi}(M)=0$.

Another important invariant of a 3-manifold is its Reidemeister Torsion, which is
obtained as the torsion of a  chain complex over $\Q [t,t^{-1}]$ obtained from a 
cell decomposition of $\widetilde M$. 
The Alexander polynomial turns out to be
 almost the same as the Reidemeister Torsion of a 3-manifold.
The relation described in the next theorem was first proven for homology circles by 
Milnor and in the general case by Turaev.
\begin{theorem}[\cite{Mil61}\cite{Tur76}]\label{thm-ReidAlex} 
Let $M$ be a compact, oriented 3-manifolds,
$\varphi: H_1(M)\to\Z$ an epimorphism as above, 
$r_{\varphi}(M)$ its Reidemeister Torsion, and $\Delta_{\varphi}(M)$ its Alexander polynomial. 
\ben
\item If $\partial M\neq \emptyset$ then 
$\displaystyle r_{\varphi}(M)\,=\,\frac 1 {(t-1)}\Delta_{\varphi}(M)$
\item If $\partial M = \emptyset$ then 
$\displaystyle r_{\varphi}(M)\,=\,\frac 1 {(t-1)^2}\Delta_{\varphi}(M)$
\een
\end{theorem}

For a 3-manifold given by surgery along a framed link we will now give a
procedure to compute the Alexander polynomial (and thus also Reidemeister Torsion).

Let ${\cal Z}\sqcup{\cal L}\subset S^3$ be a framed link consisting of a 
 framed link ${\cal L}$ and 
 a curve ${\cal Z}$ which has
trivial linking number of all components of ${\cal L}$, i.e.,  
with ${\cal L}\cdot{\cal Z}=0$.  
We denote by $M_{{\cal Z},{\cal L}}^{\bullet}$ the manifold obtained by cutting
out a tubular neighborhood of ${\cal Z}$ and doing surgery along ${\cal L}$. 
Hence, $\partial  M_{{\cal Z},{\cal L}}^{\bullet}=S^1\times S^1$, with canonical
meridian and longitude (given by 0-framing).  
Also let $M_{{\cal Z},{\cal L}}$ be the closed manifold obtained 
by doing 0-surgery  along ${\cal Z}$ so that 
$M_{{\cal Z},{\cal L}}=M_{{\cal Z},{\cal L}}^{\bullet}\cup D^2\times S^1$. 
The special component $\cal Z$ defines an epimorphism 
$\varphi_{\cal Z}: H_1(M^{(\bullet)})\to\Z$, for example via intersection
numbers with a Seifert surface. We write  
$\Delta_{{\cal Z},{\cal L}}=\Delta_{\varphi_{\cal Z}}(M_{{\cal Z},{\cal L}})=
\Delta_{\varphi_{\cal Z}}(M_{{\cal Z},{\cal L}}^{\bullet})$ for the associated 
reduced Alexander polynomial,  which is the same in both cases. 

Consider a general Seifert surface $\Sigma^{\bullet}\subset S^3$ with 
$\partial \Sigma^{\bullet} ={\cal Z}$ and 
$\Sigma^{\bullet}\cap{\cal L}=\emptyset$.  By removing a neighborhood
of the surface we obtain a relative cobordism 
$C_{\Sigma}^{\bullet}=M_{{\cal Z},{\cal L}}^{\bullet}-\Sigma^{\bullet}
\times(-\epsilon,\epsilon)$ from $\Sigma^{\bullet}$ to itself. 
Similarly, $C_{\Sigma}=M_{{\cal Z},{\cal L}}-\Sigma
\times(-\epsilon,\epsilon)$, where $\Sigma$ is the closed capped off
surface $\Sigma^{\bullet}\cup D^2$. The cobordism $C_{\Sigma}$ is obtained
from $C_{\Sigma}^{\bullet}$ by gluing in a full cylinder $D^2\times[0,1]$.

Denote by 
$\psi^{(\bullet)}_{\pm}:\Sigma^{\bullet}_{\pm}\hookrightarrow  C_{\Sigma}$
the inclusion maps of the bounding surfaces, and by 
$$
A_{\pm}=H_1(\psi^{(\bullet)}_{\pm}):\,H_1(\Sigma)\to 
H_1(C_{\Sigma}^{(\bullet)})\to H_1^{\it free}(C_{\Sigma}^{(\bullet)})\;, 
$$
the maps on the free part of homology, where the free part is
$G^{\it free}=\frac G {Tors(G)}\,$. As 
$H_1(\widetilde M)\cong 
H_1^{\it free}(\widetilde M)\oplus\,Tors(H_1(M))\otimes\Z[t,t^{-1}]$
we will consider  the first elementary ideal for the free part, which differs only by
a factor of $|Tors(H_1(M))|$.

Suppose first that $C$ does not have interior homology. This means the $A_{\pm}$
can be presented as square matrices, and $A_+-tA_-$ is  a presentation matrix.  
Consequently 
$\Delta_{{\cal Z},{\cal L}}=\pm t^p det(A_+-tA_-)$. By some linear algebra
\cite{FroNic92}
this is the same as the Lefschetz polynomial 
$$
det(A_+-tA_-)\;=\;\sum_{k=0}^{2g}(-t)^{2g-k}
trace\Bigl((\ext k A_+)\circ \star^{-1}\circ
 (\ext {2g-k}A_-^*)\circ\star\Bigr)
$$
In  \cite{FroNic92} it is also shown that the expression inside the trace is
the same as ${\cal V}^{FN}(C_{\Sigma}^{\bullet})_k$ or 
${\cal V}^{FN}(C_{\Sigma})_k$ depending on context. Hence, we have (multiplying by a unit $(-t)^{-g}$)
that 
\begin{eqnarray}
\Delta_{{\cal Z},{\cal L}} & = &\sum_{k=0}^{2g}(-t)^{g-k} trace({\cal V}^{FN}(C_{\Sigma})_k)
\label{eq-leftrace}
\\
&=& trace((-t)^{-H}{\cal V}^{FN}(C_{\Sigma}))
\label{eq-leftrace2}
\\
&=& \sum_{j=1}\;[j]_{-t}\,trace({\cal V}^{(j)}(C_{\Sigma}))\;=\;
\sum_{j=1}\;[j]_{-t}\,\Delta_{{\cal Z},{\cal L}}^{(j)}\;,
\label{eq-leftrace3}
\end{eqnarray}
where $[n]_q=\frac{q^n-q^{-n}}{q-q^{-1}}$. In (\ref{eq-leftrace2}) we used the generator $H$ of
the $\SL(2,\R)$-Lefschetz action. Formula (\ref{eq-leftrace3}) is a consequence of the 
Hard-Lefschetz decomposition from (\ref{eq-Lefsch}). We call the invariant
$\Delta_{{\cal Z},{\cal L}}^{(j)}$ the {\em $j$-th Alexander Character} of the 
Alexander polynomial.

In case $C$ does have interior rational homology the dimension of 
$H_1^{\it free}(C_{\Sigma}^{(\bullet)})$ is bigger than $H_1(\Sigma)$ so that
$H_1(\widetilde M)$ has $\Z[t,t^{-1}]$ as a direct summand. Consequently, the Alexander
polynomial vanishes. At the same time ${\cal V}^{FN}(C_{\Sigma})$ is zero since it is
a non-semisimple TQFT. Hence, (\ref{eq-leftrace3}) holds for all cases.

Suppose that in our presentation ${\cal Y}\subset S^3$ is the unknot. In this case we 
can isotop the diagram ${\cal L}\sqcup {\cal Y}\subset S^3$ into the  form shown on the 
right side of Figure~\ref{fig-stan}. Specifically, we arrange 
it that the strands of one link component alternate orientations
as we go from left to right. By application of the connecting annulus 
moves, see for example \cite{Ker99}, we can modify the link further such that the resulting
tangle ${\cal T}$ in the indicated box is admissible without through pairs
 as described in the beginning of 
Section~5 or, again, \cite{Ker99}. There is a canonical Seifert surface $\Sigma_{\cal T}$
associated to a diagram as in Figure~\ref{fig-stan} obtained by surgering the disc
bounded by ${\cal Z}$ along the framed 
components of ${\cal L}$ emerging at the bottom side. By construction ${\cal T}$ is
then a tangle presentation of $C_{\Sigma_{\cal T}}$.

\begin{figure}[ht]
\begin{center}
\leavevmode
\epsfbox{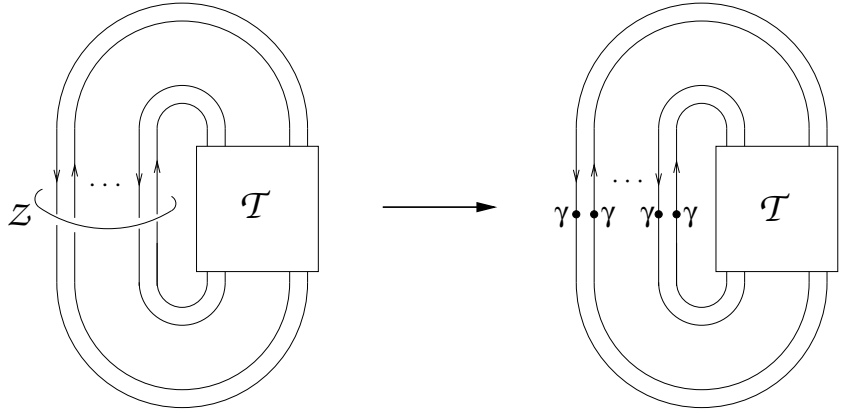}
\end{center}
\caption{Standard Presentation }\lbl{fig-stan}
\end{figure}

For the evaluation of this diagram it
 is convenient to introduce an extension of ${\cal N}$ over $\Z[t,t^{-1}]$,
given by $\Z[\gamma^{\pm 1}]\ltimes {\cal N}$. The extra generator $\gamma$ is
group like with $S(\gamma)=\gamma^{-1}$ and it acts on ${\cal N}$ by
$\gamma x\gamma^{-1}=t^Hx=t^{deg(x)}x$ for $x\in{\cal N}$ and $deg(x)$ the degree
for homogenous elements. 

In order to evaluate the diagram we apply the Hennings substitutions for crossing 
(\ref{fig-Dec-cross})  and rules (\ref{eq-elem-mv})   through (\ref{pic-Rclos}) 
 to the  ${\cal T}$ part to obtain a combination of ${\cal N}$-decorated 
arcs as in (\ref{pic-Rbot}) and (\ref{pic-Rtop}). Furthermore, we
remove the circle ${\cal Y}$ at the expense of introducing a $\gamma$-decoration
on each strand. The Hennings procedure is continued with the extended algebra over
$\Z[t,t^{-1}]$. It is easy to see that the elements that have to be evaluated against
the integral all lie in $\Z[t,t^{-1}]\otimes{\cal N}$ and that $\mu$ is cyclic also 
with respect to $\gamma$. Hence, the evaluation is well defined. 

\begin{lemma}
The evaluation procedure for a diagram as in Figure~\ref{fig-stan} yields the 
Alexander polynomial.
\end{lemma}

{\em Proof:} The standard evaluation of ${\cal T}$ yields a sum of diagrams with top and bottom 
arcs, where the $j$-th bottom arc is decorated by $b_j$ and the $j$-th top arc by $c_j$ 
as in (\ref{pic-Rbot}) and (\ref{pic-Rtop}). Hence, ${\cal V}_{\cal N}(C_{\Sigma})$ is the sum 
over all diagrams of linear maps
$\bigotimes_j^g (b_j\otimes \mu(S(\_)c_j))$.
  The extended evaluation yields closed curves, each of
which is decorated with four elements $b_j$, $c_j$, $\gamma$, and $\gamma^{-1}$. Using the 
antipodal sliding rule from (\ref{pic-Rclos}) we collect them at one side of a circle so that the 
evaluation becomes
$$
\mu(S^{-1}(b_j)\gamma c_j\gamma^{-1})\;=\;(-1)^{deg(b_j)}t^{deg(c_j)}\mu(S(b_j)c_j)
\;=\; (-t)^{-deg(b_j)}trace(b_j\otimes \mu(S(\_)c_j))\;.
$$
Note here, that $S^2(b_j)=(-1)^{deg(b_j)}$ and that the evaluation is non zero only if 
$deg(c_j)+deg(b_j)=0$. The sum (over all decorations) of the products (over $j$)
of these individual traces is thus just the trace of $(-t)^{-H}{\cal V}_{\cal N}(C_{\Sigma})$. 
Since this is (up to sign) identical with   $(-t)^{-H}{\cal V}^{FN}(C_{\Sigma})$ it follows
from (\ref{eq-leftrace2}) that the evaluation gives the Alexander polynomial.

\ep

The evaluation of a standard diagram can be described also more explicitly without the use
of the $\Z[\gamma]$ extension. Let ${\cal T}^{\#}: 2g\to 0$ be the diagram consisting of 
the tangle ${\cal T}: g\to g$ and the lower arcs. That is, 
${\cal T}=(1^g\otimes {\cal T}^{\#})\circ (X_g\otimes 1^g)$ and 
${\cal T}^{\#}= (X_g^{\dagger})\circ(1^g \otimes {\cal T})$, where $X_g^{\dagger}$ is the
upside down reflection of $X_g\,$. We define $A^{\gamma}\in
{\cal N}_0^{\otimes 2}\otimes\Z[t,t^{-1}]$ as 
\beq\label{eq-defAg}
A^{\gamma}\;=\;({\gamma}\otimes 1)A({\gamma^{-1}}\otimes 1)\;=\;
\frac 1 i \Bigl( \rho\otimes 1\,+\,1\otimes\rho\,-\,t^{-1}\tb\otimes\th\,+\,t
\th\otimes\tb\Bigr)\;\;\;.
\eeq
Moreover, we define 
$A_{\{g\}}^{\gamma}\in{\cal N}_0^{\otimes 2g}\otimes\Z[t,t^{-1}]$  from $A^{\gamma}$
as $A_{\{g\}}$ in (\ref{eq-Xn}) is defined from $A$ in (\ref{eq-defAtens}) and (\ref{eq-X1}), or,
equivalently, by
$$
A_g^{\gamma}\;=\;({\gamma}^{\otimes g}\otimes 1^{\otimes g})\circ A_{\{g\}}\circ 
({(\gamma^{-1})}^{\otimes g}\otimes 1^{\otimes g})\;. 
$$
This tensor is assigned to the upper arcs and the $\gamma$ elements in the standard diagram.
Hence, by the extended Hennings evaluation procedure the  Alexander polynomial is given by the
composition 
$$
\Delta_{{\cal Z},{\cal L}}\;=\;{\cal V}^{FN}({\cal T}^{\#}) (A^{\gamma}_g)\;,
$$
where we think of ${\cal V}^{FN}({\cal T}^{\#}):{\cal N}_0^{\otimes 2g}\to \Cc$ as being
naturally extended to a $\Z[t,t^{-1}]$-map from ${\cal N}_0^{\otimes 2g}\otimes\Z[t,t^{-1}]
\,\to\, \Cc[t,t^{-1}]\,$.

For further evaluation we use  Theorem~\ref{thm-solve} to write 
${\cal V}^{FN}({\cal T}^{\#})=\sum_{\nu}{\cal V}^{FN}(E_{\nu})$
as a combination of elementary tangles 
$E_{\nu}=(M_{k_1}\otimes\ldots M_{k_r}\otimes \epsilon^{\otimes N})\circ B\,$ so that
the Alexander polynomial is the sum of polynomials $E_{\nu}(A^{\gamma}_g)$. 
For the computation of these elementary polynomials it is convenient to use 
the following
graphical notation. As shown in (\ref{eq-symhopfeva}) we indicate the morphism $M_k$
by a tree with $k$ incoming branches. The  morphism  $X_1$ is drawn as an arc and
$X_g$ as  $g$ concentric arcs. 
\beq\label{eq-symhopfeva}
\epsfbox{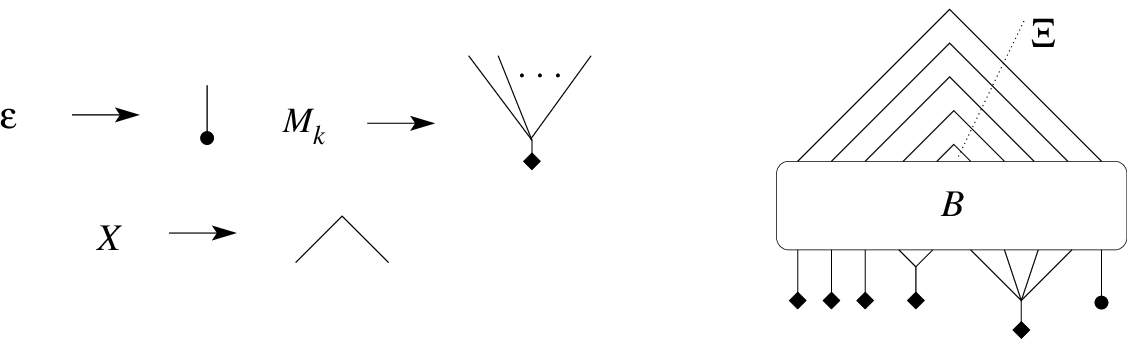}
\eeq
For $E=(M_1^{\otimes 3}\otimes M_2 \otimes M_4\otimes \epsilon)\circ B$ we obtain 
the composite shown on the right of (\ref{eq-symhopfeva}). Using relations 
$(\mu\otimes 1)A^{\gamma}=( 1\otimes \mu)A^{\gamma}= 1$, 
$(\varepsilon\otimes 1)A^{\gamma}=( 1\otimes\varepsilon)A^{\gamma}= \frac 1 i\rho$,
and $\mu(x\frac 1i \rho)=\varepsilon(x)$ we find the  graphical relations
depicted in (\ref{eq-elimeval}). 
\beq\label{eq-elimeval}
\epsfbox{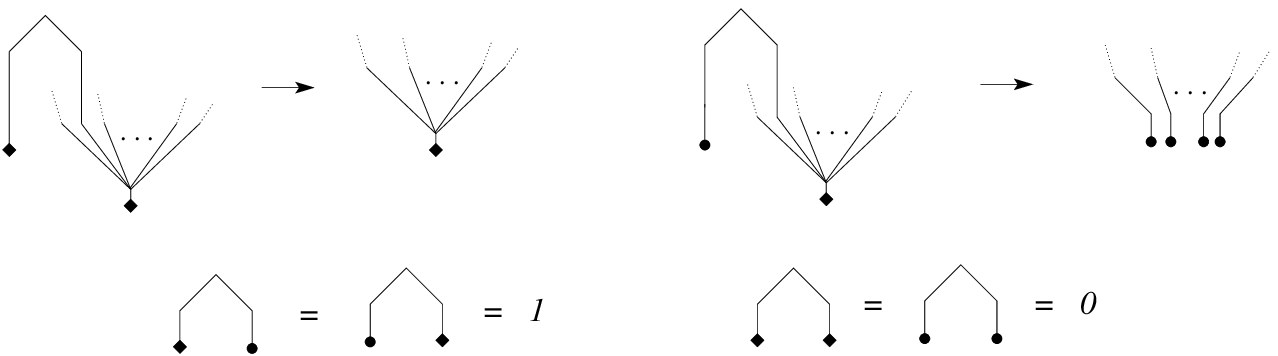}
\eeq
 Now, to each of the arcs 
the tensor $A^{\gamma}$ is associated containing the four terms $\rho\otimes 1$,
$1\otimes \rho$, $\tb\otimes\th$, and $\th\otimes\tb$ with coefficients of the form
$\pm it^m$. We represent the elementary
polynomial thus as a sum over all combinations of these terms, i.e., $4^g$ terms 
for $A^{\gamma}_{\{g\}}$. We indicate a combination in a diagram by drawing a line
with a down arrow for $\tb$, a line with an up arrow for $\th$, a line with arrows for
$\rho$ and a dashed line for $1$. Hence, (\ref{eq-defAg}) becomes the first line in 
(\ref{eq-termsarr}).
\beq\label{eq-termsarr}
\epsfbox{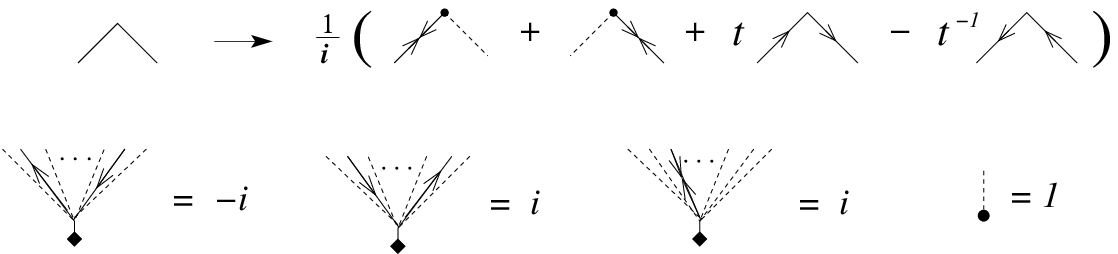}
\eeq
The tensors associated to the $M_k$ are non zero only in two cases. Namely, if one
element is $\th$, another $\tb$ and all other $1$, or if one element is $\rho$
and all others 1. In diagrams we obtain the evaluation rules as depicted. All
other configurations are evaluated to zero. 

For an elementary diagram let $N_x(=g)$ be the number of arcs at the top, 
$N_0$ the number of $\varepsilon$'s, and $N_k$ the number of $M_k$'s at 
 the bottom of the diagram for $k\geq 1$. Let us also call an elementary
diagram {\em reduced} if $N_0=N_1=0$.  We can now give the recipe for evaluating
elementary diagrams:

\begin{propos}\label{prop-HopfAlexeval}
\ 

\ben
\item We have the  relations \ \  \ \ 
$\displaystyle 
2N_x\,=\, N_0\,+\,\sum_{k\geq 1} kN_k\;,$  \ \  and  \ \ 
$\displaystyle  N_x\;=\;\sum_{k\geq 1} N_k\;.$
\item Every elementary diagram is zero or equivalent to a reduced one by application
of the moves in (\ref{eq-elimeval}). 
\item A reduced diagram is non zero only of $N_j=0$ for $j\geq 3$. That is, if
the diagram is of the form $D=M_2^{\otimes g}\circ B\circ X_g$.
\item A contributing reduced diagram $D=P_1\sqcup \ldots\sqcup P_n$ is the union of closed 
paths $P_j$, and the
polynomial $\Delta_D=\prod_j\Delta_{P_j}$ assigned to $D$ is the product 
of the polynomials assigned to the the components $P_j$. 
\item The polynomial associated to a connected component is 
$$
\Delta_P\;=\;2\,-(-1)^b\,(t^p+t^{-p})\;,
$$
where $p$ is the algebraic intersection number of the closed 
path $P$ with  a radial line segment $\Xi$ as in (\ref{eq-symhopfeva}), and $b$
is the total  number of half twists (or antipode insertions) in $B$. 
\een 
\end{propos}

{\em Proof:} {\em 1.}
In a diagram as in (\ref{eq-symhopfeva}) the number of strands entering from
the top is $2N_x$, two for each arc, and the number of strands entering from the bottom 
is $N_0+\sum_{k\geq 1}kN_k$. Obviously, 
both  numbers have to be equal. For an admissible configuration
of a contributing diagram we can also call weighted edges, where the dashed ones are
weighted 0, the ones with one arrow as 1, and those with double arrows as 2. The
top part of the diagram shows that the total weight has to be $2N_x$ since every 
admissible arc has weight 2. Also every tree has weight 2 and the $\epsilon$'s have
 weight 0 so that the total weight must also be given by $\sum_{k\geq 1} 2N_k$. 

{\em 2.} This is clear since every non-reduced one allows the application of a 
move that reduces the number of edges. 

{\em 3.} If we subtract twice the second identity in {\em i)} from the first we
find $0=N_0-N_1+N_3+2N_4+3N_5+\ldots\,$. In the reduced case with $N_0=N_1=0$ this 
implies $0=N_3=N_4=N_5=\ldots$ since these are all non negative integers. 

{\em 4.} Any graph where all vertices have valency 2 is the union of closed 
paths. Since we have a symmetric commutativity constraint we can untangle components
from each other and move them apart. The evaluation of disjoint unions of diagrams
is given by their products.

{\em 5.} There are four configurations that contribute to $\Delta_P$ for a 
closed path. Two if them are given by dashed lines alternating with double arrow
lines. This corresponds to paring factors $\frac 1 i\rho$ with integrals $\mu$ in
two different ways each evaluated as 1. Thus these two cases contribute the 2 in the 
expression. The other two configurations are given by two orientations of $P$ with
single arrows everywhere. For one given orientation we get from 
(\ref{eq-termsarr}) a factor $\frac 1 i t$ if
$P$ crosses $\Xi$ left to right and a factor $\frac 1 i (-t^{-1})$ if $P$ crosses
right to left. Thus the arcs yield a tensor $\pm (\frac 1 i)^g t^b(x_1\otimes\ldots\otimes x_{2g})$,
where each $x_i$ is either $\th$ or $\tb$. Application of $B$ yields a tensor
$\pm (\frac 1 i)^g t^b(y_1\otimes\ldots\otimes y_{g})$ where each $y_j$ is either $\th\otimes\tb$
or $\tb\otimes\th$ depending on which way the path runs through the $M_2$ piece. The 
pairwise multiplication thus yields the tensor $\pm t^b(\frac 1i\rho)^{\otimes g}$ and 
evaluation against $\mu$ the factor $\pm t^b$. For the opposite orientation the tensor
for the arcs is obtained by exchanging $t$ for $t^{-1}$ and multiplying a factor $(-1)^g$. 
The factor picked up by 
application of $B$ is unchanged, and in the evaluation against the $\mu$ we pick up a factor 
$(-1)^g$ because the orders of $\th$ and $\tb$ are exchanged canceling the one from the top. 
Hence the contribution for the opposite orientation is the same with $t$ and $t^{-1}$ exchanged.
Thus $\Delta_P=2\pm(t^b+t^{-b})$. The sign can be determined by evaluating the polynomial
at $t=1$. This is identical with the usual Hennings invariant of the 3-manifold given by 
surgery along a link associated to the connected diagram $P$ as follows. 

First choose 
over and under crossing for $P$ pushing it slightly outside the plane of projection into
a knot $P^*$. This knot is thickened to a band $N(P^*)$, which is  parallel to the
plane of projection except for half twists that are
introduced at the points where $B\subset P$ has antipodes
inserted. 

Consider the link $\partial N(P^*)$ given by the edges of the band. Generically this
link consists of  parallel strands that double cross as in Figure~\ref{fig-MC-tgl} at
simple crossings of $P^*$ and has $\Gamma$-diagram  also as in Figure~\ref{fig-MC-tgl}
for every half twist. We further modify this link  at some generic point in the band
by replacing  the parallel strands 
by a configuration with a connecting annulus as in the $\sigma$-Move of (\ref{eq-boundmove}).
We obtain a two component link ${\cal L}_P={\cal A}_P\sqcup{\cal C}_P$, where  
${\cal A}_P$ is the 0-framed annulus. The other part ${\cal C}_P$ bounds the disc obtained
by removing the small piece from the band where we applied the $\sigma$-Move and thus 
carries a natural framing.  We have by construction that 
$\Delta_P(1)=\pm\eta(M_{{\cal L}_P})$ with $\eta$ as in (\ref{eq-defeta}). For self intersection
numbers we clearly have ${\cal A}_P\cdot {\cal A}_P=0$ and 
${\cal C}_P\cdot {\cal C}_P=0$. For an even number of twists in the band $N(P^*)$ we
obtain also ${\cal A}_P\cdot {\cal C}_P=0$ and for an odd number of twists we have 
${\cal A}_P\cdot {\cal C}_P=\pm 2$. Hence $\eta(M_{{\cal L}_P})=0$ in the first case
and $\eta(M_{{\cal L}_P})=4$ in the second. 
\ep

Note, that the form of the $\Delta_P$ implies again the symmetry $\Delta(t)=\Delta(t^{-1})$
of the Alexander polynomial. In order to instill some confidence in our procedure let us
recalculate the familiar formula for the left-handed trefoil in this setting. Using the Fenn Rourke
move from  Figure~\ref{fig-FR} we  present the trefoil as an unknotted curve $\cal Z$ 
 in   a surgery diagram of Borromean rings as in (\ref{eq-trefoil}). 
\beq\label{eq-trefoil}
\epsfbox{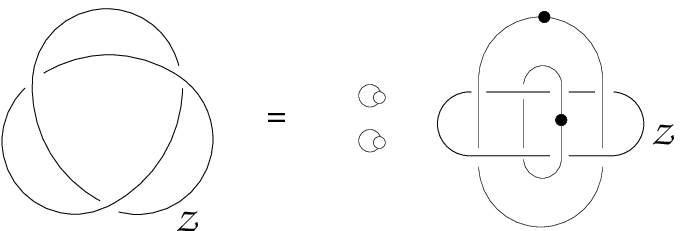}
\eeq
The standard form is obtained by moving ${\cal C}_1$ to the right off $\cal Z$ and 
letting ${\cal C}_2$ follow at the ends. The tangle ${\cal T}^{\#}$ is then as depicted on
the left of (\ref{eq-trefoiltgl}) below.  Using the framing moves from Figure~\ref{fig-frame} we 
expand it into elementary diagrams as on the right of  (\ref{eq-trefoiltgl}).
\beq\label{eq-trefoiltgl}
\epsfbox{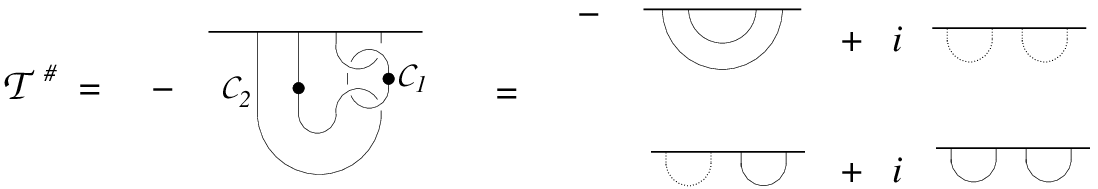}
\eeq
The translation into Hopf algebra diagrams and subsequently polynomials is indicated 
next in (\ref{eq-trefoilhopf}). 
\beq\label{eq-trefoilhopf}
\epsfbox{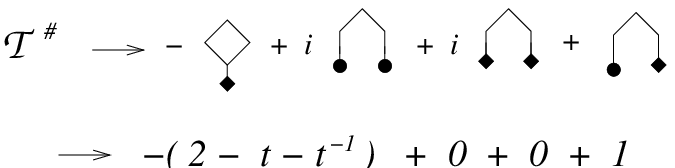}
\eeq
Thus the polynomial comes out to be  $t+t^{-1}-1$ as it had to be. 
The same calculation carries through
if we change the framings $f_j$ of the components ${\cal C}_j$ in (\ref{eq-trefoil}). The 
difference is the sign of the first summand, that is 
$\Delta_{\cal Z}=f_1 f_2(t+t^{-1}-2)+1$. Thus, if we flip both framings we obtain the right-handed
trefoil with the same polynomial. If we flip only one framing so that $f_1=-f_2$ we obtain one
of two figure-eight knots with polynomial $-t-t^{-1}+3$. Many other Alexander polynomials with
multiple twists as for example $(p,q,r)$-pretzel knots  can be computed quite conveniently
in this fashion using Fenn Rourke moves and the nilpotency of the ribbon element $v^k=1+k\rho$. 
Thus, our method proves to be 
 quite useful in the calculation of the Alexander Polynomial for knots although its primary 
application is the generalization to 3-manifolds.

We describe next a more systematic way to unknot the special strand $\cal Z$ in a general diagram 
more akin the traditional skein theory. The additional 
relations that allow us to put any diagram ${\cal L}\sqcup {\cal Z}$ into a standard form
are as follows.    

\begin{propos}\label{propos-Yskein}
 We have the following two skein relations for the special strand ${\cal Z}$
\beq\label{eq-Yskein}
\epsfbox{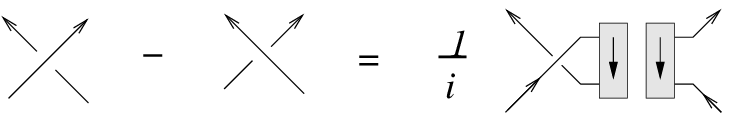}
\eeq
and
\beq\label{eq-Ycoupon}
\epsfbox{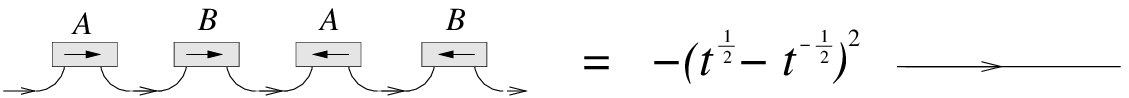}
\eeq
as well as the slide and cancellations moves analogous to (\ref{eq-1handlemove}),
and a vanishing property as in (\ref{eq-pairvanish}). 
 
These equivalences allow us to express the Alexander polynomial of
 any diagram ${\cal Y}\sqcup {\cal L}\subset S^3$ as a combination of the 
evaluations of diagrams in standard form. 
\end{propos}

{\em Proof:}  As before we change a self crossing of ${\cal Y}$ by sliding a 1-framed annulus 
${\cal A}$ over the crossing. Note, that we do not have to keep track of the framing of 
${\cal Y}$ as it is unchanged and by convention zero. 
  Using the orientation of ${\cal Y}$ we can  do this 
such that the linking numbers of ${\cal Y}$ and ${\cal A}$ remain zero. It is easy to see
that we can bring a diagram into the standard position as in  Figure~\ref{fig-stan} without
ever sliding a strand over the new component ${\cal A}$. The evaluation is obtained as the
weighted trace over the linear map associated by  ${\cal V}_{\cal N}$ to the
 cobordism represented by the tangle,
which contains ${\cal A}$. Inserting 
the relation from Figure~\ref{fig-frame} we see that this linear map, and hence
the associated polynomial, is the 
combination of the one for which ${\cal A}$ has been removed and the one for which the framing
of ${\cal A}$ has been shifted by one. In both cases the unknotting procedure can be reversed
so that we obtain the original pictures with ${\cal A}$ removed or its framing shifted by one. 
The situation in which ${\cal A}$ is removed corresponds to the opposite crossing. In the other
contribution we have a 0-framed annulus around the crossing which can be rewritten as an index-1
surgery represented by a pair of coupons. This yields (\ref{eq-Yskein}). 

The coupon combination in (\ref{eq-Ycoupon}) can be reexpressed by a tangle as in 
(\ref{eq-coupontorus}), can be isotoped into the position shown in (\ref{eq-Yext}). 
\beq\label{eq-Yext}
\epsfbox{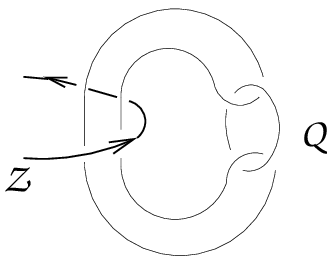}
\eeq
The extra tangle piece ${\cal Q}$ maps to the identity on a torus block. More precisely,
${\cal V}_{\cal N}({\cal Q}\sqcup{\cal T})=id_{{\cal N}_0}\otimes 
{\cal V}_{\cal N}({\cal T})$. The weighted traces thus differ by a factor
 $trace_{{\cal N}_0}((-t)^{-H})=-t+2-t^{-1}=-(t^{\frac 12}-t^{-\frac 12})^2$. 
\ep

For ordinary link and knot complements there are well known skein relations 
that uniquely characterize the Alexander-Conway polynomial of the knot,
see for example \cite{BurZie} Chapter 12.C. 

\begin{cor}
For ordinary knot complements (that is if ${\cal L}=\emptyset$) the relations 
Proposition~\ref{propos-Yskein} reduce to the ordinary Alexander-Conway skein
relations. 
\end{cor}

{\em Proof:} It is clear that with Proposition~\ref{propos-Yskein} we can resolve every diagram
into disjoint circles in the plane with coupons on them in exactly the same way as 
for the Alexander-Conway polynomial. The difference is that wherever we pick up a 
factor $(t^{\frac 12}-t^{-\frac 12})$ from the smoothening in the traditional calculus
we obtain a factor $\frac 1 i$ and a pair of coupons in our case, but all other 
numbers are the same.

Suppose now after resolving the crossings 
we have more than one circle. Since the strand $\cal Z$ has to run
though all of these components we must have coupons that are paired but on different 
circles. By (\ref{eq-pairvanish}) of Lemma~\ref{lm-1handlemoves} it follows that such a
configuration must vanish. In the Alexander-Conway calculus we also have the rule
that the link invariant for the unlinked union of an unknot with a non-trivial link
is zero. Hence we only need to compare the contributions that come from single circles.
If in the process of applying the skein relations we carried out $N$ smoothenings of
crossings the circle will carry $2N$ coupons. 

Next we claim that it is not possible to slide two paired coupons in adjacent position.
To this end note that  the coupons in the resolution of 
Proposition~\ref{propos-Yskein} stay all on one side of the
special strand. I.e., in the depicted orientation of $\cal Z$
the coupons are always on the left of $\cal Z$. Thus, 
if they become adjacent 
we would have a situation as in (\ref{eq-1handlecanc}) of  Lemma~\ref{lm-1handlemoves}. 
This is not possible since then $\cal Z$ would have at least two components. 
Thus the number $2N$  of coupons will remain the same under handle slides.

We next observe that a circle with edges that are labeled in pairs and subject to 
handle slides also occurs in the classification of compact, oriented surfaces via
their triangulations as in \cite{massey} Chapter 1. It is shown there that any
such configuration is under application of handle slides and cancellation
moves as in (\ref{eq-1handlecanc}) equivalent to a sequence of blocks as in
(\ref{eq-Ycoupon}). As before we may assume that all coupons lie on one side of 
the circle. In fact, as $\cal Z$ is connected we see from \cite{massey} that 
we can move to the configuration in standard block form without the use of 
cancellations.

Thus, we have $\frac N 2$ 4-coupon (torus) blocks as in  (\ref{eq-Ycoupon})
contributing a factor of 
$(-(t^{\frac 12}-t^{-\frac 12})^2)^{\frac N 2}=(i)^N(t^{\frac 12}-t^{-\frac 12})^N$.
Recall that in each resolution we also had a factor $\frac 1 i$ so that the 
total factor for the circle is just $(t^{\frac 12}-t^{-\frac 12})^N$ and $N$ is the
number of smoothenings. But  $(t^{\frac 12}-t^{-\frac 12})$ is precisely the factor
assigned to each smoothening by the usual Alexander-Conway calculus.

\ep 

Although we now have a systematic procedure for computing the Alexander polynomial
of a 3-manifold  it is often times efficient to use the skein relations leading up
to it directly. We illustrate this by  computing $\Delta_{{\cal C}_{k,l},{\cal Z}}$,
where ${\cal C}_{k,l}$ is the component depicted in (\ref{eq-exampleknot})
\beq\label{eq-exampleknot}
\epsfbox{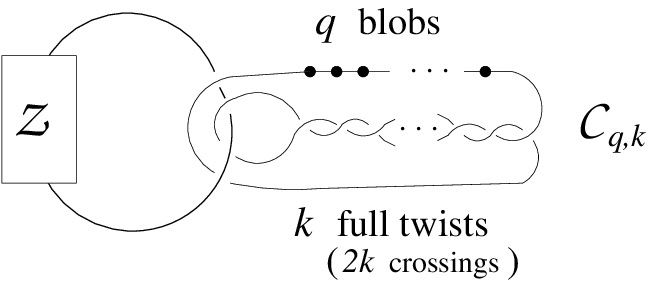}
\eeq
The two middle strands are twisted with each other $k$ times generating $2k$ crossings,
and we have $q$ full circles on the upper strand indicating shifts in the framing by -1.
The definition for $k<0$ or $q<0$ is given by choosing the opposite twistings. 

\begin{lemma} The Alexander Polynomial of $M_{{\cal C}_{k,l},{\cal Z}}$ is given  by 
the ordinary  Alexander polynomial of the knot as follows:
$$
\Delta_{{\cal C}_{k,l},{\cal Z}}\;=\;i(k(t+t^{-1})-q)\Delta_{{\cal Z}}
$$
\end{lemma}

{\em Proof:} We combine every twist with two circles so that we have $k$ twist
configurations as in Figure~\ref{fig-FR} and $l=q-2k$ remaining circles. Applying
the Fenn Rourke move to each of these we obtain a configuration in which we have a 
parallel instead of twisted pair of strands in the middle surrounded by $k$
annuli with an empty circle on them. In addition, we have $k$ separate annuli with 
full circles. Denote by $\Delta_{k,l}$ the associated Alexander Polynomial. For $k>0$
we choose one of the first annuli and apply the framing shift relation
(\ref{fig-frame}) to the empty circle on it. In the second contribution we omit the
dotted line so that we obtain the same configuration with one less annulus around
the double strands. The factor $i$ in  (\ref{fig-frame}) is canceled against one
of the separate annuli with a full circle so that the second contribution is exactly 
$\Delta_{l,k-1}$. In the first contribution we have a 0-framed annulus which, by 
Figure~\ref{fig-modif}, can be turned into a pair of coupons. The other $k-1$ coupons
can thus be slid off and canceled against $k-1$ annuli with full circles. Moreover,
the remaining $l$ full circles on the upper strand can be removed since inserting a
dotted line leaves two isolated coupons, which yields zero. The resulting 
configuration is the knot $\cal Z$ with a tangle piece ${\cal Q}$ as in 
(\ref{eq-Yext}), contributing an
extra factor $-(t^{\frac 12}-t^{-\frac 12})^2$, and an 
extra annulus with full circle with a factor $-i$. We thus obtain the recursion
relation  $P_{k,l}\,=\,i(t^{\frac 12}-t^{-\frac 12})^2\,+\,P_{k-1,l}$ so that
$P_{k,l}\,=\,ik(t^{\frac 12}-t^{-\frac 12})^2\Delta_{\cal Z}\,+\,P_{0,l}$. 
But the configuration
for $k=0$ is the separate union of $\cal Z$ and an annulus with $l$ full circles.
The latter yields a factor $-il$ so that 
$P_{k,l}\,=\,i(k(t^{\frac 12}-t^{-\frac 12})^2-l)\Delta_{\cal Z}$, which computes
to the desired formula.

\ep

\head{12. Lefschetz compatible Hopf algebra structures on $H^*(J(\Sigma))$}\lbl{S12}

It is easy to see that the natural ring structure on the cohomology 
$H^*(J(\Sigma))\cong\ext * H_1(\Sigma)$ is not compatible with the $\SL(2,\R)$
Lefschetz action as described in Section~10. For example $E(x\wedge y)=x\wedge y\wedge \omega$
but $(Ex)\wedge y +x\wedge (Ey) =2x\wedge y\wedge \omega$. The isomorphism with 
${\cal N}_0^{\otimes g}$ however induces another multiplication structure compatible 
with the $\SL(2,\R)$ action. In this section we will describe it explicitly. 

The $\Z/2$-graded Hopf algebra structure on ${\cal N}_0$ given in 
Lemma~\ref{lm-super-Hopf}  extends to a $\Z/2$-graded Hopf algebra 
structure ${\cal H}_{\cal N}$ on ${\cal N}_0^{\otimes g}$ with 
$$
(x_1\otimes \ldots\otimes x_g)(y_1\otimes \ldots\otimes y_g)
\;\;=\;\;(-1)^{\sum_{i<j} d(x_j)d(y_i)}\,x_1y_1\otimes \ldots\otimes x_gy_g\;.
$$
The formula for ${\bf \Delta}$ is the dual analog.
The precise form of ${\cal H}_{\cal N}$ is given as follows:

\begin{lemma}\lbl{lm-StrucN}  For a choice of basis of $\R^g$ there 
is  a natural isomorphism of Hopf algebras
$$
\varrho\;:\;\;\ext *(\E\otimes \R^g)\;\;\isto\;\;{\cal N}_0^{\otimes g}\;\;
$$

so that \ \ \  ${\rm Aut}({\cal N}_0^{\otimes g}, 
{\cal H}_{\cal N})\;\cong\;{\rm GL}(\E\otimes \R^g)$. 
\end{lemma}

{\em Proof:} Let $\{e_j\}$ be a basis of $\R^g$. The generating set of primitive vectors
of  $\ext *(\E\otimes \R^g)$ is given by $\E\otimes \R^g$. On this subspace we set 
$\varrho(w\otimes e_j)=1\otimes\ldots 1\otimes w \otimes 1\ldots\otimes 1$, with $w$ in
$j$-th position. We easily see  that the vectors in $\varrho(\E\otimes \R^g)$ form 
again a generating set of anticommuting, primitive vectors of ${\cal N}_0^{\otimes g}$ 
so that $\varrho$ extends to a Hopf algebra epimorphism. Equality of dimensions thus implies
that $\varrho$ is an isomorphism.\ep

 The canonical $\SL(2,\R)$-action on ${\cal N}_0^{\otimes g}$ is still compatible 
with ${\cal H}_{\cal N}$ since it preserves the degrees and factors. Under the isomorphism 
in Lemma~\ref{lm-StrucN} it is readily identified as the $\SL(2,\R)$-action on the $\E$-factor. 
The remaining action on the $\R^g$-part can be understood geometrically. Specifically, 
$\Sp(2g,\Z)$ acts on ${\cal N}_0^{\otimes g}$ since the ${\cal V}$-representation of the 
mapping class group factors through a the symplectic group with representation 
${\cal V}^{Sp}:\,\Sp(2g,\Z)\to {\rm GL}({\cal N}_0^{\otimes g})\,:\,[\psi]\mapsto 
 {\cal V}^{Sp}([\psi]):={\cal V}({\bf I}_{\psi})$. 
For a given decomposition into Lagrangian subspaces we denote the standard inclusion
\beq\lbl{eq-defkappa}
\kappa\,:\;
\SL(g,\Z)\,\into {\rm }GL(g,\Z)\,\into \,\Sp(2g,\Z)\;:
\;\;A\,\mapsto\, \kappa(A):=\,A\oplus (A^{-1})^T\;.
\eeq

\begin{lemma}\lbl{lm-slgz-invar}
The action of $\SL(g,\Z)$ on ${\cal N}_0^{\otimes g}$ induced by ${\cal V}^{Sp}\circ\kappa$
is compatible with ${\cal H}_{\cal N}$, and under the isomorphism  $\varrho$ from 
Lemma~\ref{lm-StrucN} it is identical with the $\SL(g,\Z)$-action on $\R^g$ for the given basis. 
In particular,  it commutes with the $\SL(2,\R)$-action so that we have the following
natural inclusion of the Howe pairs:
$$
\SL(2,\R)\times \SL(g,\Z) \;\;\subset\;\; 
{\rm GL}(\E\otimes \R^g)={\rm Aut}({\cal N}_0^{\otimes g},{\cal H}_{\cal N} )\;\;. 
$$
\end{lemma}

{\em Proof:} Consider the  elements
$P_j:=S_j\circ D_j^{-1}\circ S_j^{-1}$ and 
$Q_j:=S_{j+1}\circ D_j^{-1}\circ S_{j+1}^{-1}$  of $\Gamma_{g,1}$. 
From (\ref{eq-D-hom}) and (\ref{eq-S-hom}) we compute the homological
action as $[R_j]=\kappa(I_g+E_{j+1, j})$ and $[Q_j]=\kappa(I_g+E_{j, j+1})$,
with conventions again as in \cite{GooWal98}. 
The matrices $I_g+E_{j+1, j}$ and $I_g+E_{j, j+1}$ generate $\SL(g,\Z)$, and hence 
$[P_j]$ and $[Q_j]$ generate $\kappa(\SL(g,\Z))\subset \Sp(2g,\Z)\,$. 
The actions of ${\cal V}({\bf I}_{P_j})$ and 
${\cal V}({\bf I}_{Q_j})$ on ${\cal N}_0^{\otimes g}$ are given by 
placing the maps $\Pp:=(\Ss\otimes 1)\D^{-1}(\Ss^{-1}\otimes 1)$ and
$\Q:=(1 \otimes\Ss )\D^{-1}(1 \otimes\Ss^{-1})$ in the $j$-th and $j+1$-st tensor 
positions. In order to show that the actions of $P_j$ and $Q_j$ on ${\cal N}_0^{\otimes g}$
yield Hopf algebra automorphisms it thus suffices to prove this for the maps 
$\Pp$ and $\Q$ in the case $g=2$. From the tangle presentations  we find 
identities ${\bf I}_{Q_1}= ({\bf M}\otimes 1)\circ (1\otimes {\bf \Delta})$
and  ${\bf I}_{P_1}= (1\otimes {\bf M})\circ({\bf \Delta}\otimes 1)$. 
It follows that $\Pp(x\otimes y)=\Delta_0(x)(1\otimes y)$ and 
$\Q(x\otimes y)=(x\otimes 1)\Delta_0(y)$. The fact that these are Hopf automorphisms
on ${\cal N}_0\otimes {\cal N}_0$ can be verified by direct computations. For the 
multiplication this amounts to verification of equations such as 
$\Delta(w)1\otimes v=-1\otimes v\Delta(w), \forall v,w\in\E\,$, and for the 
comultiplication we use the
fact that ${\cal N}_0$ is self dual. 

From the above identities we have that  
${\cal V}({\bf I}_{Q_1})=(M_0\otimes 1)\circ (1\otimes  \Delta_0)$ so that 
${\cal V}({\bf I}_{Q_j})$ is given on a monomial by taking the coproduct of the element in
$(j+1)$-st position, multiplying the first factor of that to the element in $j$-th position 
and placing the second factor into $(j+1)$-st position. We readily infer for every $w\in\E$ that
${\cal V}({\bf I}_{Q_j})(\varrho(w\otimes e_k))\,=\,
\varrho(w\otimes e_k+\delta_{j+1,k}w\otimes e_j)\,=\,\varrho(w\otimes (I_g+E_{j+1,j})e_k)$.
The analogous relation holds for $[P_j]$ so that 
$$
{\cal V}^{Sp}(\kappa(A))(w\otimes x)\quad=\quad w\otimes (Ax)\;\;\;
\forall A\in \SL(g,\Z). 
$$
This is precisely the claim made in Lemma~\ref{lm-slgz-invar}. 
\ep

The structure ${\cal H}_{\cal N}$ is mapped by the isomorphism $\xi_g$ from (\ref{eq-defxi}) 
to a $\Z/2$-graded Hopf algebra structure ${\cal H}_{\Lambda}$ on $H^*(J(\Sigma_g))$. 
A-priori the isomorphism $\xi_g$ and thus also  ${\cal H}_{\Lambda}$ depend on the choice of
a basis of $H_1(\Sigma_g)$. However, the $\SL(g,\Z)$-invariance determined in Lemma~\ref{lm-slgz-invar}
translates to the $\SL(g,\Z)$-invariance of ${\cal H}_{\Lambda}$, where 
$\kappa(\SL(g,\Z))\subset \Sp(2g,\Z)$ acts in the canonical way on $H^*(J(\Sigma_g))$. 
Hence, ${\cal H}_{\Lambda}$ only depends on the oriented subspaces 
$\Lambda=\lz [a_1],\ldots,[a_g]\rz\subset H_1(\Sigma_g,\Z)$
and $\Lambda^*=\lz [b_1],\ldots,[b_g]\rz\subset H_1(\Sigma_g,\Z)$, 
but not the specific choice of basis within them. 
The orientations can be given by volume forms $\omega_{\Lambda}:=[a_1]\wedge\ldots\wedge[a_g]$
and $\omega_{\Lambda^*}:=[b_1]\wedge\ldots\wedge[b_g]$. 
The primitive elements $\varrho(\th\otimes e_j)$ and $\varrho(\tb\otimes e_j)$
of ${\cal N}_g^{\otimes g}$  are mapped by $\xi_g$ to
\beq\lbl{eq-primit}
\pm[a_j]\wedge\omega_{\Lambda^*}\in\ext {g+1}H_1(\Sigma_g)\qquad\mbox{and}\qquad
\pm i^*_{z_j}(\omega_{\Lambda^*})\in\ext {g-1}H_1(\Sigma_g)
\eeq
respectively, where $[a_j]\in H_1(\Sigma_g)$ and $z_j\in H^1(\Sigma_g)$, with $z_j([b_j])=1$ and
$z_j([x])=0$ on all other basis vectors. We also have $\xi_g(1)=\omega_{\Lambda^*}$ and 
$\xi_g(\rho^{\otimes g})=\omega_{\Lambda}$. 

This completes the proof of Theorem~\ref{thm-structure}. \ep
\medskip

In the remainder of this section we give a more explicit description of the 
structure ${\cal H}_{\Lambda}$ on $H^*(J(\Sigma_g))$, and relate it to an involution,
$\tau$,
 on $H^*(J(\Sigma_g))$, which acts as identity on the $\Lambda$-factor and, modulo signs,
as a Hodge star on the opposite $\Lambda^*$-factor. 

The product $\diamond$ on $(H^*(J(\Sigma_g)),{\cal H}_{\Lambda})$ 
is given on a genus one block, $\ext * \lz[a],[b]\rz$, as follows:
\beq\lbl{eq-diamond}
\begin{array}{c}
\mbox{Table for}\\
\\
\quad u\,\diamond\, t\;:=\;\phi(\phi^{-1}(u)\phi^{-1}(t))\\
\end{array}
\quad\quad
\begin{tabular}{ | c || c|c|c|c|}
\hline
${}_u \,\backslash \,{}^t$\ \  
&\ \  1\ \  &\ \  $[a]$\ \   &\ \  $[b]$\ \  &\ \  $[a]\wedge [b]$ \ \ \\
\hline\hline
1 & 0 &  0 & 1 & $[a]$ \\
\hline 
$[a]$ & 0 & 0 & $a$ & 0 \\
\hline 
$[b]$ & 1 & $[a]$ & $[b]$ & $[a]\wedge [b]$ \\
\hline 
$[a]\wedge [b]$ & $-[a]$ & 0 &  $[a]\wedge [b]$ & 0\\
\hline 
\end{tabular}
\qquad\qquad
\eeq
It extends to $\ext * H_1(\Sigma_g)$ via the formula 
\beq\lbl{eq-totprod}
(u_1\wedge\ldots\wedge u_g)\diamond (t_1\wedge\ldots\wedge t_g)
\quad=\quad (-1)^{\sum_{i<j}d_il_j}\,\,(u_1\diamond t_1)
\wedge \ldots\wedge(u_g\diamond t_g)\;,\quad
\eeq
where $u_i, t_i\in\ext *\lz [a_i],[b_i]\rz$, $d_i=1-deg(u_i)$ and  $l_j=1-deg(t_j)$. 
In particular, we have $\,u\diamond t=(-1)^{d l}\,t\diamond u\,$, 
with $d=\sum_id_i=g-deg(u)$ and $l=\sum_il_i=g-deg(t)$, which reflects the
$\Z/2$-commutativity of $H^*(J(\Sigma_g))$.

The product structure and another proof of Lemma~\ref{lm-slgz-invar} can be also
found from an involution, $\tau$, defined as follows:
 
Every cohomology
class $x\in H^*(J(\Sigma_g))$ is uniquely written as $x=\alpha\wedge\beta$, where
$\alpha\in\ext *\Lambda$ and $\beta\in \ext *\Lambda^*$. For $x$ in this form the map $\tau$
is uniquely determined by the  relations  
\beq\lbl{eq-deftau}
\tau(\alpha\wedge\beta)=\alpha\wedge\tau(\beta)\quad\mbox{and}\quad
\tau(b_1^{\epsilon_1}\wedge\ldots\wedge b_g^{\epsilon_g})=
b_1^{1-\epsilon_1}\wedge\ldots\wedge b_g^{1-\epsilon_g}\;.
\eeq
From the formulae in (\ref{eq-diamond}) and (\ref{eq-totprod}) we find that $\tau^2=1\,$, and
\beq\lbl{eq-tauprop}
\tau(u\diamond t)\;\;=\;\;\tau(t)\wedge\tau(u)\;,\quad\quad\;
\eeq
and that $\tau$ maps $\ext *\Lambda$ as well as  $\ext *\Lambda^*$ to itself. It is clear from  
(\ref{eq-deftau}) and (\ref{eq-tauprop}) that $\SL(g,\Z)$-variance
of $\diamond$ on $H^*(J(\Sigma_g))$ is equivalent to $\SL(g,\Z)$-variance
of $\diamond$ on $\ext * \Lambda^*$. Now,  for any $A\in \SL(\Lambda^*)$ the following
identity holds:
\beq\lbl{eq-tauAid}
\tau\circ(\ext * A)\circ \tau\;\;=\;\;\ext* \iota(A)\;,
\eeq
where $\iota$ is the involution on $\SL(\Lambda^*)$ defined by 
$$
\iota(A)\;:=\;D\circ(A^{-1})^T\circ D\,,\qquad\mbox{with}\;\; D[b_j]=(-1)^j[b_j]\;.
$$
This can be proven either by considering again generators of $\SL(\Lambda^*)$,
or by applying the generalized Leibniz formula for the expansion of the determinant
of a $g\times g$-matrix into products of determinants of $k\times k$ and 
$(g-k)\times(g-k)$-submatrices. See also Lemma~5.2 in \cite{FroNic94}. (\ref{eq-tauprop})
together with  (\ref{eq-tauAid}) implies now that $\diamond$ depends only on the
decomposition $H_1(\Sigma_g,\Z)=\Lambda\oplus\Lambda^*$. 

In summary, we have the following isomorphism of $\Z/2$-graded Hopf algebras:
$$
\tau':=\ext *D\circ\tau\,:\,(H^*(J(\Sigma_g)),{\cal H}_{\Lambda})\;\isto\;
(H^*(J(\Sigma_g)),{\cal H}_{ext})\;\;,
$$
The Howe pair 
$\SL(2,\R)\times \SL(g,\R)\subset GL(H_1(\Sigma_g))=
{\rm Aut}(H^*(J(\Sigma_g)),{\cal H}_{ext})$, with 
$H_1(\Sigma_g)=\E\otimes \Lambda$, is conjugated by $\tau'$ to the pair 
$\SL(2,\R)_{Lefsch.}\times \kappa(\SL(g,\R))\subset 
{\rm Aut}(H^*(J(\Sigma_g)),{\cal H}_{\Lambda})$.  
\medskip

\newpage

\head{13. More Examples of Homological TQFT's and Open Questions}\lbl{S13}


\emptystuff{
\paragraph{A. Homology TQFT's over $\Z/r$ and cut numbers.}\

Although the TQFT's  of Reshetikhin and Turaev are semisimple and 
non trivial on the Torelli groups they contain homological TQFT's in an 
indirect manner. Specifically, if we consider the TQFT for $U_q({\mathfrak s}{\mathfrak l}_2)$
for $q$ a primitive $r$-th root of unity and $r$ is an odd prime 
Gilmer \cite{Gil01} shows that it can be defined essentially as a 
theory ${\cal W}^*_r$
over the ring over cyclotomic integers $\Z[q]$. This generalizes
the integrality results in \cite{MasRob97} and \cite{Mur94} for the 
invariants of closed manifolds. 

Of particular interest are expansions in $(q-1)$ which on the level of
invariants of closed manifolds lead to the Ohtsuki invariants in $\Z/r$ 
\cite{Oh96} which in 0-th order coincides with the invariant from (\ref{eq-defeta})
and in next order is identical to the Casson invariant \cite{Mur94}. 

A candidate for a useful homological TQFT is the lowest order of the TQFT
over the cyclotomic integers. It is given by the extending the
 trace function $\Z[q]\to\Z/r$ to a transformation  
${\cal W}^*_r\,\to\,{\cal W}_r$, where ${\cal W}_r$ is the respective TQFT
defined over the finite field $\Z/r$.

In \cite{Kerr=5} we consider the first non-trivial prime $r=5$ and
find an explicit basis for ${\cal W}^*_5$ and hence a description of
${\cal W}_5$. We find that the Torelli group is not entirely in the
kernel of ${\cal W}_5$ but factors through the Johnson homomorphism.
It does, however, contain a sub-TQFT $ {\cal U}_5\subset {\cal W}_5$
which is homological, meaning does not see the Torelli group, 
such that also the quotient TQFT 
${\cal Q}_5={\cal W}_5/ {\cal U}_5$ is homological. 

Explicit computations strongly suggest that these homological TQFT's are  
also of the general form (\ref{eq-defpolyTQFT}). More precisely,
we define the following (linear)
polynomials in $\Z[x_0,x_1,\ldots]$. 
$$
Q\,=\,x_0-x_8+x_{10}-x_{18}+x_{20}-x_{28}+\ldots
$$
$$\mbox{and} 
\qquad\qquad 
U\,=\,x_3-x_5+x_{13}-x_{15}+x_{23}-x_{25}+\ldots \;\;. 
\qquad\qquad \qquad
$$
Explicit verification for low genera and comparison of dimensions leads
us to the following. 
\begin{conjecture}[see\cite{Kerr=5}]\label{conj-5}
 We have isomorphisms of TQFT's defined 
over $\Z/5$:
$$
{\cal Q}_5\,\,\cong\,\,{\cal V}^{(Q)}
\qquad\qquad \mbox{and}\qquad \qquad 
{\cal U}_5\,\,\cong\,\,{\cal V}^{(U)}
$$
\end{conjecture}
Note that the polynomials $Q$ and $U$ also have negative integers so that
we need to make sense of subtracting vector spaces or TQFT's. To this end note that
the $\Sp(2g,\Z)$-representations $W_{g,j}$ are irreducible over $\Z$ but
become reducible if we take them, for a given basis, over $\Z/5$. For instance the 
$\Z/5$-reduction of $W_{6,3}$ contains a subrepresentation isomorphic to
the $\Z/5$-reduction of  $W_{6,5}$. This explains the meaning of the 
difference $x_3-x_5$ in the expression for $U$.

A fascinating topological application is the determination of so called
{\em cut numbers}, which is investigated in joint work with Gilmer \cite{GilKer}. 
Let us denote by $cut(M)$ the maximal number rank $n$ of a (non-abelian) free group
$F_n=\Z*\ldots*\Z$ such that there is an epimorphism $\phi:\pi_1(M)\to F_n$. 
This is also the maximal number of surfaces that can be removed from $M$ without
disconnecting the manifold. For a given epimorphism $\varphi:H_1(M)\to \Z$ we
also define the relative cut number $cut(M,\varphi)$ as the maximal $n$ such that
there is an epimorhpism $\phi:\pi_1(M)\to F_n$ which factors through $\varphi$,
meaning there is a map $\tau:F_m\to\Z$ such that $\varphi=\tau\circ\phi$. This 
counts non separating surfaces with the constraint that one represents $\varphi$, \cite{GilKer}. 

Clearly we  have $cut(M)\geq cut(M,\varphi)\geq 1$ if defined. Aside from these 
constraints the
absolute and relative cut number are independent.  For example
let $M=S^1\times \Sigma_g$ with canonical projection $\pi:M\to S^1$. Then 
$cut(M,\pi_*)=1$ but $cut(M)\geq g$.

As we remarked in the beginning of Section~11 an additional non-separating
surface in the cut cobordism $C$ used to define the Alexander polynomial
implies ${\cal V}^{FN}(C)=0$ by non-semisimplicity. Thus for a 3-manifold $M$ 
with epimorphism $\varphi:H_1(M)\to\Z$  as before we have the implication: 
\beq\label{eq-cutrel}
\mbox{If}\qquad\quad cut(M,\varphi) >1\qquad\quad{\rm then}
\qquad\quad \Delta_{M,\varphi}\,=\,0.
\eeq
In \cite{GilKer} we manage to obtain  a criterium on the bare cut number
{\em independent} of  a choice of $\varphi$:  
\beq\label{eq-cut}
\mbox{If}\qquad\quad cut(M) >1\qquad\quad{\rm then}
\qquad\quad trace({\cal W}_5(C))\,\equiv\,0\;{\rm mod}\,5.
\eeq
Note that the expression on the right only depends on the homological
functors $\cal Q$ and $\cal U$. It turns out that  under the assumption of 
Conjecture~\ref{conj-5} the respective traces are easily computed form the Alexander polynomial.
In fact, under this assumption,
 the trace expression
in (\ref{eq-cut}), which is just the sum of the traces for $\cal Q$ and $\cal U$,  
  comes out to be equal to 
the unique coefficient $T_{M,\varphi}$
of the Alexander polynomial when written as follows. 
$$
\Delta_{M,\varphi}(q)\;=\;T_{M,\varphi}\;+\;B(q+q^{-1})\;\;\;\in\;\;\;\Z[q]\qquad\mbox{with}\quad 
T_{M,\varphi}, B\in \Z\;\;.
$$
The contrapositive of (\ref{eq-cut}) under the assumption of the conjecture thus becomes
\beq\label{eq-contracut}
\mbox{If}\qquad\quad T_{M,\varphi}\neq 0 \,\, {\rm mod}\,5\qquad\quad{\rm then}
\qquad\quad cut(M)=1\;.\qquad\mbox{(for any choice of $\varphi$)}
\eeq
See \cite{GilKer} for more details and applications. 
}


\paragraph{A. Relations to Gauge Theories and the TQFT-Ring \Qq Generated by  ${\cal V}$.}\ 

We begin by collecting the ingredients that imply Theorem~\ref{thm-relations}.
The first identity (\ref{eq-AlexChar}) has already been computed in 
(\ref{eq-leftrace3}). 

The invariants $I^{DC}$ and $I^{SW}_d$ are obtained by Donaldson in \cite{Don99} 
from TQFT's  ${\cal V}^{DC}$ and ${\cal V}^{SW}_d$ respectively. For both TQFT's the 
vector spaces associated to a surface $\Sigma$ are the homologies of  natural 
moduli spaces. In the case of ${\cal V}^{DC}$ this is the  moduli space ${\cal M}(\Sigma)$ 
of flat connections on a non-trivial $SO(3)$ bundle. For  ${\cal V}^{SW}_d$
the moduli space of solutions to certain vortex equations is considered, which turn
identified with the  symmetric products of the surface. The action of the mapping
class group on the resulting homologies also factors through the symplectic groups (with 
the familiar ${\mathbb F}_2$-ambiguity). Donaldson thus derives the following 
isomorphisms between $\Sp(2g,\Z)$-modules.
\beq\label{eq-Donaldson}
\begin{array}{rcccl}
\displaystyle {\cal V}^{DC}(\Sigma_g)&\;\;=\;\;&H_*\bigl({\cal M}(\Sigma_g)\bigr)
&\;\;\cong\;\;&\displaystyle 
\bigoplus_{j=0}^g\,{\mathbb Q}^{j^2}\otimes \ext {g-j}H_1(\Sigma_g)\;\;\mbox{and}\\
\displaystyle {\cal V}^{SW}_d(\Sigma_g)&\;\;=\;\;&H_*\bigl(Sym^k(\Sigma_g)\bigr)
&\;\;\cong\;\;&\displaystyle 
\bigoplus_{j=1}^{g-d}\,{\mathbb Q}^{j}\otimes \ext {g-d-j}H_1(\Sigma_g)\;,\\
\end{array}
\eeq
where $k=g-1-d$ is the degree the holomorphic line bundle of which the vortex solutions
are sections. The $\Sp(2g,\Z)$-representations can be further identified with the 
irreducible parts, which takes in our notation the form
\beq\label{eq-extexpand}
\ext {g-j}H_1(\Sigma_g)\;\;=\;\;{\cal V}^{(j+1)}(\Sigma)\;\oplus\;
{\cal V}^{(j+3)}(\Sigma)\;\oplus\;{\cal V}^{(j+5)}(\Sigma)\;\oplus\;\ldots
\eeq
Inserting (\ref{eq-extexpand}) into the isomorphisms in (\ref{eq-Donaldson}) we
find  that the ${\cal V}^{DC}(\Sigma_g)$ and ${\cal V}^{SW}(\Sigma_g)$ are direct
sums of the ${\cal V}^{(j)}(\Sigma_g)$ with $g$-independent 
multiplicities given by precisely the
non-negative coefficients in (\ref{eq-DCChar}) and (\ref{eq-SWChar}). In Chapter~5 of
\cite{Don99} 
Donaldson exploits this fact to show that the decomposition thus extends to the 
entire TQFT's, meaning that cobordisms act trivially on the mulitiplicity spaces and
have block-wise actions on the  ${\cal V}^{(j)}$ components equivalent to those in 
the ${\cal V}^{FN}$ case. Summarily, we have the following isomorphisms 
{\em of TQFT's}.
\beq\label{eq-TQFTidentities}
\begin{array}{rcl}
\displaystyle {\cal V}^{DC}
&\;\;\cong\;\;&\displaystyle 
\bigoplus_{j\geq 2}\,{\mathbb Q}^{{j+1}\choose 3}\otimes {\cal V}^{(j)}\;\\
\displaystyle {\cal V}^{SW}_d
&\;\;\cong\;\;&\displaystyle 
\bigoplus_{j\geq d+2}\,{\mathbb Q}^{\left[\!\!{\left[
\Bigl({\frac{j-d}2}\Bigr)^2\right]}\!\!\right]}\otimes {\cal V}^{(j)}\;.\\
\end{array}
\eeq
Identities (\ref{eq-DCChar}) and (\ref{eq-SWChar}) are now immediate. For the
last equation  (\ref{eq-LesChar}) in Theorem~\ref{thm-relations}
we refer to \cite{KerKyoto}.

\ep

In an effort  to find new knot invariants Frohman and Nicas generalized their
approach in  \cite{FroNic94} to higher rank Lie algebras. They construct a TQFT
${\cal V}^{PSU(n)}_k$, 
whose vector spaces are
given as intersection homology groups of certain restricted moduli spaces 
of $PSU(n)$-representations and derive from these by similar trace formulae 
 invariants $\lambda_{n,k}$ depending  on the 
rank $n$ and weight $k$. In \cite{Fro93} Frohman finds a recursive 
procedure to compute the invariants $\lambda_{n,k}$ and shows that they 
are determined by the polynomial expressions in the
coefficients of the Alexander polynomial. Consequently, they are also
polynomial in the  Alexander Characters so that we can write
\beq\label{eq-laRnk}
\lambda_{n,k}\;=\; R_{n,k}(\Delta^{(1)}, \Delta^{(2)}, \ldots)\;,
\eeq
with $R_{n,k}\in{\mathbb Z}[x_1,x_2,\ldots]$. A general, closed
formula and some integrality issues for the $R_{n,k}$ are still unresolved 
though, see also \cite{BodNic00}. This relation in (\ref{eq-laRnk}) is more 
general than those expressed in  Theorem~\ref{thm-relations} as it is no longer
linear.

More precisely, define the space of invariants 
$\Qq^{[0]}\;=\;\left\{{n_1\Delta^{(1)}\,+\,n_2\Delta^{(2)}\,+\,\ldots\,|\,
n_i\in{\mathbb Z}^{+,0}}\right\}$. Then is clear that invariant that descends
from a TQFT that is homological must be in $\Qq^{[0]}$ where the $n_i\geq 0$
are the multiplicities of the irreducible summands. Thus $I^{DC},\,I^{SW}_d\,
\in\,\Qq^{[0]}$, but we also have ${\lambda}_L\not\in\Qq^{[0]}$ since some of
the coefficients are negative. ${\lambda}_L$ is, nevertheless, related to the
quantum TQFT's, but the derivations use slightly more subtle $p$-modular
interpretations, see \cite{KerKyoto}. 

\bigskip

Similar, to sums we can derive the 
invariant given by the product 
of two Alexander Characters, say $\Delta^{(i)}\cdot \Delta^{(j)}$,
from the tensor product of the corresponding TQFT's, namely
${\cal V}^{(i)}\otimes {\cal V}^{(j)}$. Thus is the coefficients of all
the $R_{n,k}$ were non-negative integers we could easily produce a 
homological TQFT by taking corresponding direct sums and tensor products
of the ${\cal V}^{(i)}$ in order to reproduce $\lambda_{n,k}$. This 
invariant, indeed, descends from the TQFT ${\cal V}^{PSU(n)}_k$, however, 
the coefficients of the $R_{n,k}$. The point to observe here is that, for
example, ${\cal V}^{(i)}\otimes {\cal V}^{(j)}$ is generally not an 
irreducible TQFT and can be decomposed. 

Denote by ${\cal V}^{(\vec \lambda)}$ the irreducible TQFT's obtained as summands of
quotients of multiple tensor products of the ${\cal V}^{(i)}$. The superscript label,
$\vec\lambda\in\vec\Lambda$, may be roughly thought of as a semi-infinite branching
path for $\Sp(2)\subset \Sp(4)\subset \Sp(6)\subset\ldots\,$. The space of TQFT's
\beq
\Qq^{[+]}\;=\;
\Biggl\{{\,\bigoplus_{{\vec\lambda}\in{\vec\Lambda}} \,{\mathbb Q}^{n_{\vec\lambda}}
\otimes {\cal V}^{(\vec\lambda)}\,\,\Bigg |\,\,\,n_{\vec\lambda}\in {\mathbb N}\cup
\{0\}}\Biggr\}\; 
\eeq
thus has a natural ring structure with operations $\oplus$ and $\otimes$ and can be 
thought of as a type of  Grothendieck $K_0$-ring for a homological subquotient of $\Cob$. 
We denote the corresponding set of higher Alexander Characters abusively in the same
way, since it possesses the same ring structure under usual addition and multiplication.
Clearly, $\Qq^{[0]}\subset\Qq^{[+]}$. The following conjecture together with an understanding
of the ring structure of $\Qq^{[+]}$ should shed light
onto the general structure of the polynomials $R_{n,k}$. 
\begin{conjecture}
$$
{\cal  V}^{PSU(n)}_k\;\;\;\in\;\;\;\Qq^{[+]}
$$
\end{conjecture}

\emptystuff{


\bigskip

\bigskip

Using the coefficients of
the Conway polynomial $\nabla (z)=\sum_j c_j z^{2j}$ with 
$z=t^{\frac 12}-t^{-\frac 12}$ this yields polynomials
$\lambda_{n,k}=q_{n,k}(c_1,c_2,\ldots)$, which appear to
have  non-negative
integer coefficients, that is $q_{n,k}\in{\cal P}^+$ as defined in
Section~10. See also \cite{BodNic00} for more explicit formulae.
Changing the basis of the polynomial ring from $z^{2j}$ to
the $[j]_{-t}$ we are similarly able to express the higher 
rank invariants in terms of the Alexander Characters $\Delta^{(j)}$ defined
in (\ref{eq-leftrace3}). We can thus write
$$
\lambda_{n,k}\;=\; R_{n,k}(\Delta^{(1)}, \Delta^{(2)}, \ldots)
$$
for some polynomial $R_{n,k}$ with integral coefficients.
Note further that if $C$ is the cobordism on a Seifert surface of a knot 
and $P\in{\cal P}^+$ then 
$trace({\cal V}^{(P)}(C))=P(\Delta^{(1)}, \Delta^{(2)}, \ldots)$

The natural
question tied to these observations is whether the TQFT's constructed 
in  \cite{FroNic94} for general gauge groups are related or  equivalent
to a TQFT of the polynomial form ${\cal V}^{(R_{n,k})}$ as defined in 
(\ref{eq-defpolyTQFT}). If the coefficients of the $R_{n,k}$ are not
all non-negative we may have to consider two theories 
${\cal V}^{(R_{n,k}^{\pm})}$ with $R_{n,k}=R_{n,k}^+-R_{n,k}^-$ and
$R_{n,k}^{\pm}\in{\cal P}^+$ and make sense of their difference.

In \cite{Don99} Donaldson describes a slightly different TQFT, ${\cal V}^{DF}$, 
modeled on moduli spaces ${\cal M}_g$ of connections on a non-trivial $SO(3)$ bundle. 
This TQFT leads up to a Casson type invariant for homology circles $Y$, which is
determined by the coefficients of the   Alexander polynomial $\Delta_Y$. 
The vector spaces are given as 
\beq\label{eq-Donaldson}
{\cal V}^{DF}(\Sigma_g)\;\;=\;\;H_*\bigl({\cal M}_g\bigr)
\;\;\cong\;\;\bigoplus_{j=0}^g\,\R^{j^2}\otimes \ext {g-j}H_1(\Sigma_g)\;.
\eeq
The morphisms ${\cal V}^{DF}(M)$ are similarly constructed via the intersection
theory of representation varieties, using also  a dimension
reduction of the Floer-cohomology on $\hat M\times S^1$.  

Now from Corollary~5.1.9 in \cite{GooWal98} we see that
 $\ext {g-j}H_1(\Sigma_g)\,\cong\,W_{g,j+1}\oplus W_{g,j+3}\oplus W_{g,j+5}\ldots$
as $\Sp(2g,\Z)$ modules. Inserting this decomposition into (\ref{eq-Donaldson}) we
obtain the multiplicities stated in the following conjecture.
\begin{conjecture} Let $\displaystyle D=\sum_{k\geq 1}
{\scriptstyle { k+1 \choose  3}} x_k\;\in\;{\cal P}^+$. Then 
$$
{\cal V}^{DF}\;\cong\;{\cal V}^{(D)}\;.
$$
\end{conjecture}
Note that on the level of vector spaces and invariants we do in fact have equality.

The theories in \cite{FroNic94} and \cite{Don99} are all inherently
$\Z/2$-projective, and have the vanishing properties from Lemma~\ref{lm-vanish}.
This indicates that they also belong into the class  of
half-projective or non-semisimple TQFT's.

Another conjecture that is independent of a particular gauge theory may be stated
as follows.
\begin{conjecture} Suppose ${\cal V}$ is a non-semisimple TQFT in which the kernel
of the mapping class group representations are precisely the Torelli groups. The
 ${\cal V}$ is isomorphic to a sub-TQFT of some ${\cal V}^{(P)}$ for some 
$P\in{\cal P}^+$.
\end{conjecture}
To say that the Torelli group is {\em precisely} the kernel implies that we
have faithful $\Sp(2g,\Z)$-representations so that by Margulis' rigidity
these lift to algebraic  $\Sp(2g,\R)$-representations. Classifying homological 
TQFT's as described will thus involve exercises in branching rules as given
for example in Section~8.3.4 of  \cite{GooWal98}.

Another approach to describing TQFT's behind the higher rank Frohman Nicas
$PSU(n)$-theories or the Donaldson $SO(3)$-construction is to try to extract
a categorical Hopf algebra ${\cal A}_{\cal V}$ for the given TQFT by evaluation
of the cobordisms in Figure~\ref{fig-MC-tgl}. The problem with this approach,
however, is that in the non-abelian TQFT's we do not seem to have a nicely
defined tensor structure arising from gluing two one holed surfaces over a 
pair of pants. Specifically we need an isomorphism 
${\cal V} (\Sigma_{2,1})=
{\cal V} (\Sigma_{1,1})\otimes{\cal V} (\Sigma_{1,1})$ which is generally
not true because of gauge constraints over the pair of pants. 

It is also not clear whether higher rank theories such as those in \cite{FroNic94}
exhibit symmetries similar to the $\SL(2,\R)$-equivariance
 that yields a type of Lefschetz
decomposition. Particularly,  the non-abelian moduli spaces have
no canonical K\"ahler structure. They do, however, admit useful Poisson structures
\cite{FocRos97}.
}


\paragraph{B. Homology TQFT's from the  Reshetikhin-Turaev Theory:}\

Recall that the TQFT ${\cal V}^{(j)}$ is in fact a functor to 
the category of free  $\Z$-modules rather than just the category
of vector spaces over ${\mathbb Q}$. Now, for any prime $p\geq 3$, 
by taking all lattices modulo $p$ 
this in turn maps to the category of vector spaces over the finite
field ${\mathbb F}_p=\Z/p\Z$. The resulting TQFT  ${\cal V}^{(j)}_p$
over ${\mathbb F}_p$
is now no longer irreducible, but it has a {\em unique}
irreducible subquotient, which we denote by $\dov {\cal V}^{(j)}_p$,
see \cite{KerRes}.

Another way of generating TQFT's over ${\mathbb F}_p$ is to consider
the Reshetikhin Turaev Theory for quantum-$SO(3)$ at a primitive
$p$-th root of unity $\zeta_p$. As shown in \cite{Gil01} this can
be regarded as a TQFT over  the ring of cyclotomic integers
$\Z[\zeta_p]$.  The TQFT obtained from the ring reduction
$\Z[\zeta_p]\onto{12}{\mathbb F}_p\,:\,\zeta_p\mapsto 1$ is denoted
${\cal V}^{RT}_p$. The example $p=5$, which is in some sense a
fundamental case, is analyzed in \cite{Kerr=5}. We obtain an exact,
but non-split sequence of TQFT's as follows
\beq\label{eq-ext}
0\;\;\to\;\;\dov {\cal V}^{(4)}_5
\;\;\longrightarrow\;\;
{\cal V}^{RT}_5
\;\;\longrightarrow\;\;
\dov {\cal V}^{(1)}_5
\;\;\to\;\;0
\eeq
As an extension of the mapping class group $\Gamma_g$ (\ref{eq-ext})
involves a Johnson-Morita subquotient of $\Gamma_g$. The precise modular
structure of the ${\cal V}^{(j)}_p$ and $\dov {\cal V}^{(j)}_p$
TQFT's is unraveled in \cite{KerRes}. There we find resolutions
of the $\dov {\cal V}^{(j)}_p$ in terms of the ${\cal V}^{(j)}_p$,
which lead to important identities between the $p$-modular versions
of the invariants from Theorem~\ref{thm-relations} and the 
Reshetikhin Turaev Invariants.

It is easy to see that the irreducible factors of ${\cal V}^{RT}_p$ for
$p\geq 7$ can no longer be reductions of the ${\cal V}^{(j)}$. There is,
however, evidence that suggests that the irreducible factors are 
reductions of summands in the symmetric powers of the fundamental ones.
That is, TQFT's of the form
\beq
{\cal V}^{\vec\lambda}
\;\subseteq\;S^{\frac {p-3}2}{\cal V}^{FN}\,\;\in\,\;\Qq^{[+]}\,.
\eeq 
This is closely related to the conjecture that the Lescop invariant
for a closed 3-manifold $M$ with $\beta_1(M)\geq 1$ relates to the 
Reshetikhin Turaev Invariant as follows.
\beq
{\cal V}^{RT}_{\zeta_p}(M)\;\;=\;\;C_p\cdot 
\Bigl((\zeta_p-1)\lambda_L(M)\Bigr)^{\frac {p-3}2}\;\;+\;\;{\cal O}
\bigl((\zeta_p-1)^{\frac {p-1}2}\bigr)\;\;.
\eeq
This has been verified for $p=5$ in \cite{KerKyoto}.

\paragraph{C. Relation of Reshetikhin-Turaev and Hennings Theory:}\ 

Given a quasitriangular Hopf algebra, ${\cal A}$, we have described in 
Section~5 a procedure to construct a topological quantum field theory,
${\cal V}_{\cal A}^{H}$. In \cite{ResTur91} and \cite{Tur94} Reshetikhin and
Turaev give another procedure to construct a TQFT, ${\cal V}^{RT}_{\cal S}$, 
from a {\em semisimple} modular  category, ${\cal S}$. A more general
construction in \cite{KerLub00} allows us to construct a TQFT, 
${\cal V}_{\cal C}^{KL}$, also for modular categories, ${\cal C}$,
that are not semisimple, and we show in \cite{Ker96} that 
${\cal V}_{\cal A}^{H}= {\cal V}_{{\cal A}-mod}^{KL}$ and 
${\cal V}^{RT}_{\cal S}={\cal V}_{\cal S}^{KL}$ for semisimple
 ${\cal S}$. For a non-semisimple, quasitriangular algebra, ${\cal A}$, 
the semisimple category used in \cite{ResTur91} and \cite{Tur94}
is given as the semisimple trace-quotient 
${\cal S}({\cal A})=\overline{{\cal A}-mod}$ of the representation
category of $\cal A$. The relation between ${\cal V}_{\cal A}^{H}$ 
and ${\cal V}^{RT}_{{\cal S}({\cal A})}$ is generally unknown. 
We make the following conjecture in the case of quantum ${\mathfrak s}{\mathfrak l}_2$.

\begin{conjecture}
Let ${\cal A}=U_q({\mathfrak s}{\mathfrak l}_2)^{red}$, with $q$ an odd $p$-th 
root of unity, and  relations $E^p=F^p=0$ and $K^{2p}=1$ for the standard generators. Then there is monomorphic, natural transformation 
\beq\label{eq-HRT}
{\cal V}^{FN}\otimes 
{\cal V}^{RT}_{{\cal S}({\cal A})}\;\;\into\;\;{\cal V}_{\cal A}^H\;.
\eeq
\end{conjecture}

In the genus one case we have shown in  \cite{Ker94} and \cite{Ker96}
that the mapping class group representations and invariants of 
lens spaces of both theories in ({\ref{eq-HRT}) are in fact equal. 
The above inclusion of TQFT functors can also be phrased in the form
${\cal V}_{{\cal C}^{\#}}^{KL}\into{\cal V}^{KL}_{{\cal C}}$, 
where ${\cal C}:=U_q({\mathfrak s}{\mathfrak l}_2)^{red}-mod$ and 
${\cal C}^{\#}:=({\cal N}-mod)\otimes \overline{\cal C}$. 
The categories ${\cal C}$ and ${\cal C}^{\#}$ are in fact rather
similar as linear abelian categories. From \cite{Ker98} it follows 
that there an isomorphism of {\em abelian} categories
\beq
\Hh\;\;:\;\;
{\cal C}^{\#}\;\oplus\;2\cdot {\rm Vect}({\mathbb C}) \;\;\;
\stackrel{\cong}{\longrightarrow}\;\;\;{\cal C}\;,
\eeq
where the two extra ${\rm Vect}({\mathbb C})$'s account for the two $p$-dimensional,
irreducible Steinberg modules. This, however, is {\em not} a 
monoidal functor. Instead we have a natural set of monomorphisms of
the form $\Hh(X)\otimes\Hh(Y)\,\into\,\Hh(X\otimes Y)$. As a result
the braidings, integrals, and coends, that enter in a crucial way
the construction of the TQFT's  \cite{KerLub00}
can no longer be na\"\i vely
identified. Strategies of proof would include a basis 
of ${\cal A}$ as worked out in \cite{Ker94} and the use of the special
central, nilpotent element ${\sf Q}$ defined in \cite{Ker96}. 


\emptystuff{

Specifically, we know the following: 
\begin{theorem}[\cite{Ker89}]\ Let ${\cal A}=U_q({\mathfrak s}{\mathfrak l}_2)^{red}$ and 
$\cal N$ as in Section~6. 
\begin{enumerate}
\item 
For any generic Casimir value, $c\in ({\mathfrak z}({\cal A}))^*$, 
the corresponding subcategory ${\cal C}_c\subset{\cal A}-mod$ 
of representations is isomorphic to ${\cal N}-mod$. 
\item The representations with non-generic Casimir values are 
sums of the two irreducible Steinberg modules of dimension $r$
and quantum dimension 0. 
\item An indecomposable representation of ${\cal N}$ is either 
 one of the two 4-dim
projective representations in  ${\cal N}={\cal N}^+\oplus{\cal N}^-$,
or an indecomposable representation of one of the two Kronecker quivers 
$\bullet \raisebox{-.5ex}{\hbox{$\stackrel{\mbox{$\longrightarrow$}}{\longrightarrow}$}}\bullet$ and 
$\bullet \raisebox{-.5ex}{\hbox{$\stackrel{\mbox{$\longleftarrow$}}{\longleftarrow}$}}\bullet$,
where the $\bullet$'s stand for an eigenspaces of $K$.  
\end{enumerate}
\end{theorem}
The generic Casimir values are in a two to one correspondence with the
admissible irreducible representations, and we have 
${\cal C}=\bigoplus_c{\cal C}_c$ and  
${\cal C}^{\#}=\bigoplus_j {\cal N}-mod$, where $j$ runs over  
irreducible representations. Thus we have a close correspondence between
the modules in both categories. They differ, however, more strongly 
as tensor categories. Strategies of proof would include a basis 
of ${\cal A}$ as worked out in \cite{Ker94} and the use of the special
central, nilpotent element ${\sf Q}$ defined in \cite{Ker96}. 

}


\emptystuff{

\paragraph{C. Universality of ${\cal V}$ and Casson Type Gauge Theories.}\ 

In order to find new knot invariants Frohman and Nicas generalized their
approach in  \cite{FroNic94} to higher rank Lie algebras. They construct a TQFT
whose vector spaces are
given as intersection homology groups of certain restricted moduli spaces 
of $PSU(n)$-representations and derive from these by similar trace formulae 
 invariants $\lambda_{n,k}$ depending  on the 
rank $n$ and weight $k$. In \cite{Fro93} Frohman finds a recursive 
procedure to compute the invariants $\lambda_{n,k}$ and shows that they 
are determined by the polynomial expressions in the
coefficients of the Alexander polynomial. Using the coefficients of
the Conway polynomial $\nabla (z)=\sum_j c_j z^{2j}$ with 
$z=t^{\frac 12}-t^{-\frac 12}$ this yields polynomials
$\lambda_{n,k}=q_{n,k}(c_1,c_2,\ldots)$, which appear to
have  non-negative
integer coefficients, that is $q_{n,k}\in{\cal P}^+$ as defined in
Section~10. See also \cite{BodNic00} for more explicit formulae.
Changing the basis of the polynomial ring from $z^{2j}$ to
the $[j]_{-t}$ we are similarly able to express the higher 
rank invariants in terms of the Alexander Characters $\Delta^{(j)}$ defined
in (\ref{eq-leftrace3}). We can thus write
$$
\lambda_{n,k}\;=\; R_{n,k}(\Delta^{(1)}, \Delta^{(2)}, \ldots)
$$
for some polynomial $R_{n,k}$ with integral coefficients.
Note further that if $C$ is the cobordism on a Seifert surface of a knot 
and $P\in{\cal P}^+$ then 
$trace({\cal V}^{(P)}(C))=P(\Delta^{(1)}, \Delta^{(2)}, \ldots)$

The natural
question tied to these observations is whether the TQFT's constructed 
in  \cite{FroNic94} for general gauge groups are related or  equivalent
to a TQFT of the polynomial form ${\cal V}^{(R_{n,k})}$ as defined in 
(\ref{eq-defpolyTQFT}). If the coefficients of the $R_{n,k}$ are not
all non-negative we may have to consider two theories 
${\cal V}^{(R_{n,k}^{\pm})}$ with $R_{n,k}=R_{n,k}^+-R_{n,k}^-$ and
$R_{n,k}^{\pm}\in{\cal P}^+$ and make sense of their difference.

In \cite{Don99} Donaldson describes a slightly different TQFT, ${\cal V}^{DF}$, 
modeled on moduli spaces ${\cal M}_g$ of connections on a non-trivial $SO(3)$ bundle. 
This TQFT leads up to a Casson type invariant for homology circles $Y$, which is
determined by the coefficients of the   Alexander polynomial $\Delta_Y$. 
The vector spaces are given as 
\beq\label{eq-Donaldson}
{\cal V}^{DF}(\Sigma_g)\;\;=\;\;H_*\bigl({\cal M}_g\bigr)
\;\;\cong\;\;\bigoplus_{j=0}^g\,\R^{j^2}\otimes \ext {g-j}H_1(\Sigma_g)\;.
\eeq
The morphisms ${\cal V}^{DF}(M)$ are similarly constructed via the intersection
theory of representation varieties, using also  a dimension
reduction of the Floer-cohomology on $\hat M\times S^1$.  

Now from Corollary~5.1.9 in \cite{GooWal98} we see that
 $\ext {g-j}H_1(\Sigma_g)\,\cong\,W_{g,j+1}\oplus W_{g,j+3}\oplus W_{g,j+5}\ldots$
as $\Sp(2g,\Z)$ modules. Inserting this decomposition into (\ref{eq-Donaldson}) we
obtain the multiplicities stated in the following conjecture.
\begin{conjecture} Let $\displaystyle D=\sum_{k\geq 1}
{\scriptstyle { k+1 \choose  3}} x_k\;\in\;{\cal P}^+$. Then 
$$
{\cal V}^{DF}\;\cong\;{\cal V}^{(D)}\;.
$$
\end{conjecture}
Note that on the level of vector spaces and invariants we do in fact have equality.

The theories in \cite{FroNic94} and \cite{Don99} are all inherently
$\Z/2$-projective, and have the vanishing properties from Lemma~\ref{lm-vanish}.
This indicates that they also belong into the class  of
half-projective or non-semisimple TQFT's.

Another conjecture that is independent of a particular gauge theory may be stated
as follows.
\begin{conjecture} Suppose ${\cal V}$ is a non-semisimple TQFT in which the kernel
of the mapping class group representations are precisely the Torelli groups. The
 ${\cal V}$ is isomorphic to a sub-TQFT of some ${\cal V}^{(P)}$ for some 
$P\in{\cal P}^+$.
\end{conjecture}
To say that the Torelli group is {\em precisely} the kernel implies that we
have faithful $\Sp(2g,\Z)$-representations so that by Margulis' rigidity
these lift to algebraic  $\Sp(2g,\R)$-representations. Classifying homological 
TQFT's as described will thus involve exercises in branching rules as given
for example in Section~8.3.4 of  \cite{GooWal98}.

Another approach to describing TQFT's behind the higher rank Frohman Nicas
$PSU(n)$-theories or the Donaldson $SO(3)$-construction is to try to extract
a categorical Hopf algebra ${\cal A}_{\cal V}$ for the given TQFT by evaluation
of the cobordisms in Figure~\ref{fig-MC-tgl}. The problem with this approach,
however, is that in the non-abelian TQFT's we do not seem to have a nicely
defined tensor structure arising from gluing two one holed surfaces over a 
pair of pants. Specifically we need an isomorphism 
${\cal V} (\Sigma_{2,1})=
{\cal V} (\Sigma_{1,1})\otimes{\cal V} (\Sigma_{1,1})$ which is generally
not true because of gauge constraints over the pair of pants. 

It is also not clear whether higher rank theories such as those in \cite{FroNic94}
exhibit symmetries similar to the $\SL(2,\R)$-equivariance
 that yields a type of Lefschetz
decomposition. Particularly,  the non-abelian moduli spaces have
no canonical K\"ahler structure. They do, however, admit useful Poisson structures
\cite{FocRos97}.

\paragraph{D. Milnor Torsion and Seiberg Witten invariants}\

The Milnor Torsion of a 3-manifold $M$ is defined from a cell complex of a 
simplicial representation of the cyclic covering space $\widetilde M$. 
The relation between the Alexander polynomial and Milnor Torsion as stated
in Theorem~\ref{thm-ReidAlex} suggests that there should be a quantum topological
description of the invariants obtained from a simplicial complex. In fact 
the Turaev Viro and Kuperberg invariant as examples of quantum invariants 
that start from a cell decomposition of $M$. Given the non-semisimple nature
of our theory the Kuperberg invariant \cite{Kup96} is a more natural candidate.  
The basic Hopf algebra is likely to be similar to the Borel subalgebra 
generated by $K$ and $\th$ but not $\tb$ following a conjecture that the
Hennings and Kuperberg are related by Drinfel'd double construction. The
difficulties, however, consist in describing a cell decomposition of the
cyclic covering space $\widetilde M$ from a Heegaard diagram for $M$. 
The main problem being that the Kuperberg theory has no easy extension 
to a TQFT. 

Nevertheless, we propose as a problem to find a direct description of
Milnor-Reidemeister torsion via the picture developed by Kuperberg for
the construction of 3-manifold invariants. 

In \cite{Tur98} Turaev shows that an extensions of
the  Milnor torsion is equal to the Seiberg Witten invariant 
for 3-manifolds equipped with $Spin^C$-structures or, equivalently, Euler
structures and with $\beta_1(M)>0$. A weaker version of such an equivalence 
without additional structures was shown by Meng and Taubes 
\cite{MengTaub}. 
 The proof in \cite{Tur98} uses the fact that
both invariants follows the same recursion formulae under surgery.
It should be interesting to relate these formulae to the skein theory
developed here and find ways of including the additional structures
in our context.

In general our procedure is also limited to either the reduced Torsion
or Alexander polynomial if $\beta_1(M)\geq 2$. Additional generators 
of homology can be represented as additional  0-framed surgery links with
zero linking numbers. It is, however, not as obvious in this case how 
to generate a TQFT picture that would allow us to describe the full
invariants with values in $\Z[H_1^{(free)}(M)]$.

\paragraph{E. Relations to quantum field theories:}\ 

Let us mention here only briefly interpretations of the homological TQFT's in
a physics context. For small genera Rozansky and Saleur
\cite{RosSal93} find the same vector spaces and representations of the
mapping class groups from the $U(1,1)$ Wess-Zumino-Witten theory. 

The exterior product spaces may also be interpreted as fermionic Fock spaces.
Ideas of constructing such ferminoc topological $U(1)$-theories in general
have been suggested for example by Louis Crane.  
}


{\sc\footnotesize  The Ohio State University, 

Department of Mathematics,

231 West 18th Avenue,

         Columbus, OH 43210, U.S.A. }

 {\em E-mail: }{ \tt kerler@math.ohio-state.edu}
\end{document}